\documentclass[sn-chicago]{sn-jnl}

\usepackage{graphicx}%
\usepackage{multirow}%
\usepackage{amsmath,amssymb,amsfonts}%
\usepackage{amsthm}%
\usepackage{mathtools}
\usepackage{mathrsfs}%
\usepackage{mathalpha}
\usepackage[title]{appendix}%
\usepackage[table]{xcolor}%
\usepackage{textcomp}%
\usepackage{manyfoot}%
\usepackage{booktabs}%
\usepackage{algorithm}%
\usepackage{algorithmicx}%
\usepackage{algpseudocode}%
\usepackage{listings}%
\usepackage{relsize} %
\usepackage{array}
\usepackage{bbm}

\usepackage{multirow,makecell}
\usepackage[font=small,skip=0pt]{caption}
\usepackage{subcaption}
\hypersetup{
	colorlinks=true,
	linkcolor=blue,
	filecolor=magenta,      
	urlcolor=cyan,
}

\theoremstyle{thmstyleone}%
\theoremstyle{thmstyletwo}%

\theoremstyle{thmstylethree}%

\newtheorem{observation}{Observation}

\newtheoremstyle{rem}
{6pt}
{6pt}
{\small}
{}
{\bf}
{:}
{.5em}
{}
\theoremstyle{rem}
\newtheorem{remark}{Remark}[section]

\begin{document}

\title[Article Title]{An objective isogeometric mixed finite element formulation for nonlinear elastodynamic beams with incompatible warping strains}


\author*[1]{\fnm{M.-J.} \sur{Choi}}\email{choi@lbb.rwth-aachen.de}

\author[1]{\fnm{S.} \sur{Klinkel}} 

\author[1]{\fnm{S.} \sur{Klarmann}} 

\author[2,3,4]{\fnm{R. A.} \sur{Sauer}} 


\affil[1]{\orgdiv{Chair of Structural Analysis and Dynamics}, \orgname{RWTH Aachen University}, \orgaddress{\street{Mies-van-der-Rohe Str.\,1}, \city{Aachen}, \postcode{52074}, \country{Germany}}} 
\affil[2]{\orgdiv{Institute for Structural Mechanics}, \orgname{Ruhr University Bochum}, \orgaddress{\street{Universit{\"a}tsstra{\ss}e 150}, \city{Bochum}, \postcode{44801}, \country{Germany}}}
\affil[3]{\orgdiv{Faculty of Civil and Environmental Engineering}, \orgname{Gda\'{n}sk University of Technology}, \orgaddress{\street{ul. Narutowicza 11/12}, \city{Gda\'{n}sk}, \postcode{80-233}, \country{Poland}}}
\affil[4]{\orgdiv{Department of Mechanical Engineering}, \orgname{Indian Institute of Technology Guwahati}, \orgaddress{\street{Assam 781039}, \country{India}}}


\abstract{
	We present a stable mixed isogeometric finite element formulation for geometrically and materially nonlinear beams in transient elastodynamics, where a Cosserat beam formulation with extensible directors is used. The extensible directors yield a linear configuration space incorporating constant in-plane cross-sectional strains. Higher-order (incompatible) strains are introduced to correct stiffness, whose additional degrees-of-freedom are eliminated by an element-wise condensation. Further, the present discretization of the initial director field leads to the objectivity of approximated strain measures, regardless of the degree of basis functions. \textcolor{black}{For physical stress resultants and strains, we employ a global patch-wise approximation using B-spline basis functions, whose higher-order continuity enables to use much less degrees-of-freedom, compared to element-wise approximation.} For time-stepping, we employ an implicit energy--momentum consistent scheme, which exhibits superior numerical stability in comparison to standard trapezoidal and mid-point rules. Several numerical examples are presented to verify the present method.
}




\maketitle

\section{Introduction}
\label{intro}
A Cosserat beam represents a slender object described by a spatial curve with two attached directors at each material point. Those directors span the (planar) cross-section, and are often assumed to be orthonormal, so that they can be parameterized by three rotational parameters. It has been shown that a standard finite element approximation of three-dimensional rotation parameters can lead to a loss of \textit{frame-invariance} of the underlying continuous strains, due to the non-additive nature of finite rotations, see \citet{crisfield1999objectivity} and \citet{jelenic1999geometrically} for the relevant investigation and its remedies using the concept of local rotation. \textcolor{black}{For a comprehensive review of different types of rotational parameterizations for nonlinear beams, one may refer to \citet{romero2004interpolation}.} It has also been shown that a direct finite element approximation of the director vectors can circumvent this issue, which, however, requires an additional constraint for the directors' orthonormality. For example, this constraint can be imposed at the nodal points by using a Lagrange multiplier method \citep{betsch2002frame}, or introducing nodal rotational parameters \citep{gruttmann2000theory,romero2002objective}. On the other hand, this orthonormality constraint can be simply abandoned, so that one has two independent extensible (unconstrained) directors, which represent in-plane shear and transverse normal strains of the cross-section. This increases the number of degrees-of-freedom (DOFs) per cross-section from six to nine, but the configuration space becomes linear, $\Bbb{R}^9$. The relevant approaches can be found, e.g., in \citet{rhim1998vectorial}, \citet{coda2009solid}, \citet{durville2012contact}, and \citet{choi2021isogeometric}. \textcolor{black}{Note that these approaches based on a direct finite element approximation of the director field may have a singularity in coarse meshes, when the directors' orientation changes abruptly along the length \citep[Section 4.3]{romero2004interpolation}.} The direct finite element approximation of (total) directors implies that the initial director field is also approximated. A continuous (exact) representation of the initial director field, inconsistent from the finite element approximation, may lead to a loss of objectivity. In \citet[Sec. 6.2]{choi2023selectively} it turned out that using a reduced degree of basis $p_\mathrm{d}$ for the change of directors lower than that of the initial (continuous) director field may result in the inability to represent rigid body rotations. In the present paper, we present a finite element approximation of the initial director field such that the objectivity is satisfied for any degree $p_\mathrm{d}$. For the relevant discussion on the objectivity under Kirchhoff beam kinematics, we refer interested readers to \citet{meier2014objective}. 

(Out-of-plane) warping means a cross-sectional deformation such that an initially planar cross-section does not remain plane. It can be coupled with other deformation modes like bending, shear, or torsion. In order to take this into account in Timohsenko beam kinematics, \citet{simo1991geometrically} introduced an additional DOF representing the amplitude of warping, and it is multiplied by a pre-defined warping function to produce an out-of-plane displacement, which has also been applied in the framework of unconstrained directors in \citet{coda2009solid}. This has been extended by \citet{klinkel2003anisotropic} to deal with anisotropic linear elastic materials which exhibits a bending--torsion coupling. In these previous works, the additional DOF for warping leads to non-standard stress resultants, bi--moment and bi--shear in the equilibrium equations. See also \citet{nukala2004mixed} for a relevant approach with simplified shear strains such that the transverse shear strain and the warping shear strain due to non-uniform torsion are neglected. In the framework of an enhanced assumed strain (EAS) method \citep{simo1990class}, one can enhance the cross-sectional strains that do not exist in the compatible strain field from the configuration. This requires additional constraints like orthogonality between the enhanced strain and the stress field. \citet{wackerfuss2009mixed} presented a construction of global polynomial basis functions of an arbitrary order for the enhanced cross-sectional strains. This has been extended to incorporate arbitrary (open) cross-section shapes by using a local polynomial basis in \citet{wackerfuss2011nonlinear}, \textcolor{black}{which has also been employed in \citet{klarmann2022coupling} for a coupling between beam and brick elements.} In the present work, the construction of global polynomial basis functions need to be more generalized, in order to account for the additional deformation modes in the compatible strain field such as transverse normal and in-plane shear strains, due to the extensible (unconstrained) directors. An alternative way to incorporate warping, which, however, requires increased number of global DOFs, is to use higher-order directors \citep{moustacas2019enrichissement,choi2022isogeometric}. One may also consider using multiple brick (solid beam) elements in the cross-section, e.g., as in \citet{frischkorn2013solid}.

The total linear momentum, total angular momentum, and total energy are fundamental integral quantities of motion in nonlinear elastodynamics, and their conservation (under suitable boundary conditions) manifests the stability of a numerical time-stepping scheme. For nonlinear problems, the standard trapezoidal rule, corresponding to $\beta=1/4$ and $\gamma={1}/{2}$ in the Newmark family, preserves the total linear momentum but neither the total angular momentum nor the total energy \citep{simo1995non}. The standard mid-point rule preserves both of the total linear and angular momentums but not the total energy \citep{gonzalez2000exact}. It has been shown that time-stepping schemes which preserve the conservation laws underlying the time-continuous form show superior numerical stability, compared with standard ones like the trapezoidal and mid-point rules. A class of these energy--momentum consistent (EMC) algorithms can be obtained by modifying the (implicit) mid-point rule. For example, one may refer to the algorithmic treatment of the stress resultant in \citet{romero2002objective}, and \citet{betsch2016energy} for geometrically nonlinear beam and shell problems, respectively, and \citet{gonzalez2000exact} for general (nonlinear) hyperelastic materials. 

The remainder of this paper is organized as follows. In Section \ref{bkin}, we introduce the beam kinematics based on the extensible directors. Section \ref{sect_var_form_mix_dyn} presents a mixed variational formulation for nonlinear dynamics of beams, which includes an EAS method for enriching higher-order cross-sectional strains. Section \ref{sec_emc_tstep} presents an implicit time-stepping scheme for energy--momentum consistency. Section \ref{spat_disc_iga} presents an isogeometric spatial discretization. In Section \ref{num_ex}, several numerical examples are presented. Section \ref{conclusions_sec} concludes the paper. 
\section{Beam kinematics} 
\label{bkin}
\begin{figure}[h]
	\centering
	\includegraphics[width=0.6\linewidth]{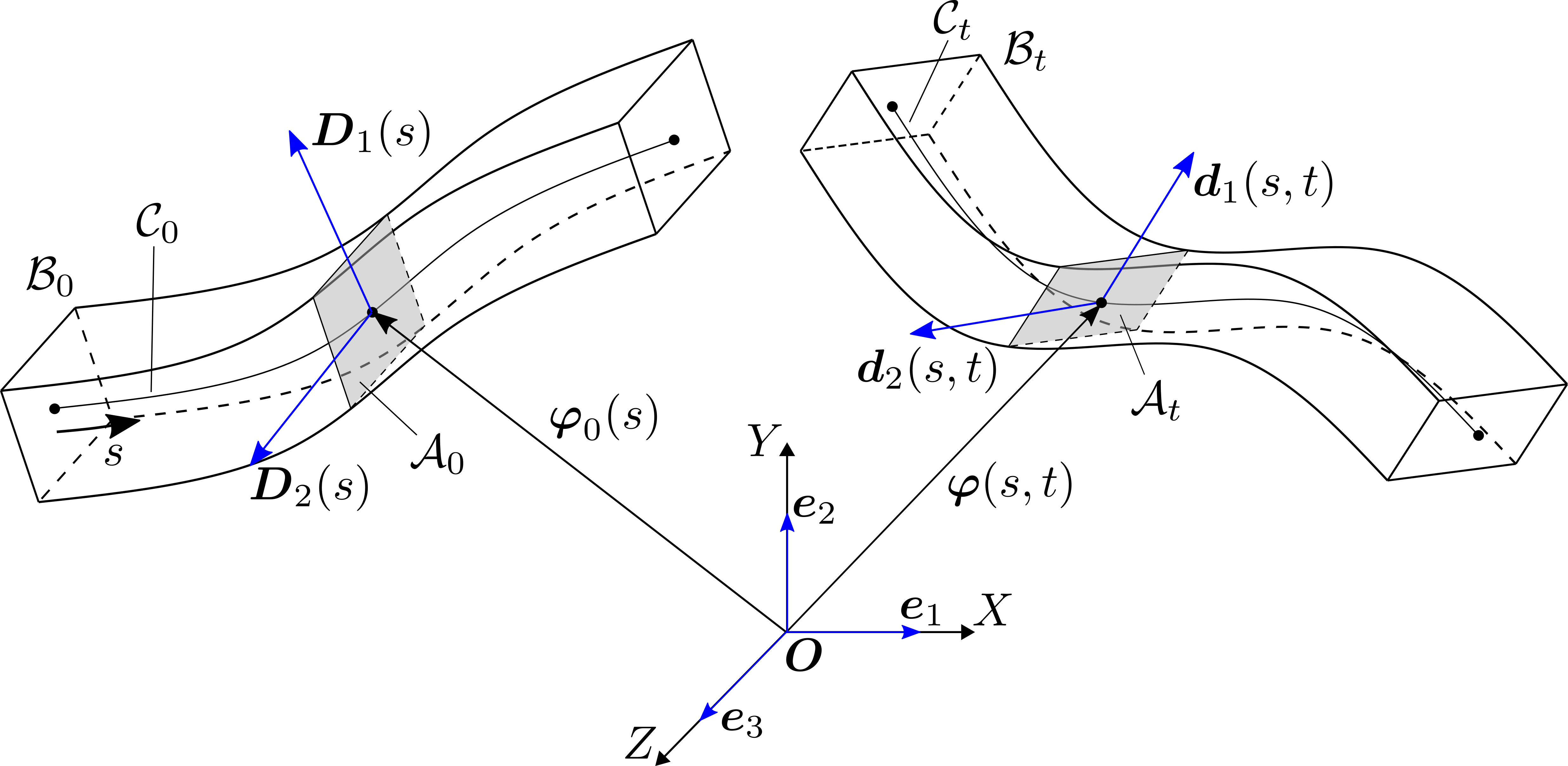}
	\caption{A schematic illustration of the beam kinematics. $\mathcal{B}_0$ and $\mathcal{B}_t$ denote the (open) domains of the initial (undeformed) and current configurations, respectively. ${\boldsymbol{e}}_1$, ${\boldsymbol{e}}_2$, and ${\boldsymbol{e}}_3$ represent the global Cartesian base vectors. This figure is redrawn with modifications from \citet{choi2023selectively}.}
	\label{redraw_beam_kin_glo}	
\end{figure}

\begin{figure}[h]
	\centering
	\includegraphics[width=0.4\linewidth]{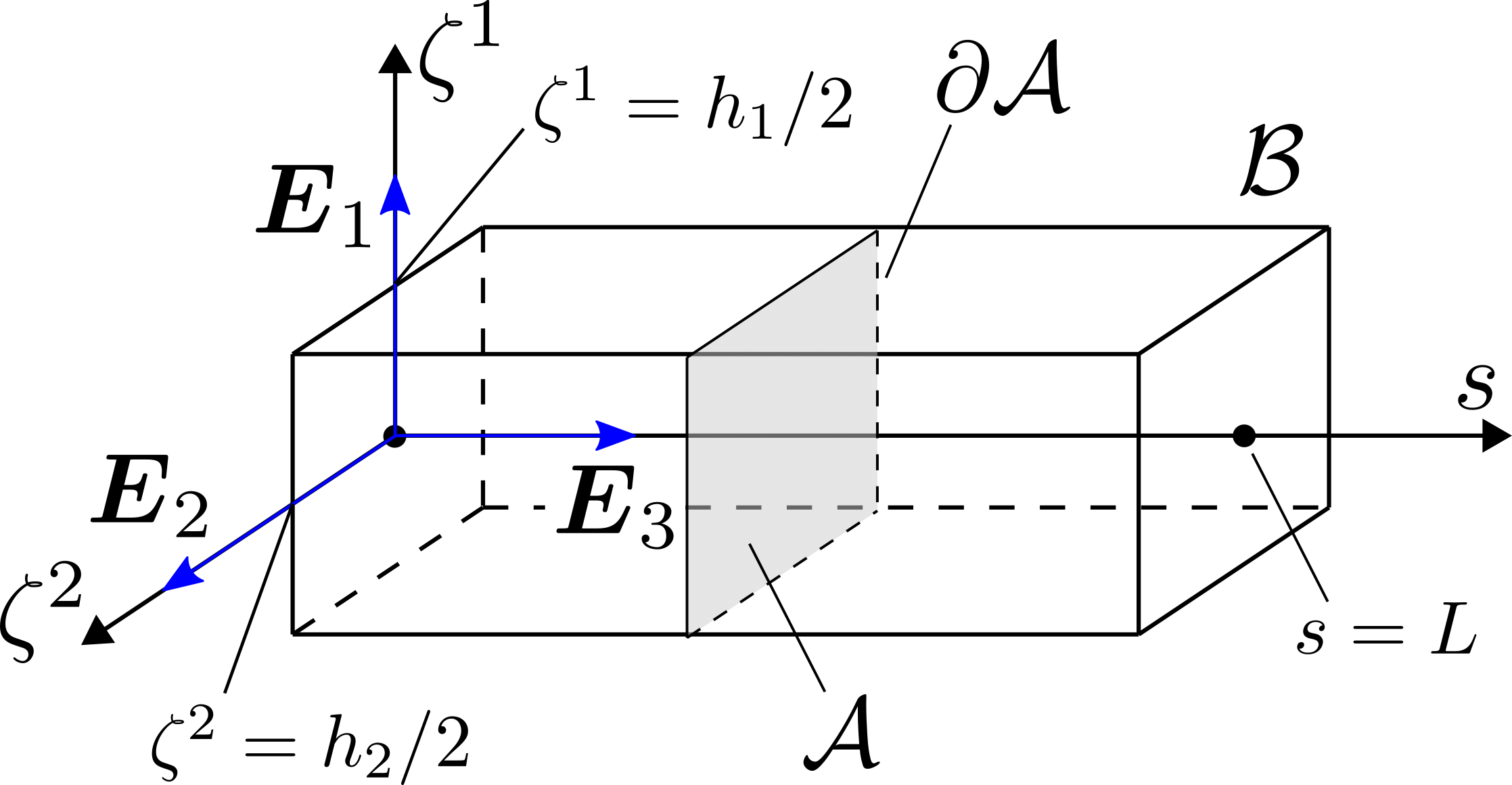}
	\caption{A schematic illustration of the reference domain of a beam having a rectangular cross-section with dimension $h_1\times{h_2}$. $\mathcal{B}$ denotes the (open) domain of the reference configuration. This figure is redrawn with modifications from \citet{choi2023selectively}.}
	\label{redraw_beam_kin_ref}	
\end{figure}

\textcolor{black}{This section presents the kinematical description of Cosserat beams, based on extensible directors, and the resulting configuration--strain relation.} A beam characterizes a three-dimensional body whose dimension along the longitudinal direction is much larger than the others such that a suitable kinematical assumption on the cross-sectional deformation can be employed. The origin of the transverse coordinates $\zeta^\alpha\,(\alpha\in\left\{1,2\right\})$ is given by the geometrical center of the initial cross-section, which coincides with the mass center, assuming constant mass density in the initial configuration. The line connecting these center points is called an \textit{initial center axis} whose position is represented by the vector ${\boldsymbol{\varphi }}_0({s})$, $s\in\left[0,L\right]$, where ${s}$ denotes its arc-length coordinate, and $L$ denotes its length. Then, the initial beam configuration can be expressed by
\begin{equation}
	\label{bkin_init}
	{{\boldsymbol{x}}}_0({\zeta ^1},{\zeta ^2},s) = {\boldsymbol{\varphi }}_0(s) + {\zeta ^\alpha}{{\boldsymbol{D}}_\alpha}(s),
\end{equation}
where the two orthonormal vectors ${\boldsymbol{D}}_\alpha(s) \in \mathbb{R}^3$ $(\alpha\in\left\{1,2\right\})$ are called \textit{initial directors}, and they span the initial cross-sectional plane, $\mathcal{A}_0$ in Fig.\,\ref{redraw_beam_kin_glo}. Those initial directors are defined along two principal directions of the cross-section's second moment of inertia. Here and hereafter, unless stated otherwise, repeated Greek indices like $\alpha$, $\beta$, and $\gamma$ imply summation over $1$ to $2$, and repeated Latin indices like $i$ and $j$ imply summation over $1$ to $3$. The initial directors are unit vectors and dimensionless, such that the cross-sectional dimension is associated with the range of the \textit{transverse coordinates} $\left(\zeta^1,\zeta^2\right)\in\mathcal{A}\,\cup\,{\partial {\mathcal{A}}}$, where $\mathcal{A}$ denotes the (open) domain of the (initial) cross-section and $\partial {\mathcal{A}}$ denotes its boundary, see Fig.\,\ref{redraw_beam_kin_ref}. It should be noted that the transverse coordinates $\zeta^\alpha$ are independent of time. In the initial configuration, we define the covariant basis ${\boldsymbol{G}}_i\coloneqq{\partial{\boldsymbol{x}}_0}/{\partial{\zeta^i}}$ $(i\in\left\{1,2,3\right\})$ with $\zeta^3\equiv{s}$, and the contravariant basis $\left\{{\boldsymbol{G}}^1,{\boldsymbol{G}}^2,{\boldsymbol{G}}^3\right\}$ can then be uniquely determined by the orthogonality condition $\boldsymbol{G}_i\cdot\boldsymbol{G}^j=\delta_i^j$, where $\delta_i^j$ denotes the Kronecker-delta.
\begin{remark} 
	\label{rem_init_curv_dir}
	\textit{Initial curvature}. The initial covariant base vectors are obtained by \citep{choi2021isogeometric}
	\begin{equation}\label{beam_th_str_init_cov_base}
		\left\{ \begin{array}{l}
			\begin{aligned}
				{{\boldsymbol{G}}_1} &= {{\boldsymbol{D}}_1}(s),\\
				{{\boldsymbol{G}}_2} &= {{\boldsymbol{D}}_2}(s),\\
				{{\boldsymbol{G}}_3} &= {{\boldsymbol{\varphi}}_{0,s}}(s) + {{\zeta ^\alpha }{{\boldsymbol{D}}_{\alpha,s}}(s)}.\\
			\end{aligned}
		\end{array} \right.
	\end{equation}
	The infinitesimal volume in the initial configuration can be decomposed into 
	\begin{align}
		\label{init_inf_vol_0}
		{\rm{d}}{\mathcal{B}_0} = {j_0}\,{\rm{d}}\mathcal{A}\,{\rm{d}}s,
	\end{align}	 	
	where the Jacobian is expressed by 
	\begin{align}
		{j_0} \coloneqq \left( {{{\boldsymbol{G}}_1} \times {{\boldsymbol{G}}_2}} \right) \cdot {{\boldsymbol{G}}_3} =1 - {\zeta ^\alpha }{\kappa_\alpha },
	\end{align}
	with
	\begin{align}
		{\kappa_\alpha} \coloneqq {{\boldsymbol{\varphi }}_{0,ss}} \cdot {{\boldsymbol{D}}_\alpha },\,\,\alpha\in\left\{1,2\right\},
	\end{align}
	which is associated with the initial curvature of the center axis, for a given initial center axis curve $\boldsymbol{\varphi}_0\equiv\boldsymbol{\varphi}_0(s)$ and the attached initial directors $\boldsymbol{D}_\alpha\equiv\boldsymbol{D}_\alpha(s)$.  Note that, for initially straight beams or linear elements, $\kappa_\alpha=0$, so that $j_0=1$.
\end{remark}
We introduce the first order beam kinematics, which expresses the current position of a material point, as
\begin{equation}
	\label{bkin_x_d_zta_ph}
	{{\boldsymbol{x}}}(\zeta^1,\zeta^2,s,t) = {\boldsymbol{\varphi}}(s,t) + {\zeta ^\alpha }{{\boldsymbol{d}}_\alpha }(s,t),
\end{equation}
where ${\boldsymbol{\varphi }}(s,t)$ denotes the current position of the center axis, and the current cross-sectional plane is spanned by two \textit{current directors} ${\boldsymbol{d}}_\alpha(s,t) \in \mathbb{R}^3$ $(\alpha\in\left\{1,2\right\})$, see Fig.\,\ref{redraw_beam_kin_glo}. Here, $t\in\left[{0},{{T}}\right]$ denotes time, where ${T}$ represents the terminal time. We obtain the velocity by differentiating Eq.\,(\ref{bkin_x_d_zta_ph}) with respect to time, as
\begin{equation}
	\label{bkin_velocity}
	{\boldsymbol{\dot x}}(\zeta^1,\zeta^2,s,t) = {\boldsymbol{\dot \varphi}}(s,t) + {\zeta ^\alpha }{{\boldsymbol{\dot d}}_{\alpha}}(s,t),
\end{equation}
where \textcolor{black}{the upper dot denotes the material time derivative.} Hereafter, for brevity, we often omit the arguments, e.g., ${\boldsymbol{x}}\equiv{\boldsymbol{x}}(\zeta^1,\zeta^2,s,t)$, $\boldsymbol{\varphi}\equiv\boldsymbol{\varphi}(s,t)$ and $\boldsymbol{d}_\alpha\equiv{\boldsymbol{d}_\alpha}(s,t)$. Further, we define a configuration variable $\boldsymbol{y}\coloneqq {\left[ {{{\boldsymbol{\varphi }}^{\rm{T}}},{{\boldsymbol{d}}^\mathrm{T}_1},{{\boldsymbol{d}}^{\rm{T}}_2}} \right]^{\!\rm{T}}}\in\Bbb{R}^d$, where $d=9$ denotes the number of independent cross-sectional components.
\subsection{Geometric strain: strain-configuration relation}
We may represent a three-dimensional strain state of a beam by so called
\textit{beam strains} evaluated at the center axis. In the current configuration, we define the covariant basis ${\boldsymbol{g}}_i\coloneqq{\partial{\boldsymbol{x}}}/{\partial{\zeta^i}}$ $(i\in\left\{1,2,3\right\})$, and then the contravariant basis $\left\{{\boldsymbol{g}}^1,{\boldsymbol{g}}^2,{\boldsymbol{g}}^3\right\}$ can be uniquely determined by the orthogonality condition $\boldsymbol{g}_i\cdot\boldsymbol{g}^j=\delta_i^j$. Then, the deformation gradient tensor $\boldsymbol{F}$ can be expressed by \citep[page 31]{wriggers2008nonlinear}
\begin{equation}
	\boldsymbol{F}\coloneqq{\partial{\boldsymbol{x}}/{\partial{\boldsymbol{x}}_0}}=\boldsymbol{g}_i\otimes{\boldsymbol{G}^i}.
\end{equation}
The Green-Lagrange strain tensor ${\boldsymbol{E}}\equiv{\boldsymbol{E}}(\boldsymbol{y})$ is defined by
\begin{align}
	{\boldsymbol{E}}(\boldsymbol{y}) \coloneqq \dfrac{1}{2}\left( {{{\boldsymbol{F}}^{\rm{T}}}{\boldsymbol{F}} - {\boldsymbol{1}}} \right)={E_{ij}}\,{\boldsymbol{G}^i}\otimes{\boldsymbol{G}^j},
\end{align}
where $\boldsymbol{1}$ denotes the identity tensor. In Voigt notation, we have 
\begin{align}
	\label{E_vgt_decomp_comp}
	\underline{\boldsymbol{E}}  &\coloneqq {\left[ {{E_{11}},{E_{22}},{E_{33}},2{E_{12}},2{E_{13}},2{E_{23}}} \right]^{\rm{T}}}= {\boldsymbol{A}}({\zeta ^1},{\zeta ^2})\,{\boldsymbol{\varepsilon }}(\boldsymbol{y}),
\end{align}
where the underline $\underline {\left(\bullet\right)}$ denotes Voigt notation of a second order tensor, and we have defined the matrix ${\boldsymbol{A}}\equiv{{\boldsymbol{A}}(\zeta^1,\zeta^2)}$ \citep{choi2021isogeometric}
\begin{align}
	\label{mat_A_poly_basis}
	{{\boldsymbol{A}}(\zeta^1,\zeta^2)}\coloneqq\left[{\setlength{\arraycolsep}{2.25pt}\renewcommand{\arraystretch}{1.45}\begin{array}{*{20}{c}}
			0&0&0&0&0&0&0&0&0&0&0&0&1&0&0\\
			0&0&0&0&0&0&0&0&0&0&0&0&0&1&0\\
			1&{{\zeta^1}}&{{\zeta ^2}}&{{\zeta ^1}{\zeta ^1}}&{{\zeta ^2}{\zeta ^2}}&{{\zeta ^1}{\zeta ^2}}&0&0&0&0&0&0&0&0&0\\
			0&0&0&0&0&0&0&0&0&0&0&0&0&0&1\\
			0&0&0&0&0&0&1&0&{{\zeta ^1}}&{{\zeta ^2}}&0&0&0&0&0\\
			0&0&0&0&0&0&0&1&0&0&{{\zeta ^1}}&{{\zeta ^2}}&0&0&0
	\end{array}} \right],
\end{align}
and the beam strain components 
\begin{subequations}
	\label{beam_th_strn_comp_basis_d123}
	\begin{alignat}{2}
		\varepsilon  &\coloneqq \frac{1}{2}({\left\| {{{\boldsymbol{\varphi }}_{\!,s}}} \right\|^2} - 1),\label{axial_strn_comp}\\ 
		{{\rho }_\alpha } &\coloneqq {{\boldsymbol{\varphi}}_{\!,s}}\!\cdot\!{{\boldsymbol{d}}_{\alpha ,s}} - {{\boldsymbol{\varphi}}_{0,s}}\!\cdot\!{{\boldsymbol{D}}_{\alpha ,s}},\label{bend_strn_comp}\\ 
		{{\delta }_\alpha } &\coloneqq {{\boldsymbol{\varphi}}_{\!,s}}\!\cdot\!{{\boldsymbol{d}}_\alpha }- {{\boldsymbol{\varphi }}_{0,s}}\!\cdot\!{{\boldsymbol{D}}_\alpha },\label{tshear_strn_comp}\\ 
		{{\gamma }_{\alpha \beta }} &\coloneqq {{\boldsymbol{d}}_\alpha }\!\cdot\!{{\boldsymbol{d}}_{\beta ,s}} - {{\boldsymbol{D}}_\alpha }\!\cdot\!{{\boldsymbol{D}}_{\beta ,s}},\label{cshear_strn_comp}\\ 
		{{\chi }_{\alpha \beta }} &\coloneqq \frac{1}{2}({{\boldsymbol{d}}_\alpha }\!\cdot\!{{\boldsymbol{d}}_\beta } - \boldsymbol{D}_{\alpha}\!\cdot\!\boldsymbol{D}_{\beta}),\\ 
		{{{\kappa}} _{\alpha \beta }} &\coloneqq \frac{1}{2}\left( {{{\boldsymbol{d}}_{\alpha ,s}}\!\cdot\!{{\boldsymbol{d}}_{\beta ,s}} - {{\boldsymbol{D}}_{\alpha ,s}}\!\cdot\!{{\boldsymbol{D}}_{\beta ,s}}} \right), \label{strn_high_bend} 
	\end{alignat}
	$\alpha,\beta\in\left\{1,2\right\}$,
\end{subequations}
which are arrayed in the form
\begin{align}
	\label{geom_b_strn_form_arr}
	{\boldsymbol{\varepsilon }}(\boldsymbol{y}) &\coloneqq \left[ \varepsilon,{\rho_1},{\rho_2},{{{\kappa }}_{11}},{{{\kappa }}_{22}},2{{{\kappa }}_{12}},{{\delta }}_1,{{\delta }}_2,{{{\gamma }}_{11}},{{{\gamma }}_{12}},{{{\gamma }}_{21}},{{{\gamma }}_{22}},{\chi}_{11},{\chi}_{22},2{\chi}_{12}\right]^{\rm{T}}\in\Bbb{R}^{d_\mathrm{p}},
\end{align}
where $d_\mathrm{p}=15$ denotes the number of independent strain components. Those beam strain components are called
\begin{alignat*}{3}
	&\varepsilon&&: \text{axial strain},\nonumber\\
	&{\rho _1},{{\rho _2}}&&: \text{bending strain},\nonumber\\
	&{\kappa _{11}},{{\kappa _{22}}},{{\kappa _{12}}}&&: \text{higher-order bending strain},\nonumber\\
	&{{\delta _1}},{{\delta _2}}&&: \text{transeverse shear strain},\nonumber\\	
	&{{\gamma _{11}}},{{\gamma _{12}}},{{\gamma _{21}}},{{\gamma _{22}}}&&: \text{couple shear strain},\nonumber\\
	&{{\chi_{11}}},{{\chi_{22}}},{{\chi_{12}}}&&: \text{in-plane stretching and shear strains}.\nonumber
\end{alignat*}
Note that $\kappa_{12}=\kappa_{21}$ and $\chi_{12}=\chi_{21}$. We call $\boldsymbol{\varepsilon}(\boldsymbol{y})$, expressed in terms of the configuration variable $\boldsymbol{y}$, a \textit{geometric strain}, which represents the strain-configuration relation, whose objectivity is proved in Section\,\ref{verif_objectivity_strn_cont}.
\section{Mixed variational formulation}
\label{sect_var_form_mix_dyn}
This section presents a mixed variational formulation for nonlinear transient dynamics of beams, based on the kinematics in Section\,\ref{bkin}.
\subsection{Variational forms}
\subsubsection{Total strain energy: Hu-Washizu principle}
\label{vform_hw_princ}
Based on the Hu-Washizu variational principle, the total strain energy, ${U}_{{\mathop{\rm HW}} } \equiv {U}_{{\mathop{\rm HW}}}({\boldsymbol{y}},{\boldsymbol{E}}_\mathrm{p},{\boldsymbol{S}}_\mathrm{p})$, can be expressed, using Eq.\,(\ref{init_inf_vol_0}), as
\begin{equation}
	\label{hw_func_init}
	{U}_{{\mathop{\rm HW}} } \coloneqq\int_0^L{\int_\mathcal{A} {\bigl[{\Psi\!\left( {{{\boldsymbol{E}}_{\rm{p}}}} \right)+\left\{ {{\boldsymbol{E}}(\boldsymbol{y})-{{\boldsymbol{E}}_{\rm{p}}}} \right\}:{{\boldsymbol{S}}_{\rm{p}}}} \bigr]\,{j_0}\,{\rm{d}}\mathcal{A}\,{\rm{d}}s}},
\end{equation}
where $\Psi\!\left( {{{\boldsymbol{E}}_{\rm{p}}}} \right)$ denotes the strain energy density per unit undeformed volume, expressed in terms of the \textit{physical} (independent) Green-Lagrange strain tensor, ${\boldsymbol{E}}_\mathrm{p}$. The compatibility condition,
\begin{equation}
	\label{GL_strn_compat}
	{{\boldsymbol{E}}(\boldsymbol{y})={{\boldsymbol{E}}_{\rm{p}}}}
\end{equation}
is enforced by a Lagrange multiplier method, where the Lagrange multiplier ${\boldsymbol{S}}_\mathrm{p}$ can be interpreted as the \textit{physical} second Piola-Kirchhoff stress tensor. This type of variational formulation is motivated by the following:
\begin{itemize}
	\item \textcolor{black}{The introduction of an independent strain ${{\boldsymbol{E}}_{\rm{p}}}$ aims at an alleviation of numerical locking due to the parasitic terms in the compatible (geometric) strain $\boldsymbol{E}(\boldsymbol{y}^h)$ from the finite element approximation of the configuration variable, $\boldsymbol{y}\approx\boldsymbol{y}^h$.}
	\item The introduction of an independent stress field ${\boldsymbol{S}}_\mathrm{p}$ alleviates the overestimation of internal force at inequilibrium configurations, which turns out to allow much larger load increments, compared with the displacement-based formulation.
\end{itemize}
\textcolor{black}{The relevant investigations can be found in}, e.g., \citet{klinkel2006robust} for elastostatics, and \citet{betsch2016energy} for elastodynamics. Further, we may enhance the physical strain field, as \citep{simo1990class,buchter1994three}
\begin{equation}
	\label{dec_phy_str_k_eh}
	{{\boldsymbol{E}}^\mathrm{tot}} = \underbrace{{\boldsymbol{E}}_{\rm{p}}}_{\textrm{kinematic}} + \underbrace {{\boldsymbol{\widetilde E}}}_{{\rm{enhanced}}},
\end{equation}
where `tot' means `total', and the `kinematic' part means its compatibility with the geometric (compatible) strain, $\boldsymbol{E}(\boldsymbol{y})$ in Eq.\,(\ref{GL_strn_compat}). Substituting Eq.\,(\ref{dec_phy_str_k_eh}) into Eq.\,(\ref{hw_func_init}), the total strain energy in Eq.\,(\ref{hw_func_init}) can be rewritten as
\begin{equation}
	\label{hw_func_mod}
	{U}_{{\mathop{\rm HW}} } \!=\!\int_0^L\!{\int_\mathcal{A} {\left[\renewcommand{\arraystretch}{1.25}\begin{array}{l}\Psi\!\left( {{{\boldsymbol{E}}^\mathrm{tot}}} \right)
				+\left\{ {{{\boldsymbol{E}}}(\boldsymbol{y})-{{\boldsymbol{E}}_{\rm{p}}}} \right\}:{{\boldsymbol{S}}_{\rm{p}}}\end{array}\right]{j_0}\,{\rm{d}}\mathcal{A}\,{\rm{d}}s}},
\end{equation}
where the following orthogonality condition is utilized 
\begin{equation}
	\label{S_E_ortho_cond}
	{\int_\mathcal{A} {{{\boldsymbol{\widetilde E}}}:{{\boldsymbol{S}}_{\rm{p}}}\,{j_0}\,{\rm{d}}\mathcal{A}} = 0}.
\end{equation}
The compatibility in Eq.\,(\ref{GL_strn_compat}) enables to represent the three-dimensional strain state in terms of the \textit{physical strain} ${\boldsymbol{\varepsilon}}_\mathrm{p}\equiv{\boldsymbol{\varepsilon}}_\mathrm{p}(s,t)\in\Bbb{R}^{{d}_\mathrm{p}}$, as
\begin{align}
	\label{phy_kin_dec}
	{\underline{{\boldsymbol{E}}_\mathrm{p}}} = {\boldsymbol{A}}({\zeta ^1},{\zeta ^2})\,{\boldsymbol{\varepsilon}}_\mathrm{p}(s,t),
\end{align}
as in Eq.\,(\ref{E_vgt_decomp_comp}). This, so called an enhanced assumed strain (EAS) method is introduced to correct stiffness, which, especially, takes the following into account:
\begin{itemize}
	\item Correction of the bending stiffness, considering the fact that a bending deformation may lead to (at least) linear in-plane cross-sectional strain, due to non-zero Poisson's ratio,
	\item Correction of the torsional stiffness, considering the fact that torsion may induce cross-sectional warping,
	\item Correction of the transverse shear stiffness. 
\end{itemize}

\noindent The enhanced strain can also be represented by the parameters $\boldsymbol{\alpha}\in\Bbb{R}^{d_\mathrm{a}}$, as
\begin{equation}
	\label{phy_enh_dec}
	{\underline{\boldsymbol{\widetilde E}}} = {\boldsymbol{\Gamma }}({\zeta ^1},{\zeta ^2})\,{\boldsymbol{\alpha }}(s,t),
\end{equation}
with the following matrix of basis functions \citep{wackerfuss2009mixed}
\begin{align}
	\label{wg_basis_matrix_gamma}
	{\boldsymbol{\Gamma }}\!\!\:\left( {{\zeta ^1},{\zeta ^2}} \right) = \left[ {\begin{array}{*{20}{c}}
			{{{\bf{w}}_1}}&{}&{}&{}&{}\\
			{}&{{{\bf{w}}_1}}&{}&{}&{}\\
			{}&{}&{{{\bf{w}}_3}}&{}&{}\\
			{}&{}&{}&{{{\bf{w}}_2}}&{}\\
			{}&{}&{}&{}&{{\bf{w}}_{4,1}}\\
			{}&{}&{}&{}&{{{\bf{w}}_{4,2}}}
	\end{array}} \right]_{\textcolor{black}{6\times{d_\mathrm{a}}}},
\end{align}
where the blank entries represent zeros, and $(\bullet)_{,\alpha}$ denotes the derivative with respect to $\zeta^\alpha$ ($\alpha\in\left\{1,2\right\}$). \textcolor{black}{Here, ${\bf{w}}_1$, ${\bf{w}}_2$, ${\bf{w}}_3$, and ${\bf{w}}_4$ denote row arrays of the \textit{global} polynomial basis functions for enriching the transverse normal (${\widetilde E}_{11}$, ${\widetilde E}_{22}$), in-plane shear (${\widetilde E}_{12}$), axial normal (${\widetilde E}_{33}$), and transverse shear (${\widetilde E}_{13}$,${\widetilde E}_{23}$) strains, respectively. 
\vspace{2mm}
\begin{remark}	
For the transverse shear strains (${\widetilde E}_{13}$, ${\widetilde E}_{23}$), we employ the derivatives of ${\bf{w}}_{4}$. However, we distinguish ${\bf{w}}_{4}$ from ${\bf{w}}_{3}$, which is in contrast to the formulation of \citet{wackerfuss2009mixed}. This is due to the quadratic basis for $E_{33}$ in Eq.\,(\ref{mat_A_poly_basis}), which does not exist for $E_{13}$ and $E_{23}$. We verify the erroneous results from using ${\bf{w}}_4={\bf{w}}_3$ in Figs.\,\ref{sv_tor_cs_strn_E33_w0.3_h0.4_ortho_up_to_lin} and \ref{sv_torsion_w0.3_h0.4} in Section \ref{num_ex_twist_straight}.
\end{remark}
\vspace{2mm}
\noindent The detailed process of constructing ${\bf{w}}_i\in\Bbb{R}^{1\times{d_i}}$ ($i\in\left\{1,2,3,4\right\}$) is given in Appendix\,\ref{app_construct_warp_basis}. Then, the total number of independent components in the enhanced strain field is $d_\mathrm{a}=2{d_1}+{d_2}+{d_3}+{d_4}$, where the dimension $d_i$ is given by Eq.\,(\ref{dimension_w_i}).} \textcolor{black}{It is noted that} `global' implies that the basis functions span the whole domain of the cross-section, ${\mathcal{A}}\ni\!\left(\zeta^1,\zeta^2\right)$. This, however, is not suitable to represent localized warping, e.g., for a beam having an open cross-section. This limits the applicability of the present formulation to convex-shaped cross-sections like a rectangular one. In order to overcome this, one may consider the local polynomial basis functions in \citet{wackerfuss2011nonlinear}.  Substituting Eqs.\,(\ref{phy_kin_dec}) and (\ref{phy_enh_dec}) into Eq.\,(\ref{hw_func_mod}) gives
\begin{equation}
	\label{hw_int_ener}
	{U_{{\rm{HW}}}} = \int_0^L {[ {\psi \left( {{{\boldsymbol{\varepsilon }}_{\rm{p}}},{\boldsymbol{\alpha }}} \right) + \underbrace{\left\{{{\boldsymbol{\varepsilon }}(\boldsymbol{y}) - {{\boldsymbol{\varepsilon }}_{\rm{p}}}}\right\}}_{\textrm{compatibility}} \cdot\:{{\boldsymbol{r}}_{\rm{p}}}}]\,{\rm{d}}s},
\end{equation}
where we define the \textit{line energy density} (strain energy per undeformed unit length),
\begin{equation}
	\label{line_e_density_ref}
	\psi({\boldsymbol{\varepsilon}}_\mathrm{p},\boldsymbol{\alpha})\coloneqq\int_\mathcal{A} {{\Psi}\left({{\boldsymbol{E}}^\mathrm{tot}}(\boldsymbol{\varepsilon}_\mathrm{p},\boldsymbol{\alpha})\right){j_0}\,{\rm{d}}\mathcal{A}}.
\end{equation}
Here, the Lagrange multiplier
\begin{equation}
 {{\boldsymbol{r}}_{\rm{p}}}\coloneqq{\int_\mathcal{A} {{{\boldsymbol{A}}^{\rm{T}}}\,\underline{{\boldsymbol{S}}_{\rm{p}}}\,{j_0}\,{\rm{d}}\mathcal{A}}} \in\Bbb{R}^{d_\mathrm{p}}
\end{equation} 
for the compatibility condition, ${{\boldsymbol{\varepsilon }}(\boldsymbol{y}) = {{\boldsymbol{\varepsilon }}_{\rm{p}}}}$, in Eq.\,(\ref{hw_int_ener}) represents
an array of \textit{physical} stress resultants. Note that each component of the array ${{\boldsymbol{r}}_{\rm{p}}}$ is an independent variable, a work-conjugate to the corresponding component in the array ${\boldsymbol{\varepsilon}}_\mathrm{p}$. Taking the first variation of Eq.\,(\ref{hw_int_ener}), we obtain the internal virtual work \citep{choi2023selectively}
\begin{align}
	\label{vform_hw}
	G_{{\mathop{\mathrm {int}}} }^{{\mathrm{HW}}} \equiv \delta {{U}_{{\mathrm{HW}}}} &= G^{\rm{y}}_\mathrm{int}(\boldsymbol{y},{\boldsymbol{r}}_\mathrm{p},\delta{\boldsymbol{y}}) + G^{\rm{r}}_\mathrm{int}(\boldsymbol{y},{\boldsymbol{\varepsilon}}_\mathrm{p},\delta{\boldsymbol{r}}_\mathrm{p}) + G^\varepsilon_\mathrm{int}({\boldsymbol{r}}_\mathrm{p},{\boldsymbol{\varepsilon}}_\mathrm{p},{\boldsymbol{\alpha}},\delta{\boldsymbol{\varepsilon}}_\mathrm{p})\nonumber\\
	&+ G^{\rm{a}}_\mathrm{int}(\boldsymbol{\varepsilon}_\mathrm{p},{\boldsymbol{\alpha}},\delta{\boldsymbol{\alpha}}),
\end{align} 
where we have defined
\begin{subequations}
	\label{sub_vform_hw}
	\begin{align}
		G^{\rm{y}}_\mathrm{int}(\boldsymbol{y},{\boldsymbol{r}}_\mathrm{p},\delta{\boldsymbol{y}}) &\coloneqq \int_0^L {\delta {\boldsymbol{y}} \cdot {\mathbb{B}}({\boldsymbol{y}})^{\rm{T}}\,{{\boldsymbol{r}}_{\rm{p}}}\,{\rm{d}}s},\label{def_int_vw_y}\\
		G^{\rm{r}}_\mathrm{int}(\boldsymbol{y},{\boldsymbol{\varepsilon}}_\mathrm{p},\delta{\boldsymbol{r}}_\mathrm{p}) &\coloneqq \int_0^L {\delta {{\boldsymbol{r}}_{\rm{p}}} \cdot \left\{ {{\boldsymbol{\varepsilon }}({\boldsymbol{y}}) - {{\boldsymbol{\varepsilon }}_{\rm{p}}}} \right\}{\rm{d}}s},\label{def_int_vw_r}\\
		G^\varepsilon_\mathrm{int}({\boldsymbol{r}}_\mathrm{p},{\boldsymbol{\varepsilon}}_\mathrm{p},{\boldsymbol{\alpha}},\delta{\boldsymbol{\varepsilon}}_\mathrm{p})  &\coloneqq \int_0^L {\delta {{\boldsymbol{\varepsilon }}_{\rm{p}}}\!\cdot\!\left\{ {{\partial _{{\boldsymbol{\varepsilon}_{\rm{p}}}}}\psi ({{\boldsymbol{\varepsilon }}_{\rm{p}}},{\boldsymbol{\alpha }}) - {{\boldsymbol{r}}_{\rm{p}}}} \right\}{\rm{d}}s},\label{def_int_vw_e}\\
		G^{\rm{a}}_\mathrm{int}(\boldsymbol{\varepsilon}_\mathrm{p},{\boldsymbol{\alpha}},\delta{\boldsymbol{\alpha}}) &\coloneqq \int_0^L {\delta {\boldsymbol{\alpha }} \cdot {\partial _{\boldsymbol{\alpha}}}\psi ({{\boldsymbol{\varepsilon }}_{\rm{p}}},{\boldsymbol{\alpha }})\,{\rm{d}}s}.\label{def_int_vw_a}
	\end{align}
\end{subequations}
Here $\delta(\bullet)$ denotes the first variation, and we have defined
\begin{subequations}
	\begin{align}
		{\partial _{{{\boldsymbol{\varepsilon }}_{\rm{p}}}}}\psi \left( {{{\boldsymbol{\varepsilon }}_{\rm{p}}},{\boldsymbol{\alpha }}} \right) \coloneqq \int_\mathcal{A} {{{\boldsymbol{A}}^{\rm{T}}}\underline {{{\left. {{\partial _{\boldsymbol{E}}}\Psi } \right|}_{{\boldsymbol{E}} = {\boldsymbol{E}}^{{\rm{tot}}}}}} \,{j_0}\,{\rm{d}}\mathcal{A}},
	\end{align}
	\begin{align}
		{\partial _{\boldsymbol{\alpha }}}\psi \left( {{{\boldsymbol{\varepsilon }}_{\rm{p}}},{\boldsymbol{\alpha }}} \right) \coloneqq \int_\mathcal{A} {{{\boldsymbol{\Gamma }}^{\rm{T}}}\underline {{{\left. {{\partial _{\boldsymbol{E}}}\Psi } \right|}_{{\boldsymbol{E}} = {\boldsymbol{E}}^{{\rm{tot}}}}}} \,{j_0}\,{\rm{d}}\mathcal{A}},
	\end{align}
	with ${{\partial_{\boldsymbol{E}}}\Psi }\coloneqq{{\partial \Psi(\boldsymbol{E}) }}/{{\partial {\boldsymbol{E}}}}$.
\end{subequations}
The detailed expression of the operator $\mathbb{B}({\boldsymbol{y}})$, defined by
\begin{align}
	\delta {\boldsymbol{\varepsilon }} = {\mathbb{B}}({\boldsymbol{y}})\,\delta {\boldsymbol{y}},
\end{align}
is given in \citet[Eq.\,(A.4.5)]{choi2021isogeometric}.
\subsubsection{Work function due to external loads}
We define a \textit{work function} for \textit{time-dependent} external loads,  
\begin{align}
	\label{pot_ext_load}
	{W_{{\rm{ext}}}}(\boldsymbol{y},t)\!&\coloneqq\! \!\int_0^L {{\boldsymbol{y}} \cdot {\boldsymbol{\bar R}}(t)\,\,{\rm{d}}s}  + {\left[ {{\boldsymbol{y}} \cdot {{{\boldsymbol{\bar R}}}_0}(t)} \right]_{s \in {\Gamma _{\rm{N}}}}},\,\,t\in\left[0,T\right],
\end{align}
where \citep{choi2021isogeometric}
\begin{equation}
	{\boldsymbol{\bar R}} \coloneqq \left\{\renewcommand{\arraystretch}{1.25}{\begin{array}{*{20}{c}}
			{{\boldsymbol{\bar n}}}\\
			{{{{{{{\boldsymbol{\bar {\tilde m}}}}}}}\,\!^1}}\\
			{{{{{{{\boldsymbol{\bar {\tilde m}}}}}}}\,\!^2}}
	\end{array}} \right\},\,\,{{\boldsymbol{\bar R}}_0} \coloneqq \left\{\renewcommand{\arraystretch}{1.25} {\begin{array}{*{20}{c}}
			{{{{\boldsymbol{\bar n}}}_0}}\\
			{{{{ { {\boldsymbol{\bar{\tilde m}}}}}}}\,\!^1_0}\\
			{{{{{ {\boldsymbol{\bar{\tilde m}}}}}}}\,\!^2_0}
	\end{array}} \right\}.
\end{equation}
Here, $\boldsymbol{\bar{n}}$ and $\boldsymbol{\bar{\tilde m}}\:\!^\alpha$ ($\alpha\in\left\{1,2\right\}$) denote the given external (distributed) force and director couple along the center axis, respectively, that follow from the body force and the lateral-surface tractions. Further, $\boldsymbol{\bar{n}}_0$ and $\boldsymbol{\bar{\tilde m}}\,\!^\alpha_0$ denote the prescribed force and director couple at the Neumann boundary, $\Gamma_\mathrm{N}\subset{\emptyset\cup\left\{0,L\right\}}$, due to a given surface traction acting on the end faces. $\emptyset$ denotes the empty set. Then, by taking the first variation of the work function in Eq.\,(\ref{pot_ext_load}), the external virtual work at time $t\in\left[0,T\right]$ can be expressed as 
\begin{align}
	\label{ext_load_vs}
	{G_{{\rm{ext}}}}(\delta {\boldsymbol{y}},t)& = \delta{W_\mathrm{ext}}=\! \int_0^L {\delta {\boldsymbol{y}} \cdot {\boldsymbol{\bar R}}(t){\rm{d}}s} + {\left[ \delta {\boldsymbol{y}} \cdot {{{\boldsymbol{\bar R}}}_0}(t) \right]_{s \in {\Gamma _{\rm{N}}}}}.
\end{align}
\textcolor{black}{Note that, here, we assume that the external load is deformation-independent.} If the external load has no explicit time dependence, the system is \textit{conservative}, i.e., the total energy of the system is conserved \citep[Sections 1.7 and 1.8]{lanczos2012variational}, and we may define an external potential energy as
\begin{align}
	U_\mathrm{ext}(\boldsymbol{y}) = -W_\mathrm{ext}(\boldsymbol{y}).
\end{align}
\subsubsection{Kinetic energy}
%
We define an independent \textcolor{black}{(material)} velocity vector $\boldsymbol{v}$ that satisfies
\begin{align}
	\boldsymbol{v}=\boldsymbol{\dot x}. \label{compati_disp_vel}
\end{align}
From the beam kinematics in Eq.\,(\ref{bkin_velocity}), we may also decompose the independent velocity field into two parts: the center axis velocity $\boldsymbol{v}_\varphi\equiv\boldsymbol{v}_\varphi(s,t)\in\mathbb{R}^3$ and the director velocity $\boldsymbol{v}_\alpha\equiv\boldsymbol{v}_\alpha(s,t)\in\mathbb{R}^3$ ($\alpha\in\left\{1,2\right\}$) that satisfy the \textit{compatibility} condition
\begin{align}
	\label{const_compat_v_y_v_d}
	{\boldsymbol{V}} = {{\boldsymbol{\dot y}}}, 
\end{align}
where, for brevity, we define a \textit{generalized velocity} ${\boldsymbol{V}} \coloneqq {\left[ {{\boldsymbol{v}}_\varphi ^{\rm{T}},{\boldsymbol{v}}_1^{\rm{T}},{\boldsymbol{v}}_2^{\rm{T}}} \right]^{\rm{T}}}\!\in\!\mathbb{R}^d$. Then, we can express the \textit{total linear momentum} as
\begin{align}
	\label{def_tot_lin_momentum}
	{\boldsymbol{L}}(\boldsymbol{V}) =\int_0^L {{\boldsymbol{p}}(\boldsymbol{V})\,{\rm{d}}s},
\end{align}
where we have defined the \textit{local linear momentum} per undeformed unit arc-length,
\begin{equation}
	\label{vec_lin_momenta}
	{\boldsymbol{p}}(\boldsymbol{V}) \coloneqq {\rho _A}{\boldsymbol{v}}_\varphi + \,I_\rho ^\alpha {{\boldsymbol{v}}_\alpha},
\end{equation}
with the mass per undeformed unit arc-length, i.e., the line mass density,
\begin{equation}
	{\rho_{A}} \coloneqq \int_\mathcal{A} {{\rho _0}\,{j_0}\,{\rm{d}}\mathcal{A}},
\end{equation}
and the \textit{first mass moment of inertia} per undeformed unit arc-length,
\begin{equation}
	I_\rho ^\alpha  \coloneqq \int_\mathcal{A} {{\rho _0}\,{\zeta ^\alpha }{j_0}\,{\rm{d}}\mathcal{A}},\,\,\alpha\in\left\{1,2\right\}.
\end{equation}
We can then express the \textit{total angular momentum}, as
\begin{align}	
	\label{def_tot_ang_momentum}
	{\boldsymbol{J}}(\boldsymbol{y},\boldsymbol{V}) = \int_0^L {\left\{{\boldsymbol{\varphi }} \times {\boldsymbol{p}}(\boldsymbol{V}) + {{\boldsymbol{d}}_\alpha } \times {{\boldsymbol{\mu }}^\alpha }(\boldsymbol{V}) \right\}\,{\rm{d}}s},
\end{align}
where we have defined the \textit{local director momentum} per undeformed unit arc-length,
\begin{equation}
	\label{vec_dir_momenta}	
	{{\boldsymbol{\mu }}^\alpha }(\boldsymbol{V}) \coloneqq I_\rho ^\alpha {\boldsymbol{v}}_\varphi + I_\rho ^{\alpha \beta }{{\boldsymbol{v}}_\beta},\,\,\alpha\in\left\{1,2\right\},
\end{equation}
with the \textit{second mass moment of inertia} per undeformed unit arc-length,
\begin{equation}
	I_\rho ^{\alpha\beta}  \coloneqq \int_\mathcal{A} {{\rho _0}\,{\zeta ^\alpha }{\zeta ^\beta }{j_0}\,{\rm{d}}\mathcal{A}},\,\,\alpha,\beta\in\left\{1,2\right\}.
\end{equation}
For brevity, we combine Eqs.\,(\ref{vec_lin_momenta}) and (\ref{vec_dir_momenta}) into 
\begin{subequations}
	\begin{equation}
		\left\{ {\renewcommand{\arraystretch}{1.25}\begin{array}{*{20}{c}}
				{{\boldsymbol{p}}}(\boldsymbol{V})\\
				{{{\boldsymbol{\mu }}^1}}(\boldsymbol{V})\\
				{{{\boldsymbol{\mu }}^2}}(\boldsymbol{V})
		\end{array}} \right\} = {\boldsymbol{\mathcal{M}}}{\boldsymbol{V}},
	\end{equation}
	with the mass moment of inertia
	\begin{equation}
		{\boldsymbol{\mathcal{M}}} \coloneqq \left[ {\renewcommand{\arraystretch}{1.25}\begin{array}{*{20}{c}}
				{{\rho _A}}{\boldsymbol{1}}&{I_\rho ^1}{\boldsymbol{1}}&{I_\rho ^2}{\boldsymbol{1}}\\
				{}&{I_\rho ^{11}}{\boldsymbol{1}}&{I_\rho ^{12}}{\boldsymbol{1}}\\
				{{\rm{sym}}{\rm{.}}}&{}&{I_\rho ^{22}}{\boldsymbol{1}}
		\end{array}} \right],
	\end{equation}
	which is time-independent. It should be noted that $\boldsymbol{\mathcal{M}}$ is symmetric (i.e., $I^{12}_\rho\!=\!I^{21}_\rho$), regardless of the cross-section's shape and initial mass distribution.
\end{subequations}
Further, we define the total kinetic energy, as
\begin{align}
	\label{tot_kin_e_beam}
	\mathcal{K}(\boldsymbol{V}) = \dfrac{1}{2}\int_0^L {\boldsymbol{V}} \cdot {\boldsymbol{\mathcal{M}}}{\boldsymbol{V}}\,{\rm{d}}s.
\end{align}
Accordingly, we have an \textit{inertial contribution} to the virtual work,
\begin{equation}
	\label{iner_cont_vwork}
	{G_{{\rm{iner}}}}({\boldsymbol{\dot V}},\delta {\boldsymbol{y}}) \coloneqq \int_0^L {\delta {\boldsymbol{y}} \cdot {\boldsymbol{\mathcal{M}}}{\boldsymbol{\dot V}}\,{\rm{d}}s}.
\end{equation}
\subsection{Euler-Lagrange equations: Livens' theorem}
\textcolor{black}{Based on Livens' theorem \citep[Sec. 26.2]{pars1965treatise}, we define the Lagrangian with the independent velocity field $\boldsymbol{V}$, as}
\begin{align}
	\label{Lagrangian_def_f}
	\mathcal{L} &\coloneqq \mathcal{K}({\boldsymbol{V}}) - {{{U}_\mathrm{HW}}\left( {{\boldsymbol{\varepsilon }}({\boldsymbol{y}}),{{\boldsymbol{\varepsilon }}_{\rm{p}}},{{\boldsymbol{r}}_{\rm{p}}},{\boldsymbol{\alpha }}} \right) + {W_\mathrm{ext}}({\boldsymbol{y}},t)} + \int_0^L {({{\boldsymbol{\dot y}}} - {\boldsymbol{V}}) \cdot \left\{ {\renewcommand{\arraystretch}{1.25}\begin{array}{*{20}{c}}
				{{{\boldsymbol{\lambda }}_\varphi }}\\
				{{{\boldsymbol{\lambda }}^1}}\\
				{{{\boldsymbol{\lambda }}^2}}
		\end{array}} \right\}{\rm{d}}s},
\end{align}
\textcolor{black}{where the last term represents a weak enforcement of the compatibility condition in Eq.\,(\ref{const_compat_v_y_v_d}).} Here, the Lagrange multipliers ${\boldsymbol{\lambda}}_\varphi,{\boldsymbol{\lambda}}^\alpha\in\mathbb{R}^3, \alpha\in\left\{1,2\right\}$ can be identified by the stationarity of the Lagrangian with respect to the independent velocities ${\boldsymbol{v}}_\varphi$ and ${\boldsymbol{v}}_\alpha$ ($\alpha\in\left\{1,2\right\}$), as \citep[Section 3.2]{betsch2016energy}
\begin{subequations}
	\begin{align}
		{{\boldsymbol{\lambda }}_\varphi } &= {\boldsymbol{p}}({\boldsymbol{V}}),\label{lag_mult_ident_phi}\\
		{{\boldsymbol{\lambda }}^\alpha } &= {{\boldsymbol{\mu }}^\alpha }({\boldsymbol{V}}),\,\,\alpha  \in \left\{ {1,2} \right\}.\label{lag_mult_ident_alph}
	\end{align}
	\label{lag_mult_ident}
\end{subequations}
Substituting Eqs.\,(\ref{lag_mult_ident_phi})-(\ref{lag_mult_ident_alph}) into Eq.\,(\ref{Lagrangian_def_f}), the Lagrangian can be rewritten as
\begin{align}
	\mathcal{L} &= \int_0^L {\left( {{{\boldsymbol{\dot y}}} - \frac{1}{2}{{\boldsymbol{V}}}} \right) \cdot \boldsymbol{\mathcal{M}}{\boldsymbol{V}}\,{\rm{d}}s}- {{{U}_\mathrm{HW}}\left( {{\boldsymbol{\varepsilon }}({\boldsymbol{y}}),{{\boldsymbol{\varepsilon }}_{\rm{p}}},{{\boldsymbol{r}}_{\rm{p}}},{\boldsymbol{\alpha }}} \right) + {W_\mathrm{ext}}({\boldsymbol{y}},t)},
\end{align}
from which we define an action integral in the time span $\left[0,T\right]\ni{t}$,
\begin{equation}
	\mathcal{S} \coloneqq \int_{{0}}^{{T}} {{\mathcal{L}}(\boldsymbol{V},\boldsymbol{y},\boldsymbol{\varepsilon}_\mathrm{p},\boldsymbol{r}_\mathrm{p},\boldsymbol{\alpha},t)\,{\rm{d}}t}.
\end{equation}
By the stationarity condition of the action integral, we obtain the following weak form of the initial-boundary-value problem: Find the solution $\left\{\boldsymbol{V},\boldsymbol{y},\boldsymbol{r}_\mathrm{p},\boldsymbol{\varepsilon}_\mathrm{p},\boldsymbol{\alpha}\right\}\in{{\bar {\mathcal{V}}}\times\mathcal{V}\times\mathcal{V}_\mathrm{p}\times\mathcal{V}_\mathrm{p}\times\mathcal{V}_\mathrm{a}}$ at time $t\in(0,T)$ such that 
\begin{subequations}
	\label{eqs_wform_mbal_compat}
	\begin{align}
		\label{el_eq_var}
		{G_{{\rm{iner}}}}\left({\boldsymbol{\dot V}},\delta {\boldsymbol{y}}\right) + {G_\mathrm{int}^\mathrm{y}}({\boldsymbol{y}},{\boldsymbol{r}}_\mathrm{p},\delta {\boldsymbol{y}}) &= {G_{{\rm{ext}}}}(\delta {\boldsymbol{y}},t),\\
		{G_\mathrm{int}^\mathrm{r}}({\boldsymbol{y}},{\boldsymbol{\varepsilon}}_\mathrm{p},\delta {\boldsymbol{r}}_\mathrm{p}) &= 0,\\
		{G_\mathrm{int}^\varepsilon}({\boldsymbol{r}}_\mathrm{p},{\boldsymbol{\varepsilon}}_\mathrm{p},{\boldsymbol{\alpha}},\delta {\boldsymbol{\varepsilon}}_\mathrm{p}) &= 0,\\		
		{G_\mathrm{int}^\mathrm{a}}({\boldsymbol{\varepsilon}}_\mathrm{p},{\boldsymbol{\alpha}},\delta {\boldsymbol{\alpha}}) &= 0,\label{el_eq_var_a}
	\end{align}
	are satisfied for all $\left\{\delta\boldsymbol{y},\delta\boldsymbol{r}_\mathrm{p},\delta\boldsymbol{\varepsilon}_\mathrm{p},\delta\boldsymbol{\alpha}\right\} \in \bar {\mathcal{V}}\times {\mathcal{V}}_\mathrm{p}\times{{\mathcal{V}}_\mathrm{p}}\times{\mathcal{V}_\mathrm{a}}$, together with the \textit{compatibility} condition from Eq.\,(\ref{const_compat_v_y_v_d}),
\end{subequations}
\begin{equation}	
	\label{compat_integ_v_y}
	\int_0^L {\delta {\boldsymbol{V}} \cdot {\boldsymbol{\mathcal{M}}}\left({\boldsymbol{\dot y}}-{\boldsymbol{V}}\right)\,{\rm{d}}s}  = 0,\,\,\forall{\delta \boldsymbol{V}}\in{\bar{\mathcal{V}}}, 
\end{equation}
and the initial conditions at $s\in\left(0,L\right)$,
\begin{subequations}
	\begin{alignat}{5}
		{\boldsymbol{\varphi}}(s,0) &= {\boldsymbol{\varphi}}_0(s),\quad{\boldsymbol{d}}_\alpha(s,0) &&= {{\boldsymbol{D}}_\alpha}(s),\quad\alpha\in\left\{1,2\right\},\\
		{\boldsymbol{v}_{\varphi}}(s,0) &= {{\boldsymbol{\bar v}}_{\varphi}}(s),\quad{\boldsymbol{v}}_\alpha(s,0) &&= {{\boldsymbol{\bar v}}_{\alpha}}(s),\quad\alpha\in\left\{1,2\right\}.
	\end{alignat}
\end{subequations}
Here, ${{\boldsymbol{\bar v}}_{\varphi}}(s)$ and ${{\boldsymbol{\bar v}}_{\alpha}}(s)$ denote the given initial velocity fields for the center axis and directors, respectively. We define the kinematically admissible solution space for the configuration variable $\boldsymbol{y}$, as
\begin{align}
	\mathcal{V} \coloneqq &\Bigl\{ \left. {\boldsymbol{y} \in {{\left[ {{H^1}(0,L)} \right]}^d}} \right|{\boldsymbol{\varphi}} = {{{\boldsymbol{\bar \varphi}}}},\,\,{{\boldsymbol{d}}_1} = {{\boldsymbol{\bar d}}_{1}},\,{\rm{and}}\,\,{{\boldsymbol{d}}_2} = {{{\boldsymbol{\bar d}}}_{2}}\,\,{\rm{at}}\,\,s \in {\Gamma _{\rm{D}}}\, \Bigl\},
\end{align}
where $\Gamma_\mathrm{D}\subset\emptyset\cup{\left\{0,L\right\}}$ denotes the Dirichlet boundary such that $\Gamma_\mathrm{D}\cup\Gamma_\mathrm{N}=\left\{0,L\right\}$ and $\Gamma_\mathrm{D}\cap\Gamma_\mathrm{N}=\emptyset$. The variational space is defined as
\begin{align}
	\mathcal{\bar V} \coloneqq &\Bigl\{ \left. {\delta \boldsymbol{y} \in {{\left[ {{H^1}(0,L)} \right]}^d}} \right|\delta{\boldsymbol{\varphi}} = {\delta{\boldsymbol{d}}_1} = \delta{{\boldsymbol{d}}_2} = {{{\boldsymbol{0}}}}\,\,{\rm{at}}\,\,s \in {\Gamma _{\rm{D}}}\, \Bigl\}.
\end{align}
Further, we employ the solution spaces for the physical variables, ${\mathcal{V}}_\mathrm{p} \coloneqq\left[L^2(0,L)\right]^{d_\mathrm{p}}$, and ${\mathcal{V}}_\mathrm{a}\coloneqq\left[L^2(0,L)\right]^{d_\mathrm{a}}$ such that those solution functions may have discontinuity in the domain. Further comment on the selection of the approximation space is given in Sections \ref{approx_fe_phy_kin_var} and \ref{approx_fe_phy_enh_var}.
\section{Temporal discretization: an EMC scheme}
\label{sec_emc_tstep}
\textcolor{black}{We divide a given time span $[0,T]$ into $N$ intervals, and find approximate solutions at discrete times $\left\{t_0=0,t_1,\cdots,t_n,t_{n+1},\cdots,t_N = T\right\}$, where $n$ is a non-negative integer. Note that time increment $\Delta{t}\coloneqq{t_{n+1}-t_n}$ can vary adaptively.} We denote the approximation of the independent solution variables at time $t_n$ by ${}^n(\bullet)$. In the given time interval $[t_n,t_{n+1}]$, our goal is to find an approximate solution at time ${t}_{n+1}$ for a given solution at time $t_n$. We present an implicit time-stepping scheme, based on the previous work in \citet{betsch2016energy}, which exactly preserves the total energy for linear constitutive laws. A further extension to nonlinear constitutive laws, e.g., using the discrete derivative of \citet{gonzalez2000exact} remains future work. The present time-stepping scheme turns out to exhibit superior numerical stability in comparison to standard trapezoidal and mid-point rules. We first approximate the first-order time derivatives at the mid-point, as
\begin{subequations}
	\begin{align}
		\label{approx_tderiv_y}
		{{\boldsymbol{\dot y}}}\!\left({{t_n} + \frac{1}{2}\Delta t}\right) &\approx \frac{1}{{\Delta t}}\left( {{}^{n + 1}{\boldsymbol{y}} - {}^n{\boldsymbol{y}}} \right),
	\end{align}
	and
	\begin{align}
		\label{approx_tderiv_V}
		{{\boldsymbol{\dot V}}}\!\left({{t_n} + \frac{1}{2}\Delta t}\right) &\approx \frac{1}{{\Delta t}}\left( {{}^{n + 1}{\boldsymbol{V}} - {}^n{\boldsymbol{V}}} \right).
	\end{align}
\end{subequations}
Further, we define the mid-point approximations
\begin{subequations}
	\begin{alignat}{3}
		{\boldsymbol{y}}\!\left({t_n} + \frac{1}{2}\Delta t\right) &\approx \frac{1}{2}\left( {{}^{n + 1}{\boldsymbol{y}} + {}^n{\boldsymbol{y}}} \right) &&\eqqcolon {}^{n + \frac{1}{2}}{\boldsymbol{y}},\\
		{{\boldsymbol{r}}_{\rm{p}}}\!\left({t_n} + \frac{1}{2}\Delta t\right) &\approx  \frac{1}{2}\left( {{}^{n + 1}{{\boldsymbol{r}}_{\rm{p}}} + {}^n{{\boldsymbol{r}}_{\rm{p}}}} \right)&&\eqqcolon {}^{n + \frac{1}{2}}{{\boldsymbol{r}}_{\rm{p}}},\\
		{{\boldsymbol{\varepsilon }}_{\rm{p}}}\!\left({t_n} + \frac{1}{2}\Delta t\right) &\approx \frac{1}{2}\left( {{}^{n + 1}{{\boldsymbol{\varepsilon }}_{\rm{p}}} + {}^n{{\boldsymbol{\varepsilon }}_{\rm{p}}}} \right)&&\eqqcolon  {}^{n + \frac{1}{2}}{{\boldsymbol{\varepsilon }}_{\rm{p}}},\\
		{{\boldsymbol{\alpha}}}\!\left({t_n} + \frac{1}{2}\Delta t\right) &\approx \frac{1}{2}\left( {{}^{n + 1}{{\boldsymbol{\alpha}}} + {}^n{{\boldsymbol{\alpha}}}} \right)&&\eqqcolon  {}^{n + \frac{1}{2}}{{\boldsymbol{\alpha}}},
	\end{alignat}
	and
	\begin{equation}
			\label{approx_vel_mid_pt}
		{\boldsymbol{V}}\!\left( {{t_n} + \frac{1}{2}\Delta t} \right) \approx \frac{1}{2}\left( {{}^{n + 1}{\boldsymbol{V}} + {}^n{\boldsymbol{V}}} \right)\eqqcolon{}^{n + \frac{1}{2}}{\boldsymbol{V}}.		
	\end{equation}
\end{subequations}
From the compatibility condition in Eq.\,(\ref{compat_integ_v_y}) at the mid-point, combining Eqs.\,(\ref{approx_tderiv_y}) and (\ref{approx_vel_mid_pt}) gives \citep[Sec.\,3.8]{betsch2016energy}
\begin{align}
	\label{rep_cur_V_conf}
	{}^{n + 1}{\boldsymbol{V}} = \frac{2}{{\Delta t}}\left( {{}^{n + 1}{\boldsymbol{y}} - {}^n{\boldsymbol{y}}} \right) - {}^n{\boldsymbol{V}}.
\end{align}
Then, substituting Eq.\,(\ref{rep_cur_V_conf}) into Eq.\,(\ref{approx_tderiv_V}) yields
\begin{align}
	\label{rep_cur_V_conf_1}	
	{{\boldsymbol{\dot V}}}\!\left({{t_n} + \frac{1}{2}\Delta t}\right) &\approx \frac{2}{{\Delta {t^2}}}\left( {{}^{n + 1}{\boldsymbol{y}} - {}^n{\boldsymbol{y}}} \right) - \frac{2}{{\Delta t}}{}^n{\boldsymbol{V}}.
\end{align}
Substituting Eq.\,(\ref{rep_cur_V_conf_1}) into Eq.\,(\ref{el_eq_var}), a time-discrete form of the momentum balance equation at the mid-point can be stated as:
\begin{subequations}
	\label{eqs_wform_mbal_compat_mid_pt}
	\begin{align}
		&\frac{2}{{\Delta {t^2}}}{G_{{\rm{iner}}}}\left( {{}^{n + 1}{\boldsymbol{y}},\delta {\boldsymbol{y}}} \right) + {G^\mathrm{y}_{{\mathop{\rm int}} }}\left({}^{n+\frac{1}{2}}{\boldsymbol{y}},{}^{n+\frac{1}{2}}{\boldsymbol{r}}_\mathrm{p},\delta {\boldsymbol{y}}\right) \nonumber\\
		&= {G_{{\rm{ext}}}}\!\left(\delta {\boldsymbol{y}},{t_n}+\frac{1}{2}{\Delta t}\right) + \frac{2}{{\Delta {t}}}{G_{{\rm{iner}}}}\!\left( {{}^n{\boldsymbol{V}}+\frac{1}{{\Delta t}}{{}^n{\boldsymbol{y}}},\delta {\boldsymbol{y}}} \right),\label{eq_var_eq_y}
	\end{align}
	which is subject to
	\begin{align}		
		{G^\mathrm{r}_{{\mathop{\rm int}}}}\left( {{}^{n + 1}{{\boldsymbol{y}}},{}^{n + 1}{\boldsymbol{\varepsilon}}_{\rm{p}},\delta {{\boldsymbol{r}}_{\rm{p}}}} \right) &= 0,\label{t_emc_compat_cond}\\
		G_{{\mathop{\rm int}} }^\varepsilon\left( {{}^{n + \frac{1}{2}}{{\boldsymbol{r}}_{\rm{p}}},{}^{n + \frac{1}{2}}{{\boldsymbol{\varepsilon }}_{\rm{p}}},{}^{n + \frac{1}{2}}{\boldsymbol{\alpha }},\delta {{\boldsymbol{\varepsilon }}_{\rm{p}}}} \right)&=0,\label{t_emc_constitutiv_cond0}\\
		G_{{\mathop{\rm int}} }^\mathrm{a}\left( {{}^{n + \frac{1}{2}}{{\boldsymbol{\varepsilon }}_{\rm{p}}},{}^{n + \frac{1}{2}}{\boldsymbol{\alpha }},\delta {{\boldsymbol{\alpha}}}} \right)&=0,\label{t_emc_enh_a_cond0}
	\end{align}
	$\forall\left\{\delta\boldsymbol{y},\delta\boldsymbol{r}_\mathrm{p},\delta\boldsymbol{\varepsilon}_\mathrm{p},\delta\boldsymbol{\alpha}\right\} \in \bar {\mathcal{V}}\times {\mathcal{V}}_\mathrm{p}\times{{\mathcal{V}}_\mathrm{p}}\times{\mathcal{V}_\mathrm{a}}$.
\end{subequations}	
It should be noted that, in Eq.\,(\ref{t_emc_compat_cond}), we have imposed the compatibility between the displacement and the physical strain at time $t=t_{n+1}$ instead of $t=t_{n}+\frac{1}{2}\Delta{t}$, which is sufficient for the exact conservation of total energy \textcolor{black}{for linear constitutive laws} \citep{betsch2016energy}. Note that this treatment is in line with using the algorithmic stress resultant in \citet{romero2002objective}, which evaluates the stress resultant by taking an average between two subsequent discrete times, instead of evaluating it from the nonlinear strain-configuration relation at the mid-point.
\subsection{Linearization}
Due to the geometrical and material nonlinearities in the internal virtual work term of Eq.\,(\ref{vform_hw}), we need an incremental-iterative (Newton-Raphson) solution process. It can be stated as: In the $(n+1)$th time step, for the given solution at the previous $n$th time step $\left\{{}^n{\boldsymbol{V}},{}^n{\boldsymbol{y}},{}^n{\boldsymbol{r}}_\mathrm{p},{}^n{\boldsymbol{\varepsilon}}_\mathrm{p},{}^n{\boldsymbol{\alpha}}\right\}$, and that in the previous iteration
\begin{equation*}
	\left\{{}^{n+1}{\boldsymbol{y}}^{(i-1)},{}^{n+1}{\boldsymbol{r}}_\mathrm{p}^{(i-1)},{}^{n+1}{\boldsymbol{\varepsilon}}^{(i-1)}_\mathrm{p},{}^{n+1}{\boldsymbol{\alpha}}^{(i-1)}\right\},
\end{equation*}	
find the solution increment 
\begin{equation*}
	\left\{\Delta{\boldsymbol{y}},\Delta{\boldsymbol{r}}_\mathrm{p},\Delta{\boldsymbol{\varepsilon}}_\mathrm{p},\Delta{\boldsymbol{\alpha}}\right\}\in{\bar {\mathcal{V}} \times {\mathcal{V}_{\rm{p}}} \times {\mathcal{V}_{\rm{p}}} \times {\mathcal{V}_{\rm{a}}}}
\end{equation*}
such that
\begin{align}
	\label{lin_var_eq_emc}
	&\int_0^L {\left\{ {\begin{array}{*{20}{c}}
				{\delta {\boldsymbol{y}}}\\
				{\delta {{\boldsymbol{r}}_{\rm{p}}}}\\
				{\delta {{\boldsymbol{\varepsilon }}_\mathrm{p}}}\\
				{\delta {\boldsymbol{\alpha }}}
		\end{array}} \right\} \cdot {\boldsymbol{k}}\,\left\{ {\begin{array}{*{20}{c}}
				{\Delta {\boldsymbol{y}}}\\
				{\Delta {{\boldsymbol{r}}_{\rm{p}}}}\\
				{\Delta {{\boldsymbol{\varepsilon }}_\mathrm{p}}}\\
				{\Delta {\boldsymbol{\alpha }}}
		\end{array}} \right\}{\rm{d}}s} +\frac{4}{{\Delta t}^2}\int_0^L{\delta{\boldsymbol{y}}\cdot\boldsymbol{\mathcal{M}}\Delta\boldsymbol{y}}\,\mathrm{d}s\nonumber\\	
	&=2\,{G_{{\rm{ext}}}}\left( {\delta {\boldsymbol{y}},{t_n} + \frac{1}{2}\Delta t} \right) + \frac{4}{{\Delta t}}{G_{{\rm{iner}}}}\left( {{}^n{\boldsymbol{V}} - \frac{1}{{\Delta t}}\left( {{}^{n + 1}{{\boldsymbol{y}}^{(i - 1)}} - {}^n{\boldsymbol{y}}} \right),\delta {\boldsymbol{y}}} \right)\nonumber\\
	& - 2\,G_{{\mathop{\rm int}} }^{\rm{y}}\left( {{{\boldsymbol{\bar y}}},{\boldsymbol{\bar r}}_{\rm{p}},\delta {\boldsymbol{y}}} \right)-G_{{\mathop{\rm int}} }^{\rm{r}}\left( {{}^{n + 1}{{\boldsymbol{y}}^{(i - 1)}},{}^{n + 1}{\boldsymbol{\varepsilon }}_{\rm{p}}^{(i - 1)},\delta {{\boldsymbol{r}}_{\rm{p}}}} \right)\nonumber\\
	&-2\,G_{{\mathop{\rm int}}}^{\varepsilon} \left( {{\boldsymbol{\bar r}}_\mathrm{p}},{\boldsymbol{\bar \varepsilon }}_{\rm{p}},{\boldsymbol{\bar \alpha}},\delta {\boldsymbol{\varepsilon}}_\mathrm{p} \right)\nonumber\\
	&-2\,G_{{\mathop{\rm int}} }^{\rm{a}}\left({\boldsymbol{\bar \varepsilon }}_{\rm{p}},{\boldsymbol{\bar \alpha}},\delta {\boldsymbol{\alpha}}\right),
\end{align}
$\forall\left\{\delta{\boldsymbol{y}},\delta{\boldsymbol{r}}_\mathrm{p},\delta{\boldsymbol{\varepsilon}}_\mathrm{p},\delta{\boldsymbol{\alpha}}\right\}\in{\bar {\mathcal{V}} \times {\mathcal{V}_{\rm{p}}} \times {\mathcal{V}_{\rm{p}}} \times {\mathcal{V}_{\rm{a}}}}$, and update
\begin{subequations}
	\begin{alignat}{5}
		{}^{n+1}{{\boldsymbol{y}}^{(i)}} &= {}^{n+1}{{\boldsymbol{y}}^{(i - 1)}} &&+ \Delta {\boldsymbol{y}},\,\,&&{}^{n+1}{{\boldsymbol{y}}^{(0)}}&&\equiv{}^{n}{{\boldsymbol{y}}},\\
		{}^{n+1}{{\boldsymbol{r }}_\mathrm{p}^{(i)}} &= {}^{n+1}{{\boldsymbol{r }}_\mathrm{p}^{(i - 1)}} &&+ \Delta {\boldsymbol{r }}_\mathrm{p},\,\,&&{}^{n+1}{{\boldsymbol{r }}_\mathrm{p}^{(0)}}&&\equiv{}^{n}{{\boldsymbol{r }}_\mathrm{p}},\\
		{}^{n+1}{{\boldsymbol{\varepsilon}}_\mathrm{p}^{(i)}} &= {}^{n+1}{{\boldsymbol{\varepsilon }}_\mathrm{p}^{(i - 1)}}\!\!&&+ \Delta {\boldsymbol{\varepsilon }}_\mathrm{p},\,\,&&{}^{n+1}{{\boldsymbol{\varepsilon }}_\mathrm{p}^{(0)}}&&\equiv{}^{n}{{\boldsymbol{\varepsilon }}}_\mathrm{p},\\
		{}^{n+1}{{\boldsymbol{\alpha }}^{(i)}} &= {}^{n+1}{{\boldsymbol{\alpha }}^{(i - 1)}}&&+ \Delta {\boldsymbol{\alpha }},\,\,&&{}^{n+1}{{\boldsymbol{\alpha }}^{(0)}}&&\equiv{}^{n}{{\boldsymbol{\alpha }}},
	\end{alignat}
\end{subequations}
for $i=1,2,\cdots$, until a convergence criteria is satisfied. After the iteration finishes, we should also update the independent velocity field using Eq.\,(\ref{rep_cur_V_conf}). In Eq.\,(\ref{lin_var_eq_emc}), we have the tangent operator
\begin{align}
	\label{tstiff_emc_unsym}
	&{{\boldsymbol{k}}}\!\coloneqq\!\left[\renewcommand{\arraystretch}{1.7}\begin{array}{*{20}{c}}
		{{{\boldsymbol{Y}}^{\rm{T}}}{{\boldsymbol{k}}_{\rm{G}}}\!\left( {\boldsymbol{\bar r}}_{\rm{p}} \right)\!{\boldsymbol{Y}}}&{{{\mathbb{B}}}}{\left( {\boldsymbol{\bar y}} \right)}\:\!\!^{\rm{T}}&{{{\bf{0}}_{d \times {d_{\rm{p}}}}}}&{{{\bf{0}}_{d \times {d_{\rm{a}}}}}}\\
		{{{\mathbb{B}}}\!\left( {{}^{n + 1}{{\boldsymbol{y}}^{(i - 1)}}} \right)}&{{{\bf{0}}_{{d_{\rm{p}}} \times {d_{\rm{p}}}}}}&{ - {{\bf{1}}_{{d_{\rm{p}}}\times{d_{\rm{p}}}}}}&{{{\bf{0}}_{{d_{\rm{p}}} \times {d_{\rm{a}}}}}}\\
		{{{\bf{0}}_{{d_{\rm{p}}} \times d}}}&{ - {{\bf{1}}_{{d_{\rm{p}}}\times{d_{\rm{p}}}}}}&{{\mathbb{C}}_\mathrm{p}^{\varepsilon \varepsilon}}\!\left( {{{{\boldsymbol{\bar \varepsilon }}}_{\rm{p}}},{\boldsymbol{\bar \alpha }}} \right)&{{\mathbb{C}}_\mathrm{p}^{{\rm{a}}\varepsilon}}\!\left( {{{{\boldsymbol{\bar \varepsilon }}}_{\rm{p}}},{\boldsymbol{\bar \alpha }}} \right)^\mathrm{T}\\
		{{{\bf{0}}_{{d_{\rm{a}}} \times d}}}&{{{\bf{0}}_{{d_{\rm{a}}} \times {d_{\rm{p}}}}}}&{{\mathbb{C}}_\mathrm{p}^{{\rm{a}}\varepsilon}}\left( {{{{\boldsymbol{\bar \varepsilon }}}_{\rm{p}}},{\boldsymbol{\bar \alpha }}} \right)&{{\mathbb{C}}_\mathrm{p}^{{\rm{aa}}}}\left( {{{{\boldsymbol{\bar \varepsilon }}}_{\rm{p}}},{\boldsymbol{\bar \alpha }}} \right)
	\end{array} \right],\!
\end{align} 
where we have defined the constitutive matrices, as \citep{choi2023selectively}
\begin{subequations}
	\begin{alignat}{3}
		\frac{{{\partial ^2}\psi }}{{\partial {{\boldsymbol{\varepsilon }}_{\rm{p}}}\partial {{\boldsymbol{\varepsilon }}_{\rm{p}}}}} &= \int_\mathcal{A} {{{\boldsymbol{A}}^{\rm{T}}}\:\!{{\underline{\underline{{\boldsymbol{\mathcal{C}}}_{\rm{p}}}}}}\:\!{\boldsymbol{A}}\,{j_0}}\,{\rm{d}}\mathcal{A}&&\eqqcolon{\mathbb{C}_{\rm{p}}^{\varepsilon\varepsilon}}\!\left( {{{{\boldsymbol{\varepsilon }}}_{\rm{p}}},{\boldsymbol{\alpha }}} \right),\\
		\frac{{{\partial ^2}\psi }}{{\partial {\boldsymbol{\alpha }}\partial {\boldsymbol{\alpha }}}}  &= \int_\mathcal{A} {{\boldsymbol{\Gamma}}^{\rm{T}}}\:\!\underline{\underline{\boldsymbol{\mathcal{C}}_{\rm{p}}}}\:\!\boldsymbol{\Gamma}\,{j_0}\,{\rm{d}}\mathcal{A}&&\eqqcolon {{\mathbb{C}}_\mathrm{p}^{{\rm{aa}}}}\!\left( {{{{\boldsymbol{\varepsilon }}}_{\rm{p}}},{\boldsymbol{\alpha }}} \right),\\
		\frac{{{\partial ^2}\psi }}{{\partial {\boldsymbol{\alpha }}\partial {{\boldsymbol{\varepsilon }}_{\rm{p}}}}} &= \int_\mathcal{A} {{\boldsymbol{\Gamma}}^{\rm{T}}}\:\!\underline{\underline{\boldsymbol{\mathcal{C}}_{\rm{p}}}}\:\!\boldsymbol{A}\,{j_0}\,{\rm{d}}\mathcal{A}&&\eqqcolon{{\mathbb{C}}_\mathrm{p}^{{\rm{a}}\varepsilon}}\!\left( {{{{\boldsymbol{\varepsilon }}}_{\rm{p}}},{\boldsymbol{\alpha }}} \right),
	\end{alignat}
\end{subequations}
based on the material elasticity tensor
\begin{equation}
	{{\boldsymbol{\mathcal{C}}}_{\rm{p}}} \coloneqq \frac{{{\partial ^2}\Psi({{\boldsymbol{E}}_{\mathrm{p}}})}}{{\partial {{\boldsymbol{E}}_{\mathrm{p}}}\,\partial {{\boldsymbol{E}}_{\mathrm{p}}}}}.
\end{equation}
Here, ${{\underline{\underline{(\bullet)}}}}$ denotes Voigt notation of a fourth-order tensor with major and minor symmetries. The detailed expressions of the operators ${\boldsymbol{k}}_\mathrm{G}$ and $\boldsymbol{Y}$ can be found in Eq.\,(A.4.9) and (A.4.11) of \citet{choi2021isogeometric}, respectively. In Eq.\,(\ref{tstiff_emc_unsym}), for brevity, we also employ the following notations for averaged solution variables
\begin{subequations}
	\label{mid_pt_def}
	\begin{align}
		{\boldsymbol{\bar y}}&\coloneqq\frac{1}{2}\left({{}^{n + 1}{{\boldsymbol{y}}^{(i - 1)}} + {}^n{\boldsymbol{y}}}\right),\label{mid_pt_def_ybar}\\
		{\boldsymbol{\bar r}}_\mathrm{p}&\coloneqq{\frac{1}{2}\left( {{}^{n + 1}{\boldsymbol{r}}_{\rm{p}}^{(i - 1)} + {}^n{{\boldsymbol{r}}_{\rm{p}}}} \right)},\\	
		{\boldsymbol{\bar \varepsilon}}_\mathrm{p}&\coloneqq{\frac{1}{2}\left( {{}^{n + 1}{\boldsymbol{\varepsilon}}_{\rm{p}}^{(i - 1)} + {}^n{{\boldsymbol{\varepsilon}}_{\rm{p}}}} \right)},\\		
		{\boldsymbol{\bar \alpha}}&\coloneqq{\frac{1}{2}\left( {{}^{n + 1}{\boldsymbol{\alpha}}^{(i - 1)} + {}^n{{\boldsymbol{\alpha}}}} \right)}.\label{mid_pt_def_alph}	
	\end{align}
\end{subequations}
Further, ${\bf{1}}_{n\times{m}}$ and ${\bf{0}}_{n\times{m}}$ denote the identity and zero matrices of dimension $n\times{m}$, respectively. 
\vspace{2mm}
\begin{remark}
	\label{rem_unsym_b_btr}
It should be noted that the tangent stiffness operator in Eq.\,(\ref{tstiff_emc_unsym}) is unsymmetric due to ${{\mathbb{B}}}(\boldsymbol{\bar y})\ne{{\mathbb{B}}}({}^{n+1}\boldsymbol{y}^{(i-1)})$, as we enforce the compatibility condition at time $t=t_{n+1}$ instead of $t=t_{n}+\frac{1}{2}{\Delta t}$ in Eq.\,(\ref{t_emc_compat_cond}). 
\end{remark}
\textcolor{black}{
\begin{remark}
In Section\,\ref{app_verif_emc_discrete_time}, we analytically verify the energy--momentum consistency of the present EMC scheme, which holds regardless of the element sizes in temporal and spatial domains.
\end{remark}}
\section{Spatial discretization: an isogeometric approach}
\label{spat_disc_iga}
\textcolor{black}{Next, we present the spatial discretization of the configuration and physical variables in the framework of an isogeometric analysis (IGA).} Especially, we show that the finite element approximation of the (total) director field leads to the objectivity of the approximated (discrete) beam strain. This basically necessitates an approximation of the initial director field.  
\subsection{NURBS curve}
In the isogeometric approach, we employ the same spline basis functions in the Computer-Aided Design (CAD) model for the analysis. This enables (i) an exact representation of the geometry of the beam's initial center axis \textcolor{black}{including conic sections like a circle}, and (ii) higher-order inter-element continuity which leads to the superior per DOF accuracy, compared with the conventional $C^0$--finite elements, see, e.g., \citet{choi2023selectively}. Using NURBS (Non-Uniform Rational B-Spline), the beam's initial center axis can be represented by
\begin{equation}\label{beam_curve_pos_nurbs}
	\boldsymbol{\varphi}_0(s(\xi))\equiv{\boldsymbol{X}}(\xi ) = \sum\limits_{I = 1}^{n_{\mathrm{cp}}} {{N^p_I}(\xi )\,{{\bf{P}}_{\!I}}},
\end{equation}
where ${\bf{P}}_{\!I}$ denotes the position of the $I$th control point, $N^p_I(\xi)$ denotes the corresponding NURBS basis function of degree $p$, and $n_\mathrm{cp}$ denotes the total number of basis functions (control points) in the given patch. Here, $\xi\in\varXi\subset{\mathbb R}$ denotes the parametric coordinate, and the parametric domain $\varXi\coloneqq\left[\xi_{1},\xi_{{n_{\mathrm{cp}}}+p+1}\right]$ is associated with the so called \textcolor{black}{\textit{knot vector}}, which is a given sequence of non-decreasing real numbers, ${\widetilde \varXi}={\left\{{\xi_1},{\xi_2},...,{\xi_{{n_{\mathrm{cp}}}+p+1}}\right\}}$, where ${\xi_i}\in{\mathbb R}$ denotes the $i$th knot. For a more detailed description on the NURBS basis functions, and mesh refinements, one may refer to \citet{piegl1996nurbs} and \citet{hughes2005isogeometric}. The arc-length coordinate along the initial center axis can be expressed by the mapping $s(\xi):{\varXi}\to{\left[0,L\right]}$, \citep{choi2021isogeometric}
\begin{equation}\label{beam_curve_pos_nurbs_alen_map}
	s(\xi )\coloneqq \int_{{\xi _1}}^{\eta  = \xi } {\left\| {{{\boldsymbol{X}}_{\!,\eta }}(\eta )} \right\|{\mathrm{d}}\eta },
\end{equation}
whose Jacobian is
\begin{align}\label{beam_curve_pos_nurbs_alen_map_jcb}
	\tilde j\coloneqq \frac{{{\mathrm{d}}s}}{{{\mathrm{d}}\xi }}= \left\| {{{\boldsymbol{X}}_{\!,\xi }}(\xi )} \right\|.
\end{align}
In the following section, we often use the notation ${N^p_{I,s}}$ for brevity, which is defined by
\begin{align}\label{nurbs_basis_deriv_wrt_s}
	{N^p_{I,s}} \coloneqq {N^p_{I,\xi }}\,\frac{{{\mathrm{d}}\xi }}{{{\mathrm{d}}s}} = \frac{1}{{\tilde j}}\,{N^p_{I,\xi }},
\end{align}
where ${N^p_{I,\xi }}$ denotes the differentiation of the basis function ${N^p_{I}(\xi)}$ with respect to $\xi$.
\subsubsection{Approximation of the initial director field}
\label{approx_init_dir_sec}
For a given continuous initial director field $\boldsymbol{D}_\alpha(\xi)$, we approximate (reconstruct) it using the NURBS basis functions of degree $p_\mathrm{d}$, as
\begin{align}
	\label{approx_dir_init_nurbs}
	{\boldsymbol{D}}_\alpha ^h(\xi ) = \sum\limits_{I = 1}^{n^\mathrm{d}_{\mathrm{cp}}} {N_I^{{p_{\rm{d}}}}(\xi )\,{{\bf{D}}_{\alpha I}}},\,\,\xi\in{\varXi},\,\,\alpha\in\left\{1,2\right\},
\end{align}
where $(\bullet)^h$ denotes a finite element approximation in space. \textcolor{black}{Here, the knot vector for the basis functions is given by ${{\widetilde \varXi}^\mathrm{d}}={\left\{{\xi^{\mathrm{d}}_1},{\xi^{\mathrm{d}}_2},...,{\xi^{\mathrm{d}}_{{n^\mathrm{d}_{\mathrm{cp}}}+p_\mathrm{d}+1}}\right\}}$.} We determine the control coefficients $\left\{{{\bf{D}}_{\alpha 1}},{{\bf{D}}_{\alpha 2}},\cdots,{{\bf{D}}_{\alpha {n^\mathrm{d}_\mathrm{cp}}}} \right\}$ ($\alpha\in\left\{1,2\right\}$) from the set of equations,
\begin{align}
	\label{approx_dir_init_nurbs_disc}
	{\boldsymbol{D}}_\alpha ^h({{\bar\xi}_J}) = {{\boldsymbol{D}}_\alpha }({{\bar \xi}_J}),\,\,J\in\left\{1,2,\cdots,{n}^\mathrm{d}_\mathrm{cp}\right\},
\end{align}
where the coordinate ${\bar \xi}_J$ is given by the Greville-abscissae (a moving average of knots) \citep[Chapter 10]{farin2000essentials}
\begin{align}
    \label{grev_abscissae_colloc}
	{\bar \xi _J} = \frac{1}{p_\mathrm{d}}\sum\limits_{i = 1}^{{p_{\rm{d}}}} {{\xi^\mathrm{d}_{J + i}}},\,\,J\in\left\{1,2,\cdots,{n}^\mathrm{d}_\mathrm{cp}\right\}.
\end{align}
\textcolor{black}{Note that, here, the collocation points from the Greville-abscissae are merely introduced to approximate the initial director field. Further their number is the same as that of the unknown control coefficient vectors, so that Eq.\,(\ref{approx_dir_init_nurbs_disc}) can be rewritten as}
\begin{align}
	\label{sys_lin_eq_init_dir}
	{{\bf{\bar N}}}\left[ {\setlength{\arraycolsep}{2.25pt}\renewcommand{\arraystretch}{1.45}\begin{array}{*{20}{c}}
			{{\bf{D}}_{\alpha 1}^{\rm{T}}}\\
			{{\bf{D}}_{\alpha 2}^{\rm{T}}}\\
			\vdots \\
			{{\bf{D}}_{\alpha n_{{\rm{cp}}}^{\rm{d}}}^{\rm{T}}}
	\end{array}} \right] = \left[ {\setlength{\arraycolsep}{2.25pt}\renewcommand{\arraystretch}{1.45}\begin{array}{*{20}{c}}
			{{{\boldsymbol{D}}^\mathrm{T}_\alpha }({{\bar \xi }_1})}\\
			{{{\boldsymbol{D}}^\mathrm{T}_\alpha }({{\bar \xi }_2})}\\
			\vdots \\
			{{{\boldsymbol{D}}^\mathrm{T}_\alpha }({{\bar \xi }_{n_{{\rm{cp}}}^{\rm{d}}}})}
	\end{array}} \right],\,\,\alpha\in\left\{1,2\right\},
\end{align}
\textcolor{black}{with square matrix}
\begin{align}
	{{\bf{\bar N}}} \coloneqq \left[ {\setlength{\arraycolsep}{2.25pt}\renewcommand{\arraystretch}{1.45}\begin{array}{*{20}{c}}
			{N_1^{{p_{\rm{d}}}}({{\bar \xi }_1})}&{N_2^{{p_{\rm{d}}}}({{\bar \xi }_1})}& \cdots &{N_{n_{{\rm{cp}}}^{\rm{d}}}^{{p_{\rm{d}}}}({{\bar \xi }_1})}\\
			{N_1^{{p_{\rm{d}}}}({{\bar \xi }_2})}&{N_2^{{p_{\rm{d}}}}({{\bar \xi }_2})}& \cdots &{N_{n_{{\rm{cp}}}^{\rm{d}}}^{{p_{\rm{d}}}}({{\bar \xi }_2})}\\
			\vdots & \vdots & \ddots & \vdots \\
			{N_1^{{p_{\rm{d}}}}({{\bar \xi }_{n_{{\rm{cp}}}^{\rm{d}}}})}&{N_2^{{p_{\rm{d}}}}({{\bar \xi }_{n_{{\rm{cp}}}^{\rm{d}}}})}& \cdots &{N_{n_{{\rm{cp}}}^{\rm{d}}}^{{p_{\rm{d}}}}({{\bar \xi }_{n_{{\rm{cp}}}^{\rm{d}}}})}
	\end{array}} \right].
\end{align}
Note that the matrix $\bf{\bar N}$ is banded with band-width $p_\mathrm{d}+1$ due to the local support property of NURBS basis functions, and invertible due to their linear independence. The control coefficients of the initial director field are determined by solving Eq.\,(\ref{sys_lin_eq_init_dir}) in \textit{pre-processing}. A continuous director field $\boldsymbol{D}_\alpha(\xi)$, which is evaluated on the right-hand side of Eq.\,(\ref{sys_lin_eq_init_dir}) at the discrete points of ${\bar\varXi}^\mathrm{d}\coloneqq\left\{{\bar \xi _1},{\bar \xi _2},\cdots,{\bar \xi _{n^\mathrm{d}_\mathrm{cp}}}\right\}$, can be obtained in various ways like with the Frenet-Serret formula and the smallest rotation method, see, e.g., \citet{meier2014objective}, \citet{choi2019isogeometric}, and references therein. In the present approach, ${\boldsymbol{D}}_\alpha ^h(\xi) \ne {{\boldsymbol{D}}_\alpha }(\xi)$ in general, except for the chosen points $\xi\in{\bar\varXi}^\mathrm{d}$. This approximation error should decrease by increasing the number of basis functions (control coefficients) or elevating the degree $p_\mathrm{d}$. One may also refer to the approach to enforce the condition ${\boldsymbol{D}}_\alpha ^h(\xi) = {{\boldsymbol{D}}_\alpha }(\xi)$ at the Gauss quadrature points of \citet{dornisch2013isogeometric}, which requires an additional matrix multiplication to obtain a square (invertible) system matrix. However, note that, in the present formulation, the beam strain components in Eqs.\,(\ref{axial_strn_comp})-(\ref{strn_high_bend}) still vanish at the undeformed configuration, in spite of the approximation error in the initial director field. 
\subsection{Spatial discretization of the variational forms}
\subsubsection{Patch-wise approximation of kinematic variables}
We approximate the current center axis position, as
\begin{equation}
	\label{cdisp_entire_patch0}
	{\boldsymbol{\varphi}}^{h}(\xi) = \sum\limits_{I = 1}^{n_\mathrm{cp}} {N_I^{{p}}(\xi )\,{{{\boldsymbol{\varphi}}}_{{{I}}}}},\quad\xi\in{\varXi},
\end{equation}
with the control coefficients ${{{{\boldsymbol{\varphi}}}}_{{{I}}}}\in\mathbb{R}^3$. Similarly, we approximate the director field, as
\begin{equation}
	\label{tot_dirv_disp_entire_patch0}
	{\boldsymbol{d}}_\alpha^{h}(\xi) = \sum\limits_{I = 1}^{n_\mathrm{cp}^{\mathrm{d}}} {N_I^{{p_\mathrm{d}}}(\xi )\,{{{\bf{d}}}_{\alpha{I}}}},\quad\alpha\in\left\{1,2\right\},\,\,\xi\in{\varXi}.
\end{equation}
In IGA, we define an \textit{element} by a non-zero knot span. Let ${\varXi}_e\coloneqq{\left[ {\xi _1^e,\xi _2^e} \right)}\ni{\xi}$ denote the $e$th element such that $\varXi={\varXi_1}\cup{\varXi_2}\cup\cdots\cup{\varXi}_{n_\mathrm{el}}$ where $n_\mathrm{el}$ denotes the total number of elements. For the last element (i.e., $e=n_\mathrm{el}$) of an open curve, we have a closed interval ${\varXi}_{n_\mathrm{el}}\coloneqq{\left[ {\xi _1^{n_\mathrm{el}},\xi _2^{n_\mathrm{el}}} \right]}$ to include the end point. \textcolor{black}{Further, we use the same mesh along the initial center axis for both the current center axis position and director vectors.} Then, in each element, from Eqs.\,(\ref{cdisp_entire_patch0}) and (\ref{tot_dirv_disp_entire_patch0}), we have
\begin{equation}
	\label{approx_caxis_pos}
	{{\boldsymbol{\varphi}}^{h}}(\xi) = \sum\limits_{I = 1}^{{n_e}} {N_I^p(\xi )\,{{\boldsymbol{\varphi}}^e_{{I}}}},\quad\xi\in{\varXi_e},
\end{equation}
and
\begin{equation}
	\label{tot_dirv_disp}
	{\boldsymbol{d}}_\alpha^{h}(\xi) = \sum\limits_{I = 1}^{n_e^\mathrm{d}} {N_I^{{p_\mathrm{d}}}(\xi )\,{{{\bf{d}}}^e_{\alpha{I}}}},\quad\xi\in{\varXi_e},\,\alpha\in\left\{1,2\right\},
\end{equation}
respectively, where $n_e=p+1$ and $n^\mathrm{d}_e=p_\mathrm{d}+1$ denote the numbers of local support basis functions, and ${{\boldsymbol{\varphi}}^e_{{I}}},{{{\bf{d}}}^e_{\alpha{I}}}\in\Bbb{R}^3$ denote the $I$th control coefficients in the $e$th element. Combining Eqs.\,(\ref{approx_caxis_pos}) and (\ref{tot_dirv_disp}), we obtain
\begin{equation}
	\label{approx_y_N}
	{{\boldsymbol{y}}^h} = \left\{ {\renewcommand{\arraystretch}{1.25}\begin{array}{*{20}{c}}
			{{{\boldsymbol{\varphi}}^{h}}}\\
			{\begin{array}{*{20}{c}}
					{{\boldsymbol{d}}_1^{h}}\\
					{{\boldsymbol{d}}_2^{h}}
			\end{array}}
	\end{array}} \right\} = {{\bf{N}}_e}(\xi )\,{{{\bf{y}}^e}},\quad\xi\in\varXi_e,
\end{equation}
with 
\begin{equation}
	{{\bf{N}}_e}\coloneqq\left[ {\begin{array}{*{20}{c}}
			{N_1^p{{\bf{1}}_{3\times3}}}& \cdots &{N_{{n_e}}^p{{\bf{1}}_{3\times3}}}&\vline& {}&{{{\bf{0}}_{3 \times 6n_e^\mathrm{d}}}}&{}\\
			{}&{{{\bf{0}}_{6 \times 3{n_e}}}}&{}&\vline& {N_1^{{p_\mathrm{d}}}{{\bf{1}}_{6\times6}}}& \cdots &{N_{n_e^\mathrm{d}}^{{p_\mathrm{d}}}{{\bf{1}}_{6\times6}}}
	\end{array}} \right],
\end{equation}
and the column array of the control coefficients in each element,
\begin{equation*} 		{{{\bf{y}}^e}}\coloneqq\left[{{{\boldsymbol{\varphi}}_{{1}}^{e\,\mathrm{T}}}},\cdots,{{{\boldsymbol{\varphi}}_{{n_e}}^{e\,\mathrm{T}}}},{{\bf{d}}_{1}^{e\,\mathrm{T}}},\cdots,{{\bf{d}}_{n^\mathrm{d}_e}^{e\,\mathrm{T}}}
	\right]^\mathrm{T},
\end{equation*}
where we have defined ${{\bf{d}}^e_{I}} \coloneqq \left[{\bf{d}}_{1I}^{e\,\mathrm{T}},{\bf{d}}_{2I}^{e\,\mathrm{T}}\right]^\mathrm{T}\in {{\mathbb{R}}^6}$. \textcolor{black}{Further, using the same basis functions as in Eq.\,(\ref{approx_y_N}), we discretize the independent velocity field, as
\begin{align}
	\label{fe_approx_phy_vel}
	{\boldsymbol{V}}^h = {{\bf{N}}_e}(\xi )\,{\bf{V}}^e,\quad\xi\in\varXi_e.
\end{align}
}
\subsubsection{Patch-wise approximation of physical stress resultants and strains}
\label{approx_fe_phy_kin_var}
It has been shown that the element-wise approximation of the physical stress resultant and strain in an isogeometric mixed finite element formulation may still suffer from numerical locking due to the parasitic terms from the higher-order inter-element continuity in the displacement field, see \citet{kikis2022two}, \citet{casquero2022removing}, \citet{choi2023selectively}, \citet{sauer2024simple}, and references therein. To circumvent this, we introduce the following patch-wise approximation of the physical variables.
\begin{subequations}
	\label{patch_approx_phy_strs_strn}
	\begin{align}
		{\boldsymbol{r}}_{\rm{p}}^h(\xi ) &= {{\bf{B}}_e}(\xi )\,{\bf{r}}^e, \\
		{\boldsymbol{\varepsilon }}_{\rm{p}}^h(\xi ) &= {{\bf{B}}_e}(\xi )\,{\bf{e}}^e,
	\end{align}
\end{subequations}
$\xi\in{\varXi}_e$, with
\begin{align}
	{{\bf{B}}_e}(\xi ) \coloneqq \left[ {\begin{array}{*{20}{c}}
			{B_1^{{p_{\rm{p}}}}(\xi ){{\bf{1}}_{{d_\mathrm{p}} \times {d_\mathrm{p}}}}}& \cdots &{B_{n_e^{\rm{p}}}^{{p_{\rm{p}}}}(\xi ){{\bf{1}}_{{{d_\mathrm{p}} \times {d_\mathrm{p}}}}}}
	\end{array}} \right],
\end{align}
where $B^{p_\mathrm{p}}_I(\xi)$ denotes the $I$th B-spline basis functions of degree $p_\mathrm{p}$, and $n^\mathrm{p}_e=p_\mathrm{p}+1$ denotes the number of basis functions per element. Here, we choose $p_\mathrm{p}=p-1$ for numerical stability. 
\vspace{2mm}
\begin{remark} 
\label{rem_patch_nofs_compare}
\textcolor{black}{Even though this \textit{global} approach turns out to be numerically more stable compared to the \textit{local} (element-wise) approach with further reduced degree $p_\mathrm{p}$ of \citet{choi2023selectively}, the condensation of those control coefficients of the physical solution fields is computationally prohibitive. Thus, here, we include those control coefficients in the global system of equations. However, compared to the local approaches, for the same number of elements, this approach still requires a much smaller number of control coefficients, due to the higher-order continuity of the B-spline functions $B^{p_\mathrm{p}}_I$.}
\end{remark}
\vspace{2mm}
\noindent\textcolor{black}{In Section \ref{num_ex_twist_straight}, we present a numerical example to compare DOF numbers between the global (`mix.glo') and local (`mix.loc-ur', `mix.loc-sr') approaches.} 
\subsubsection{Element-wise approximation of enhanced strains}
\label{approx_fe_phy_enh_var}
It should be noted that, in contrast to the physical strain ($\boldsymbol{\varepsilon}_\mathrm{p}$), the enhanced strain ($\boldsymbol{\alpha}$) has no associated geometric strain from the configuration variable. This may simplify the selection of approximate solution space for the enhanced strain field, without any special care about numerical locking. In order to facilitate an element-wise elimination (condensation) of the control coefficients, we allow inter-element discontinuity. We utilize Lagrange polynomial basis functions in the parametric domain $\left[-1,1\right)\ni{\bar \xi}$ which can be mapped to the interval ${{\varXi}}_e\coloneqq\left[\xi^e_1,\xi^e_2\right)\ni{\xi}$ via 
\begin{equation}
	{\bar \xi}=1-2\left(\frac{\xi^e_2-\xi}{\xi^e_2-\xi^e_1}\right).
\end{equation}
Then, we approximate the enhanced strain as
\begin{equation}
	\label{disc_alph_L}
	{{\boldsymbol{\alpha }}^h}\!\left({\bar \xi}\,\right) = {{\bf{L}}_e}\!\left({\bar \xi}\,\right){{\boldsymbol{\alpha}}^e},
\end{equation}
with
\begin{equation}
	{{\bf{L}}_e}\!\left({\bar \xi}\,\right) \coloneqq \left[ {\begin{array}{*{20}{c}}
			{L_1^{{p_{\rm{a}}}}\!\left({\bar \xi}\,\right)\!{{\bf{1}}_{{d_{\rm{a}}}\times{d_{\rm{a}}}}}}& \cdots &{L_{n_e^{\rm{a}}}^{{p_{\rm{a}}}}\!\left({\bar \xi}\,\right)\!{{\bf{1}}_{{d_{\rm{a}}}\times{d_{\rm{a}}}}}}
	\end{array}} \right],
\end{equation}
and the nodal coefficients $\boldsymbol{\alpha}^e$ in the $e$th element. Here, $L^{p_\mathrm{a}}_I$ denotes the $I$th basis function of degree $p_\mathrm{a}$, and $n^e_\mathrm{a}=p_\mathrm{a}+1$ denotes the number of basis functions per element. Here, without any special care about numerical instability, we simply choose the lowest degree of the basis functions, $p_\mathrm{a}=0$.
\subsection{Discretization of linearized variational equation}
We first discretize the external virtual work in Eq.\,(\ref{lin_var_eq_emc}), using Eq.\,(\ref{approx_y_N}), as
\begin{equation}
	{G_{{\rm{ext}}}}\left( {\delta {\boldsymbol{y}},{t_n} + \frac{1}{2}\Delta t} \right) \approx \delta {{\bf{y}}^{\rm{T}}}{{\bf{F}}_{{\rm{ext}}}}\left(t_n+\frac{1}{2}{\Delta t}\right),
\end{equation}
where we have defined the global external load vector
\begin{equation}
	{{\bf{F}}_{{\rm{ext}}}}(t) \coloneqq \mathop {\mathlarger{\mathlarger{{\bf{A}}}}}\limits_{e = 1}^{{n_{{\rm{el}}}}} {\bf{f}}_{{\rm{ext}}}^e(t) + {\mathlarger{\mathlarger{{\bf{A}}}}}\,{\bigl[ {{{{\boldsymbol{\bar R}}}_0}}(t) \bigr]_{s \in {\Gamma _{\rm{N}}}}},
\end{equation}
from the element external load vector due to the distributed load,
\begin{align}
	{\bf{f}}_{{\rm{ext}}}^e(t) \coloneqq \int_{{\varXi _e}} {{\bf{N}}_e^{\rm{T}}{\boldsymbol{\bar R}}(t)\,{\tilde j}\,{\rm{d}}\xi }.
\end{align}
Here, $\mathop{\mathlarger{\mathlarger{\mathlarger{{\bf{A}}}}}}$ denotes the finite element assembly operator, and ${\bf{y}}$ denotes the global array of the control coefficients for the approximated configuration variable $\boldsymbol{y}^h$. Then, the linearized equation\,(\ref{lin_var_eq_emc}) is discretized as
\begin{subequations}
	\begin{align}
		&\sum\limits_{e = 1}^{{n_{{\rm{el}}}}} \delta {{\bf{y}}^{e{\rm{T}}}}\left\{{\left({{\bf{k}}_{{\rm{yy}}}^e}+\frac{4}{\Delta t^2}{\bf{m}}^e\right)\Delta {{\bf{y}}^e} + {\bf{k}}_{{\rm{yr}}}^e\Delta {{\bf{r}}^e}} \right\}= \nonumber\\
		&\delta {{\bf{y}}^{\rm{T}}}2\,{{\bf{F}}_{{\rm{ext}}}}\!\left( {{t_n} \!+ \!\frac{1}{2}\Delta t} \right)\!-\!\sum\limits_{e = 1}^{{n_{{\rm{el}}}}} {\delta {{\bf{y}}^{e{\rm{T}}}}2\!\!\:\left({ {\bf{f}}_{\rm{y}}^e}+\frac{2}{\Delta t}{\bf{f}}^e_\rho\right)},\!\label{disc_var_eq_emc_y_sig}
	\end{align}
	and
	\begin{align}
		&\sum\limits_{e = 1}^{{n_{{\rm{el}}}}} \delta {{\bf{r}}^{e{\rm{T}}}}\left({\bf{k}}_{{\rm{ry}}}^e\Delta {{\bf{y}}^e} + {\bf{k}}_{\mathrm{r}\varepsilon }^e\Delta {{\bf{e}}^e}\right) =  -\sum\limits_{e = 1}^{{n_{{\rm{el}}}}} \delta {{\bf{r}}^{e{\rm{T}}}}{\bf{f}}_{\rm{r}}^e,\label{elem_eq_disc_r}\\
		&\sum\limits_{e = 1}^{{n_{{\rm{el}}}}}\delta {{\bf{e}}^{e{\rm{T}}}}\left({\bf{k}}_{\mathrm{r}\varepsilon }^{e\,{\rm{T}}}\Delta {{\bf{r}}^e} + {\bf{k}}_{\varepsilon \varepsilon }^e\Delta {{\bf{e}}^e} + {\bf{k}}_{{\rm{a}}\varepsilon}^{e\,\mathrm{T}}\Delta {{\boldsymbol{\alpha}}^e}\right) =  - 2\,\sum\limits_{e = 1}^{{n_{{\rm{el}}}}}\delta {{\bf{e}}^{e{\rm{T}}}}{\bf{f}}_\varepsilon ^e,\label{elem_eq_disc_e}\\
		&{\bf{k}}_{{\rm{a}}\varepsilon}^{e}\Delta {{\bf{e}}^e} + {\bf{k}}_{{\rm{aa}}}^e\Delta {{\boldsymbol{\alpha }}^e} =  - 2\,{\bf{f}}_{\rm{a}}^e,\,\,\mathrm{for}\,\,e \in \left\{ {1,\cdots,{n_{{\rm{el}}}}} \right\},\label{elem_eq_disc_a}
	\end{align}
\end{subequations}
where we have defined the following elemental matrices for the inertial contribution 
\begin{align}
	{{\bf{m}}^e} \coloneqq \int_{{\varXi _e}} {{\bf{N}}_e^{\rm{T}}\,{\boldsymbol{\mathcal{M}}}\,{{\bf{N}}_e}\,\tilde j\,{\rm{d}}\xi }\,,
\end{align}
\begin{align}
	{{\bf{f}}^e_\rho}&\coloneqq\int_{{\varXi _e}}{{{\bf{N}}_e^{{\rm{T}}}}{\boldsymbol{\mathcal{M}}}\left\{{\frac{1}{{\Delta t}}\!\left( {{}^{n + 1}{{\boldsymbol{y}}^{h\,(i - 1)}} - {}^n{\boldsymbol{y}}^h} \right) - {}^n{\boldsymbol{V}}^h} \right\}\tilde j\,{\rm{d}}\xi }\,.
\end{align}
We also have
\begin{subequations}
	\begin{align}
		{\bf{k}}_{{\rm{yy}}}^e &\coloneqq \int_{{\varXi _e}} {\mathbb{Y}}_e^{\rm{T}}{{{\boldsymbol{k}}}_{\rm{G}}({\bar{\boldsymbol{r}}}^h_\mathrm{p})}{{\mathbb{Y}}_e}\,\tilde j\,{\rm{d}}\xi,\\
		{\bf{k}}_{{\rm{yr}}}^e &\coloneqq \int_{{\varXi _e}} {{\mathbb{B}}_e^{{\rm{T}}}\!\left({\boldsymbol{\bar y}}^h\right){{\bf{B}}_e}\,\tilde j\,{\rm{d}}\xi },\\
		{\bf{k}}_{{\rm{ry}}}^e &\coloneqq \int_{{\varXi _e}} {{\bf{B}}_e^{\rm{T}}{\mathbb{B}}_e\!\left( {{{}^{n+1}{\boldsymbol{y}}^{h(i - 1)}}} \right)\tilde j\,{\rm{d}}\xi },\\
		{\bf{k}}_{{\rm{r}}\varepsilon }^e &\coloneqq  - \int_{{\varXi _e}} {{\bf{B}}_e^{\rm{T}}{\bf{B}}_e^{}\,\tilde j\,{\rm{d}}\xi },\\
		{\bf{k}}_{\varepsilon \varepsilon }^{e} &\coloneqq \int_{{\varXi _e}} {{\bf{B}}_e^{\rm{T}}\,{{{\mathbb{C}}}_\mathrm{p}^{\varepsilon \varepsilon }}\!\left( {\boldsymbol{\bar \varepsilon}}_{\rm{p}}^{\,h},{{\boldsymbol{\bar \alpha }}^{h}}\right){\bf{B}}_e\,\tilde j\,{\rm{d}}\xi },\\
		{\bf{k}}_{{\rm{a}}\varepsilon}^e &\coloneqq \int_{{\varXi _e}} {{\bf{L}}_e\!\!\!\!\;^{\rm{T}}\,{{{\mathbb{C}}}_\mathrm{p}^{{\rm{a}}\varepsilon}}\!\left( {\boldsymbol{\bar\varepsilon }}_{\rm{p}}^{\,h},{{\boldsymbol{\bar\alpha }}^{h}}\right){\bf{B}}_e\,\tilde j\,{\rm{d}}\xi },\\
		{\bf{k}}_{{\rm{aa}}}^e &\coloneqq \int_{{\varXi _e}} {{\bf{L}}_e\!\!\!\!\;^{\rm{T}}\,{{{\mathbb{C}}}_\mathrm{p}^{{\rm{aa}}}}\!\left( {\boldsymbol{\bar\varepsilon }}_{\rm{p}}^{\,h},{{\boldsymbol{\bar\alpha }}^{h}}\right){\bf{L}}_e^{}\,\tilde j\,{\rm{d}}\xi },
	\end{align}
\end{subequations}
and
\begin{subequations}
	\begin{align}
		{{\bf{f}}^e_\mathrm{y}}&\!\coloneqq\!\int_{{\varXi _e}} {{\mathbb{B}}_{e}^\mathrm{T}\!\left({\boldsymbol{\bar y}}^h\right){\boldsymbol{\bar r}}_\mathrm{p}^{h}\,\tilde j\,{\rm{d}}\xi},\\
		{\bf{f}}_{\rm{r}}^e\!&\!\coloneqq\!\int_{{\varXi _e}}\!{{\bf{B}}_e^{\rm{T}}\!\left\{ {\boldsymbol{\varepsilon}\!\left( {{}^{n + 1}{{\boldsymbol{y}}^{h\,(i - 1)}}} \right)\!-\!{}^{n + 1}{\boldsymbol{\varepsilon }}_{\rm{p}}^{h\,(i - 1)}}\!\right\}\tilde j\,{\rm{d}}\xi},\label{approx_f_r_e_resid}\!\\
		{\bf{f}}_\varepsilon ^e &\!\coloneqq\!\int_{{\varXi _e}} {{\bf{B}}_e^{\rm{T}}\left\{{{{\partial_{{{\boldsymbol{\varepsilon}_\mathrm{p}}}}}\psi\!\left({\boldsymbol{\bar\varepsilon }}_{\rm{p}}^{\,h},{{\boldsymbol{\bar\alpha }}^{h}} \right)} - {\boldsymbol{\bar r}}_{\rm{p}}^{h}}\right\}\tilde j\,{\rm{d}}\xi},\\
		{\bf{f}}_{\rm{a}}^e &\!\coloneqq\!\int_{{\varXi _e}} {{\bf{L}}_e\!\!\!\!\;^{\rm{T}}\,{\partial_{{{\boldsymbol{{\alpha}}}}}}\psi\!\left({\boldsymbol{\bar \varepsilon }}_{\rm{p}}^{h},{{\boldsymbol{\bar \alpha }}^{h}}\right)\,\tilde j\,{\rm{d}}\xi},
	\end{align}
\end{subequations}
\textcolor{black}{where the matrices ${\Bbb{B}}_e$ and $\Bbb{Y}_e$ are defined by Eqs.\,(A.5.2) and (A.5.4) in \citet{choi2021isogeometric}, respectively.} Here, from Eqs.\,(\ref{mid_pt_def_ybar})-(\ref{mid_pt_def_alph}), we also define
\begin{subequations}
	\begin{align*}
		{\boldsymbol{\bar y}}^h&\coloneqq\frac{1}{2}\left({{}^{n + 1}{{\boldsymbol{y}}^{h\,(i - 1)}} + {}^n{\boldsymbol{y}}^h}\right),\\
		{\boldsymbol{\bar r}}^h_\mathrm{p}&\coloneqq{\frac{1}{2}\left( {{}^{n + 1}{\boldsymbol{r}}_{\rm{p}}^{h\,(i - 1)} + {}^n{{\boldsymbol{r}}^h_{\rm{p}}}} \right)},\\	
		{\boldsymbol{\bar \varepsilon}}^h_\mathrm{p}&\coloneqq{\frac{1}{2}\left( {{}^{n + 1}{\boldsymbol{\varepsilon}}_{\rm{p}}^{h\,(i - 1)} + {}^n{{\boldsymbol{\varepsilon}}^h_{\rm{p}}}} \right)},\\		
		{\boldsymbol{\bar \alpha}}^h&\coloneqq{\frac{1}{2}\left( {{}^{n + 1}{\boldsymbol{\alpha}}^{h\,(i - 1)} + {}^n{{\boldsymbol{\alpha}}}^h} \right)}.	
	\end{align*}
\end{subequations}
Note that Eq.\,(\ref{elem_eq_disc_a}) is given element-wisely, since we allow discontinuity in the enhanced strain field $\boldsymbol{\alpha}^h$. In contrast, Eqs.\,(\ref{disc_var_eq_emc_y_sig})--(\ref{elem_eq_disc_e}) are given by a global equation, due to the higher-order continuity in the solution fields $\boldsymbol{y}^h$, $\boldsymbol{r}_\mathrm{p}^h$, and $\boldsymbol{\varepsilon}_\mathrm{p}^h$.
\subsubsection{Condensation and solution update}
\label{sec_condens_sol_update}
The unknown nodal variables of the enhanced strain field ${\boldsymbol{\alpha}}^e$ can be eliminated from the global finite element equation, via element-wise condensation. As the matrix ${\bf{k}}^e_\mathrm{aa}$ is invertible, from Eq.\,(\ref{elem_eq_disc_a}), we obtain
\begin{subequations}
	\begin{align}
		\Delta {{\boldsymbol{\alpha }}^e} =  - {\bf{k}}_{{\rm{aa}}}^{e - 1}\!\left( {2\,{\bf{f}}_{\rm{a}}^e + {\bf{k}}_{{\rm{a}}\varepsilon}^{e}\Delta {{\bf{e}}^e}} \right),\,\,\mathrm{for}\,\,e \in \left\{ {1,\cdots,{n_{{\rm{el}}}}} \right\}.\label{emc_cond_inc_a}
	\end{align}
\end{subequations}
Substituting Eq.\,(\ref{emc_cond_inc_a}) into Eq.\,(\ref{elem_eq_disc_e}), we obtain
\begin{align}
	\label{inc_rew_del_e}
	\sum\limits_{e = 1}^{{n_{{\rm{el}}}}}\delta {{\bf{e}}^{e{\rm{T}}}}\left({\bf{k}}_{\mathrm{r}\varepsilon }^{e\,{\rm{T}}}\Delta {{\bf{r}}^e} + {\mathbbm{k}}_{\varepsilon \varepsilon }^e\Delta {{\bf{e}}^e} \right) =  - 2\,\sum\limits_{e = 1}^{{n_{{\rm{el}}}}}\delta {{\bf{e}}^{e{\rm{T}}}}{\mathbbm{f}}_\varepsilon ^e,
\end{align}
with
\begin{align}
	{\mathbbm{f}}_\varepsilon ^e &\coloneqq {\bf{f}}_\varepsilon ^e - {\bf{k}}_{{\rm{a}}\varepsilon}^{e\,\mathrm{T}}{\bf{k}}_{{\rm{aa}}}^{e\, - 1}{\bf{f}}_{\rm{a}}^e,\\
	{{\Bbbk}}_{\varepsilon \varepsilon }^e &\coloneqq {\bf{k}}_{\varepsilon \varepsilon }^e - {\bf{k}}_{{\rm{a}}\varepsilon}^{e\,\mathrm{T}}{\bf{k}}_{{\rm{aa}}}^{e\, - 1}{\bf{k}}_{{\rm{a}}\varepsilon }^e.
\end{align}
Combining Eqs.\,(\ref{disc_var_eq_emc_y_sig}), (\ref{elem_eq_disc_r}), and (\ref{inc_rew_del_e}), we obtain
\begin{align}
	\label{fe_elem_eq_cond_alpha}
	\sum\limits_{e = 1}^{{n_{{\rm{el}}}}} {{{\left\{ {\begin{array}{*{20}{c}}
						{\delta {{\bf{y}}^e}}\\
						{\delta {{\bf{r}}^e}}\\
						{\delta {{\bf{e}}^e}}
				\end{array}} \right\}}^{\!\rm{T}}}} {{\mathbbm{k}}^e}\left\{ {\begin{array}{*{20}{c}}
			{\Delta {{\bf{y}}^e}}\\
			{\Delta {{\bf{r}}^e}}\\
			{\Delta {{\bf{e}}^e}}
	\end{array}} \right\} = \delta {{\bf{y}}^{\rm{T}}}\left\{2\,{{\bf{F}}_{{\rm{ext}}}}\left({t_n} + \frac{1}{2}\Delta t\right)\right\} - \sum\limits_{e = 1}^{{n_{{\rm{el}}}}} {{{\left\{ {\begin{array}{*{20}{c}}
						{\delta {{\bf{y}}^e}}\\
						{\delta {{\bf{r}}^e}}\\
						{\delta {{\bf{e}}^e}}
				\end{array}} \right\}}^{\!\rm{T}}}(2\,{{\mathbbm{f}}^e})},
\end{align}
with
\begin{align}
	\label{elem_t_stiff_kke}
	{{\mathbbm{k}}^e} \coloneqq \left[ {\renewcommand{\arraystretch}{2}\begin{array}{*{20}{c}}
			\left({{\bf{k}}_{{\rm{yy}}}^e + \dfrac{4}{{\Delta {t^2}}}{{\bf{m}}^e}}\right)&{{\bf{k}}_{{\rm{yr}}}^e}&{\bf{0}}_{{m_e}\times{m^\mathrm{p}_e}}\\
			{{\bf{k}}_{{\rm{ry}}}^e}&{\bf{0}}_{{m^\mathrm{p}_e}\times{m^\mathrm{p}_e}}&{{\bf{k}}_{{\rm{r}}\varepsilon }^e}\\
			{\bf{0}}_{{m^\mathrm{p}_e}\times{m_e}}&{{\bf{k}}_{{\rm{r}}\varepsilon }^{e\,{\rm{T}}}}&{{\mathbbm{k}}_{\varepsilon \varepsilon }^e}
	\end{array}} \right]
\end{align}
and
\begin{align}
	{{\mathbbm{f}}^e} \coloneqq \left\{ {\renewcommand{\arraystretch}{2}\begin{array}{*{20}{c}}
			{{{\bf{f}}_{\rm{y}}^e + \dfrac{2}{{\Delta t}}{\bf{f}}_\rho ^e}}\\
			\dfrac{1}{2}{{\bf{f}}_{\rm{r}}^e}\\
			{{\mathbbm{f}}_\varepsilon ^e}
	\end{array}} \right\},
\end{align}
where we have the DOF numbers $m_e=3{n_e}+6{n^\mathrm{d}_e}$, and $m^\mathrm{p}_e=d_\mathrm{p}\cdot{n^\mathrm{p}_e}$. Then, Eq.\,(\ref{fe_elem_eq_cond_alpha}) can be rewritten as
\begin{equation}
	\label{Kd_R_glob_form}
	\delta {{\bf{z}}^{\rm{T}}}{{\bf{K}}}\Delta {\bf{z}} = \delta {{\bf{z}}^{\rm{T}}}2\,{{\bf{R}}},
\end{equation}
where we define the \textit{effective} global tangent stiffness matrix 
\begin{equation}
	{{\bf{K}}} \coloneqq \mathop {\mathlarger{\mathlarger{{\bf{A}}}}}\limits_{e = 1}^{{n_{{\rm{el}}}}} {{\mathbbm{k}}^e},
\end{equation}
and the effective global residual vector
\begin{align}
	\label{fe_disc_resid_vec_condens}	
	{\bf{R}} \coloneqq \left\{ {\renewcommand{\arraystretch}{2}\begin{array}{*{20}{c}}
			{{{\bf{F}}_{{\rm{ext}}}}\!\left( {{t_n} + \dfrac{1}{2}\Delta t} \right)}\\
			{{{\bf{0}}_{2n^\mathrm{p}_\mathrm{dof} \times 1}}}
	\end{array}} \right\} - \mathop {\mathlarger{\mathlarger{{\bf{A}}}}}\limits_{e = 1}^{{n_{{\rm{el}}}}} {{\mathbbm{f}}^e}.
\end{align}
Here, $\bf{z}\coloneqq{[{\bf{y}}^\mathrm{T},{\bf{r}}^\mathrm{T},{\bf{e}}^\mathrm{T}]}^\mathrm{T}$ denotes a column array of the unknown global control coefficients, which assembles the control coefficients for the configuration variable (${\bf{y}}$), the physical stress resultant (${\bf{r}}$), and the physical strain (${\bf{e}}$), where $\bf{r}$ and $\bf{e}$ have the same number of entries, $n^\mathrm{p}_\mathrm{dof}$.
It should be noted that ${{\bf{k}}^{e\,\mathrm{T}}_\mathrm{yr}}\ne{{\bf{k}}^e_\mathrm{ry}}$ in Eq.\,(\ref{elem_t_stiff_kke}) \textcolor{black}{(see Remark\,\ref{rem_unsym_b_btr})}, which makes ${\mathbbm{k}}^e$, and eventually $\bf{K}$ unsymmetric. After applying the displacement boundary conditions to Eq.\,(\ref{Kd_R_glob_form}), we finally have the reduced system of linear equations at $i$th iteration in the $(n+1)$th time step,
\begin{equation}
	\label{Kd_R_glob_form_red}
	{{\bf{K}}}_\mathrm{r}\,\Delta {\bf{z}}_\mathrm{r} = 2\,{{\bf{R}}}_\mathrm{r},
\end{equation}
where $(\bullet)_\mathrm{r}$ denotes the reduced matrix. 
Using $\Delta {\bf{z}}$ obtained from solving Eq.\,(\ref{Kd_R_glob_form_red}), we update the global control coefficients at the $i$th iteration in the $(n+1)$th time step, as 
\begin{subequations}
	\begin{alignat}{5}
		{}^{n+1}{{\bf{y}}^{(i)}} &= {}^{n+1}{{\bf{y}}^{(i-1)}} &&+ \Delta {\bf{y}},\,\,{}^{n+1}{{\bf{y}}^{(0)}}&&\equiv{}^{n}{\bf{y}},\\
		{}^{n+1}{{\bf{r}}^{(i)}} &= {}^{n+1}{{\bf{r}}^{(i-1)}} &&+ \Delta {\bf{r}},\,\,{}^{n+1}{{\bf{r}}^{(0)}}&&\equiv{}^{n}{\bf{r}},\\	%
		{}^{n+1}{{\bf{e}}^{(i)}} &= {}^{n+1}{{\bf{e}}^{(i-1)}} &&+ \Delta {\bf{e}},\,\,{}^{n+1}{{\bf{e}}^{(0)}}&&\equiv{}^{n}{\bf{e}}.
	\end{alignat}
\end{subequations}
Further, using the element-wise control coefficients $\Delta {\bf{e}}^e$ extracted from the global one $\Delta {\bf{e}}$, we calculate the increment $\Delta{\boldsymbol{\alpha}}^e$ using Eq.\,(\ref{emc_cond_inc_a}), which is followed by the update,
\begin{align}
	{}^{n+1}{{\boldsymbol{\alpha}}^{e\,(i)}} = {}^{n+1}{{\boldsymbol{\alpha}}^{e\,(i - 1)}} + \Delta {\boldsymbol{\alpha}}^e,\,{}^{n+1}{{\boldsymbol{\alpha}}^{e\,(0)}}\equiv{}^{n}{{\boldsymbol{\alpha}}}^e,\,\,\mathrm{for}\,\,e\in\left\{1,2,\cdots,n_\mathrm{el}\right\}.
\end{align}
\textcolor{black}{After the iteration finishes, we need to update the array of global control coefficients ${\bf{V}}$ for the approximated independent velocity field in Eq.\,(\ref{fe_approx_phy_vel}), using the following update formula from Eq.\,(\ref{rep_cur_V_conf}).
\begin{align}
	{}^{n + 1}{{\bf{V}}} = \frac{2}{{\Delta t}}\left( {{}^{n + 1}{\bf{y}} - {}^n{\bf{y}}} \right) - {}^n{{\bf{V}}}.
\end{align}
}
\section{Numerical examples}
\label{num_ex}
We present six numerical examples whose objectives are as follows.
\begin{itemize}
	\item \textbf{Ex.\,1.\,Rigid rotation of a stress-free rod (quasi-statics).} We verify the invariance of the beam strain under a rigid rotation, for any degree $p_\mathrm{d}$ by employing the present discretization of the initial director field.\\ 
	\item \textbf{Ex.\,2.\,Rigid rotation of a bent rod (quasi-statics).} We investigate the invariance of the beam strain under a rigid rotation superposed to a deformed configuration.\\
	\item \textbf{Ex.\,3.\,A straight beam under twisting moment (quasi-statics).} We verify the present EAS method to correct the torsional stiffness by comparison with analytical and brick element solutions.\\
	\item \textbf{Ex.\,4.\,Twisting of an elastic ring (quasi-statics).} We verify the present EAS method to correct stiffness for an initially curved beam undergoing a large rotational motion.\\
	\item \textbf{Ex.\,5.\,Flying beam (dynamics).} We verify the superior numerical stability in time-stepping due to the conservation of the total energy from using the EMC scheme, compared with the results from using standard schemes.\\
	\item \textbf{Ex.\,6.\,Slotted ring with a rectangular cross-section (quasi-statics and dynamics).} We verify the efficiency and accuracy of the present beam solution by comparison with a brick element one. Further, we verify the energy--consistency of the present EMC scheme.\\
\end{itemize}
We consider two different ways of evaluating the initial director field:
\begin{itemize}
	\item {[{$D$--cont.}]} A continuous (exact) representation of the initial director field.\\
	\item {[$D$--disc.]} A discrete (reconstructed) initial director field from the continuous one evaluated at chosen discrete points, see Section \ref{approx_init_dir_sec}.
\end{itemize}
Further, we employ the following different finite element formulations for reference solutions:
\begin{itemize}
	\item {[disp.]} The displacement-based isogeometric beam formulation with extensible directors in \citet{choi2021isogeometric}.\\
	\item {[$\Delta\theta$]} The finite element formulation of geometrically exact beams from \citet{simo1986three}, with \textcolor{black}{Lagrange polynomial basis functions and} uniformly reduced integration.\\
	\item {[Brick]} An isogeometric brick element formulation, where we denote the degree of basis functions by `deg.$=(p_\mathrm{L},p_\mathrm{W},p_\mathrm{H})$', and the number of elements by `$n_\mathrm{el}={n^\mathrm{L}_\mathrm{el}}\times{n^\mathrm{W}_\mathrm{el}}\times{n^\mathrm{H}_\mathrm{el}}$', where `L', `W', and `H' represent the directions along the length, width, and height, respectively.
\end{itemize}
The present isogeometric mixed finite element formulation can be combined with the following approaches:
\begin{itemize}
	\item {[{mix.loc-ur}]} $p_\mathrm{d}=p-1$, and a uniformly reduced degree $p_\mathrm{p} = 1$.\\	
	\item {[{mix.loc-sr}]} $p_\mathrm{d}=p-1$, and a selectively reduced degree $p_\mathrm{p}$ from Table A.1 in \citet{choi2023selectively}.\\	
	\item {[{mix.glo}]} $p_\mathrm{d}=p$, $p_\mathrm{p}=p-1$, and the inter-element continuity of the physical stress resultant and strain is $C^{p_\mathrm{p}-1}$, see Section \ref{approx_fe_phy_kin_var}.\\	
\end{itemize}
In both [{mix.loc-ur}] and [{mix.loc-sr}], the reduction of the degrees $p_\mathrm{p}$ and $p_\mathrm{d}=p-1$ compared to [mix.glo] is due to the parasitic strains arising from the higher-order inter-element continuity of the displacement field in IGA. \citet{choi2023selectively} showed that a further (selective) reduction of $p_\mathrm{p}$ in [mix.loc-sr] can alleviate the parasitic strains more effectively in general cases than the uniform degree $p_\mathrm{p}=1$ in [mix.loc-ur] does. \textcolor{black}{Although these local approaches allow element-wise condensation, we have observed numerical instability in some cases. In contrast, the approach [mix.glo] shows numerical stability in those cases. Further, due to the global patch-wise approximation, the number of DOFs for the physical stress resultants and strains in [mix.glo] is much lower than those in the local approaches, [mix.loc-ur] and [mix.loc-sr].} In all of those three approaches, we apply element-wise condensation of the enhanced strain field $\boldsymbol{\alpha}^h$. \textcolor{black}{In the present mixed formulation, we use full Gauss integration for both longitudinal and transverse directions. In the center axis, we use $p+1$ quadrature points per element, and in the cross-section, we use $(m_\mathrm{max}+1)\times(m_\mathrm{max}+1)$ quadrature points, where $m_\mathrm{max}$ denotes the largest value in $\left\{m_1,m_2,m_3,m_4\right\}$.}
\subsection{Rigid rotation of a stress-free rod: Objectivity test 1}
\label{sec_ex_obj_test_free}
\begin{figure}[h]
	\centering
	\includegraphics[width=0.75\linewidth]{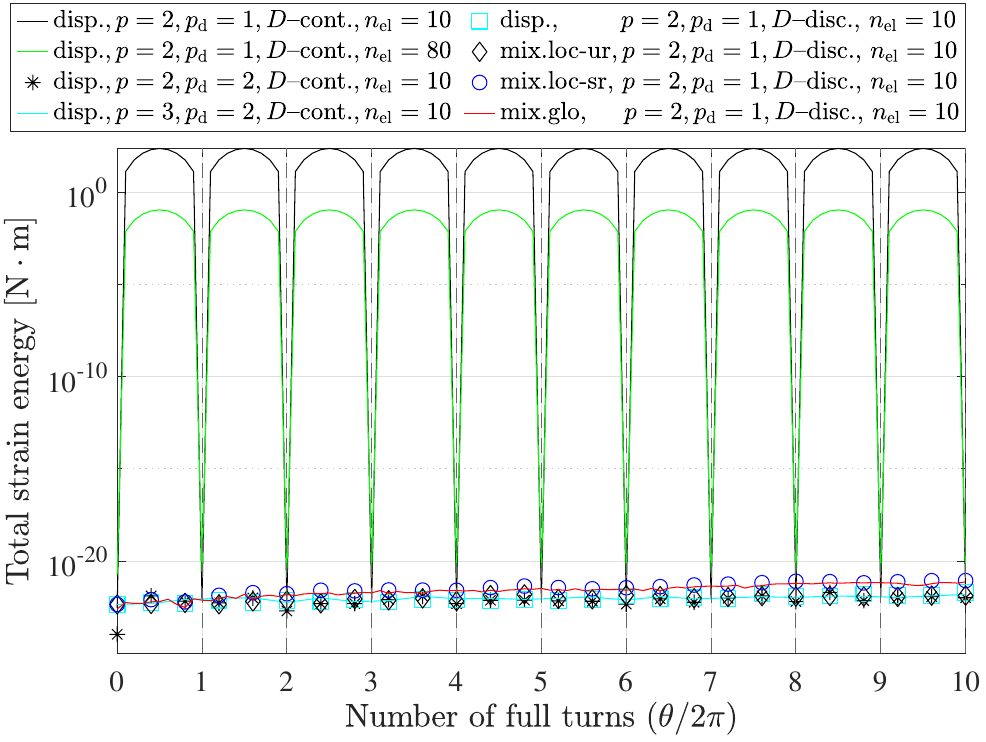}
	\caption{{Rigid rotation of a stress-free rod (objectivity test\,1): Change of the total strain energy under a rigid rotation. The dashed vertical lines indicate full turns. Note that those markers are only plotted in every four rotational increments for better visibility.}}
	\label{obj_test_free_total_strain_e}	
\end{figure}

We consider \textcolor{black}{a quasi-static beam problem where the initial center axis} is placed on the $XY$-plane and represents a quarter circle of radius $R=100\,\mathrm{m}$. It has a square cross-section of dimension $d=1\,\mathrm{m}$. We constrain one end of the beam, at which a rotation of angle $20\pi\,[\mathrm{rad}]$ (10 full turns) with respect to the $X$-axis is prescribed. We prescribe the rotational motion by the two directors at the constrained end with their extensions not allowed. We incrementally apply the total rotation in 100 uniform increments ($n_\mathrm{load}=100$), and investigate the change of total strain energy. A St.\,Venant-Kirchhoff type isotropic hyperelastic material is considered with the Young's modulus $E=1\,\mathrm{MPa}$ and the Poisson's ratio $\nu=0$. This example for testing the objectivity was proposed in \citet[Sec. 11.1]{meier2019geometrically}.
\noindent In Fig.\,\ref{obj_test_free_total_strain_e}, in cases of using the displacement-based beam formulation {[{disp.}]}, and the continuous initial director field {[{$D$--cont.}]} with $p_\mathrm{d}=1$, it is seen that the strain energy does not vanish (black and green curves), which represents the inability to represent the rigid body rotation (= non-objectivity). It is noticeable that the spurious strain energy vanishes in every full turn (= path-independence). It is also seen that the error decreases by increasing the number of elements (green curve). If we consider the continuous (exact) initial director field, we need at least quadratic NURBS basis functions in order to represent the rigid body rotation exactly, due to the circular initial geometry \citep{choi2023selectively}. This can be seen in cases of using $p_\mathrm{d}=2$ and 3 (black star and cyan curve) where the strain energy vanishes up to machine precision. In contrast, employing the discrete initial director field {{[{$D$--disc.}]}} gives objectivity even with $p_\mathrm{d}=1$ (cyan square, black diamond, blue circle, and red curve). 
\subsection{Rigid rotation of a bent rod: Objectivity test 2}
\label{num_ex_obj_test2}
\begin{figure}[H]
	\centering
	\begin{subfigure}[b]{0.4875\textwidth}\centering
		\includegraphics[width=\linewidth]{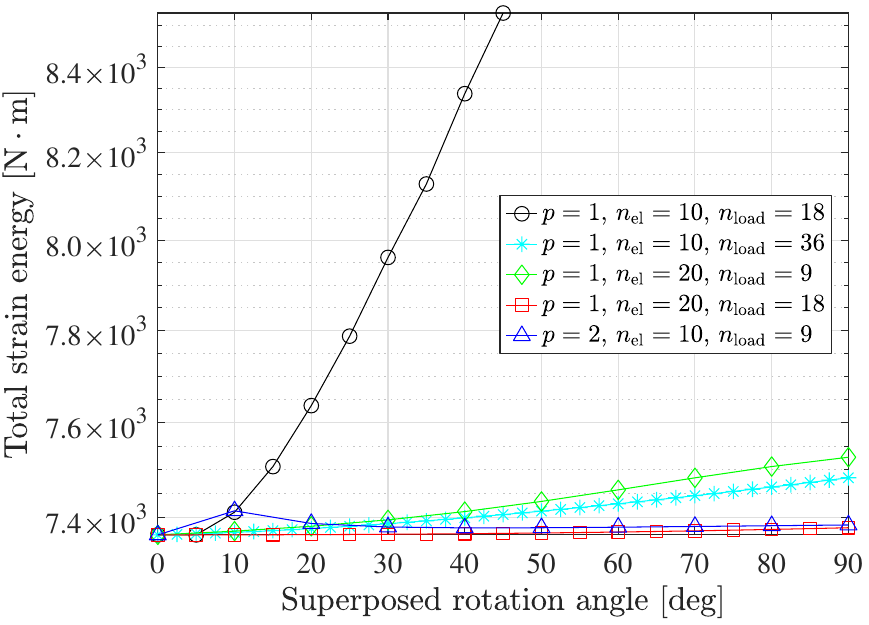}
		\caption{Beam formulation [$\Delta\theta$]}
		\label{obj_test_bent_total_strain_e_ext_dir_sv}			
	\end{subfigure}
	\begin{subfigure}[b]{0.4875\textwidth}\centering
		\includegraphics[width=\linewidth]{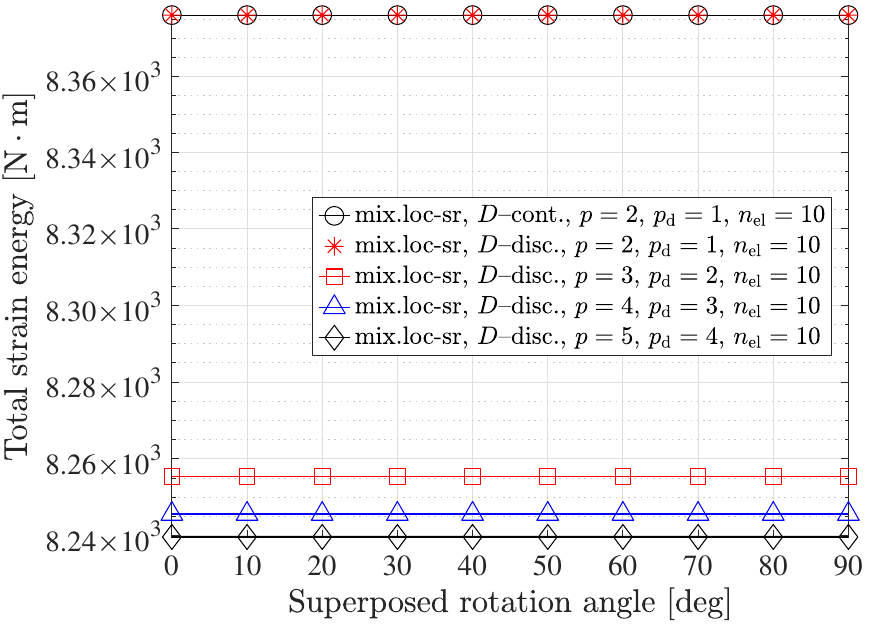}
		\caption{Beam formulation with extensible directors}
		\label{obj_test_bent_total_strain_e_ext_dir}			
	\end{subfigure}	
	\caption{Rigid rotation of a bent rod (objectivity test 2): Change of the total strain energy under a superposed rigid body rotation in the second phase. (a) Here, $n_\mathrm{load}$ represents the number of load steps in the second phase. (b) We use $n_\mathrm{load}=9$ for all cases.}
	\label{obj_test_bent_total_strain_e}	
\end{figure}

We investigate the invariance of non-vanishing strain under rigid body rotation superposed onto a deformed configuration \textcolor{black}{in quasi-statics}. This test was presented in \citet[Sec.\,5.1]{romero2002objective}. We consider a straight beam in its initial (undeformed) configuration with the position vectors of center axis' two end points A and B given by \textcolor{black}{ $\boldsymbol{\varphi}^\mathrm{A}_{0}=\left[0,0,2\,\mathrm{m}\right]^\mathrm{T}$ and $\boldsymbol{\varphi}^\mathrm{B}_{0}=\left[0,0,5\,\mathrm{m}\right]^\mathrm{T}$}, respectively, in a global Cartesian coordinate system. The application of external load is divided into two phases, whose boundary conditions are as follows:
\begin{itemize}
	\item \textbf{First phase.} The displacement of the end point B is prescribed by \textcolor{black}{${\bar {\boldsymbol{u}}}_\mathrm{B}=\left[1\,\mathrm{m},-1\,\mathrm{m},0\right]^\mathrm{T}$}, with the directors at B free, and the cross-section at point A fixed. That is, we apply the boundary conditions
	\begin{subequations}
	\begin{align}
		{\boldsymbol{\varphi }} = {{\boldsymbol{\varphi }}_0},\,\,{{\boldsymbol{d}}_\alpha } = {{\boldsymbol{D}}_\alpha}\,\,(\alpha=1,2)\,\,\text{at point A},
	\end{align}
	\begin{align}
		{\boldsymbol{\varphi }} = {{\boldsymbol{\varphi }}_0} + {{\boldsymbol{\bar u}}_{\rm{B}}}\,\,\text{at point B},
	\end{align}
	and
	\begin{align}
		{{{{{\bar {\tilde {\boldsymbol{m}}}}}}\,\!^\alpha_0}}=\boldsymbol{0}\,\,, \text{i.e.}, {{\boldsymbol{d}}_\alpha }\,\,(\alpha=1,2)\,\,\text{are free} \,\,\text{at point B}.		
	\end{align}
	\end{subequations}
	\item \textbf{Second phase.} In order to investigate the objectivity, a series of nine rigid body rotations are \textcolor{black}{quasi-statically} superposed onto the deformed configuration after the first phase. We apply the boundary conditions
	\begin{align}		
		{\boldsymbol{\varphi }} = {{\boldsymbol{\varphi }}_0},\,\,{{\boldsymbol{d}}_\alpha } = {{\boldsymbol{\bar \Lambda }}_i}{{\boldsymbol{D}}_\alpha }\,\,(\alpha=1,2)\,\,\text{at point A},
	\end{align}
	and
	\begin{align}		
		{\boldsymbol{\varphi }} = {{\boldsymbol{\bar \Lambda}}_i}\left( {{{\boldsymbol{\varphi }}_0} + {{{\boldsymbol{\bar u}}}_{\rm{B}}}} \right),\,\mathrm{with}\,\,{{\boldsymbol{d}}_\alpha }\,\,(\alpha=1,2)\,\,\text{free}\,\,\text{at point B},
	\end{align}	
	for $i=1,2,\cdots,9$, where
	\begin{equation}
		{{\bar{\boldsymbol{\Lambda }}}_i} = \exp \left[ {i\frac{\pi }{{18}}{\boldsymbol{a}}} \right]\in\mathrm{SO(3)},
	\end{equation}
	with the axis of rotation 
	\begin{equation*}
		{\boldsymbol{a}} = \frac{1}{{\sqrt 3 }}{\left[ {1,1,1} \right]^{\rm{T}}}.
	\end{equation*}
\end{itemize}
We consider a St.\,Venant-Kirchhoff type isotropic hyperelastic material with Young's modulus $E=21\,\mathrm{MPa}$ and Poisson's ratio $\nu=0.3$, and a square cross-section of dimension $d=0.1\,\mathrm{m}$. In this example, since the beam's initial geometry is straight, the two approaches [$D$--cont.] and [$D$--disc.] yield the same initial director fields. In Fig.\,\ref{obj_test_bent_total_strain_e_ext_dir_sv}, it is seen that the beam formulation [$\Delta\theta$] suffers from spurious strain energy in superposed rigid rotation, that is, it is non-objective. It is noticeable that the numerical error decreases by increasing number of load steps, $n_\mathrm{load}$ in the second phase, where we use a uniform load increment. The reason for this is that the error in the finite element approximation of an incremental rotation should decrease, as the magnitude of rotation decreases. This can be seen, in the much lower spurious energy resulting from using doubled $n_\mathrm{load}$ (cyan curve with star), compared with that of using $n_\mathrm{load}=18$ (black curve with circle). This non-objectivity can be also alleviated by mesh refinement, a degree elevation (blue curve with triangle) or a larger number of elements (red curve with square). In contrast, the beam formulation with a direct finite element approximation of the director vectors [mix.loc-sr] exhibits no spurious strain energy in the second phase. In Fig.\,\ref{obj_test_bent_total_strain_e_ext_dir}, both methods [$D$--cont.] and [$D$--disc.] give the same results and have constant strain energy under superposed rigid rotation (red star and black curve with circle). Fig.\,\ref{obj_test_bent_total_strain_e_ext_dir} also verifies the objectivity in higher degrees $p=3,4,5$ and $p_\mathrm{d}=p-1$ (red, blue, and black curves).
\subsection{Straight beam under twisting moment}
\label{num_ex_twist_straight}
\begin{figure}[H]
	\centering
	\includegraphics[width=0.65\linewidth]{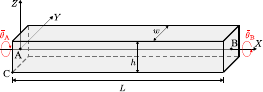}
	\caption{{Straight beam under twisting moment: Initial geometry and boundary conditions. ${\bar \theta}_\mathrm{A}$ and ${\bar \theta}_\mathrm{B}$ denote the prescribed rotation angles of the cross-sections at the points A and B with respect to the center axis, respectively.}}	
	\label{sv_torsion_init_geom}	
\end{figure}

\begin{figure}[H]
	\centering
	\begin{subfigure}[b]{\linewidth}\centering
		\includegraphics[width=\linewidth]{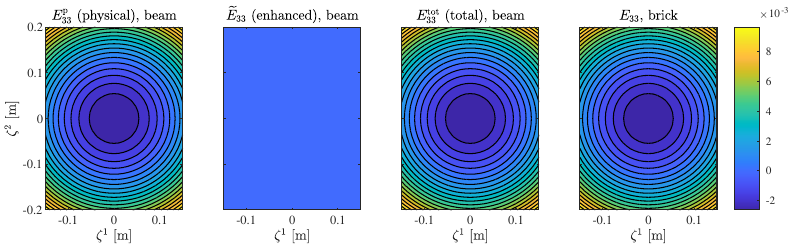}
		\caption{Strain field}
	\end{subfigure}
	\begin{subfigure}[b]{\linewidth}\centering
		\includegraphics[width=0.3\linewidth]{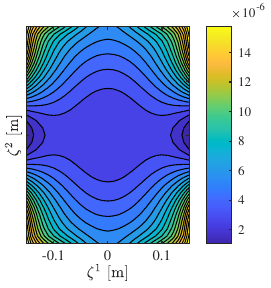}		
		\caption{Difference between the beam and brick solutions, $E^\mathrm{tot}_{33}-E_{33}$}		
	\end{subfigure}
	\caption{Straight beam under twisting moment for $w=0.3\,\mathrm{m}$ and $h=0.4\,\mathrm{m}$: Comparison of the axial normal strain, $E_{33}$ at $s=L/2$. \textcolor{black}{For the beam solution, we enhance ${\widetilde E}_{33}$ by the basis ${\bf{w}}_3$ with the maximum degree $m_3=4$ and the minimum degree ${\bar m}_3+1=3$, which is orthogonal to polynomials up to degree ${\bar m}_3=2$. The maximum value of $\left|{\widetilde E}_{33}\right|$ is $4.0\times10^{-17}\,[-]$. There is excellent agreement between the total strains for beam and brick.}}
	\label{sv_tor_cs_strn_E33_w0.3_h0.4}	
\end{figure}

\begin{figure}[H]
	\centering
	\includegraphics[width=\linewidth]{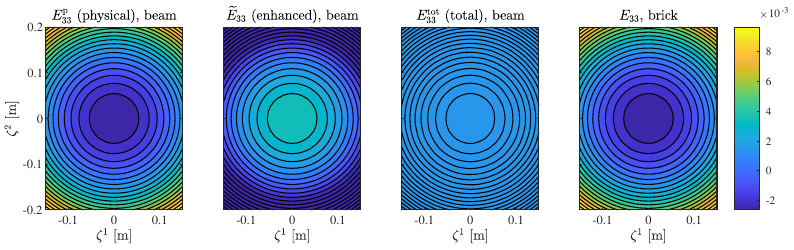}
	\caption{Straight beam under twisting moment for $w=0.3\,\mathrm{m}$ and $h=0.4\,\mathrm{m}$: Comparison of the axial normal strain, $E_{33}$ at $s=L/2$. \textcolor{black}{Brick solution is the same as in Fig.\,\ref{sv_tor_cs_strn_E33_w0.3_h0.4}. For the beam solution, we enhance ${\widetilde E}_{33}$ by the basis ${\bf{w}}_4={\bf{w}}_3$ with the maximum degree $m_3=4$ and the minimum degree ${\bar m}_3+1=2$, which is orthogonal to polynomials up to degree ${\bar m}_3=1$. The resulting enhanced strain ${\widetilde E}_{33}$ leads to erroneous total strain ($E^\mathrm{tot}_{33}$).}}
	\label{sv_tor_cs_strn_E33_w0.3_h0.4_ortho_up_to_lin}	
\end{figure}

\begin{figure}[H]
	\centering
	\begin{subfigure}[b]{\linewidth}\centering
		\includegraphics[width=\linewidth]{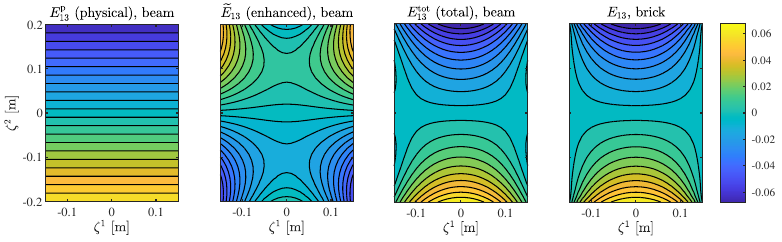}
		\caption{Strain field}		
	\end{subfigure}
	\begin{subfigure}[b]{\linewidth}\centering
		\includegraphics[width=0.3\linewidth]{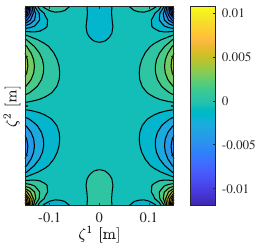}		
		\caption{Difference between the beam and brick solutions, $E^\mathrm{tot}_{13}-E_{13}$}		
	\end{subfigure}	
	\caption{Straight beam under twisting moment for $w=0.3\,\mathrm{m}$ and $h=0.4\,\mathrm{m}$: Comparison of the transverse shear strain, $E_{13}$ at $s=L/2$. \textcolor{black}{For the beam solution, we enhance ${\widetilde E}_{13}$ by the derivative ${\bf{w}}_{4,1}$, where the maximum and minimum degrees of ${\bf{w}}_4$ are ${m_4}=4$ and ${\bar m}_4+1=2$, respectively. There is excellent agreement between the total strains for beam and brick.}}
	\label{sv_tor_cs_strn_E13_w0.3_h0.4}	
\end{figure}

\begin{figure}[H]
	\centering
	\begin{subfigure}[b]{\linewidth}\centering	
		\includegraphics[width=\linewidth]{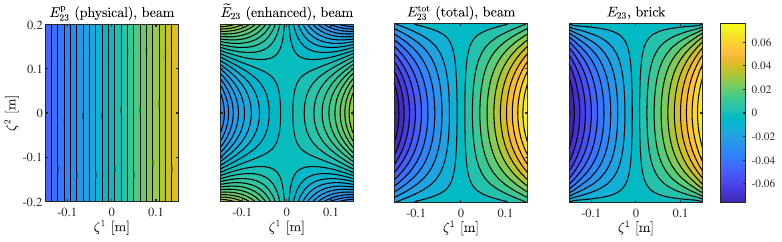}
		\caption{Strain field}
	\end{subfigure}
	\begin{subfigure}[b]{\linewidth}\centering	
		\includegraphics[width=0.3\linewidth]{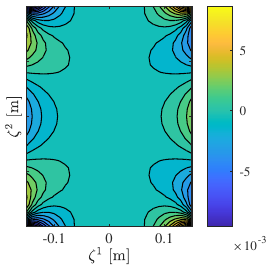}		
		\caption{Difference between the beam and brick solutions, $E^\mathrm{tot}_{23}-E_{23}$}		
	\end{subfigure}			
	\caption{Straight beam under twisting moment for $w=0.3\,\mathrm{m}$ and $h=0.4\,\mathrm{m}$: Comparison of the transverse shear strain component, $E_{23}$ at $s=L/2$. \textcolor{black}{For the beam solution, we enhance ${\widetilde E}_{23}$ by the derivative ${\bf{w}}_{4,2}$, where the maximum and minimum degrees of ${\bf{w}}_{4}$ are ${m_4}=4$ and ${{\bar m}_4+1=2}$. There is excellent agreement between the total strains for beam and brick.}}
	\label{sv_tor_cs_strn_E23_w0.3_h0.4}	
\end{figure}

\begin{figure}[H]
	\centering
	\includegraphics[width=\linewidth]{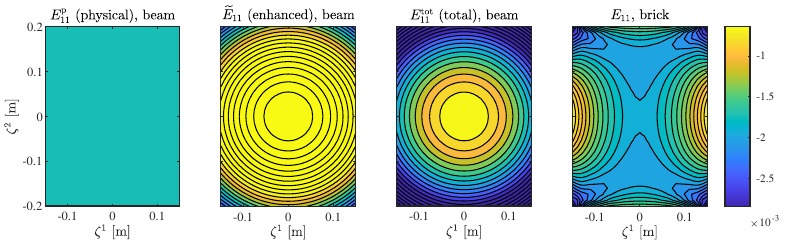}
	\caption{Straight beam under twisting moment for $w=0.3\,\mathrm{m}$ and $h=0.4\,\mathrm{m}$: Comparison of the transverse normal strain, $E_{11}$ at $s=L/2$. For the beam solution, we enhance ${\widetilde E}_{11}$ \textcolor{black}{by ${\bf{w}}_1$ with the maximum and minimum degrees $m_1=4$ and ${\bar m}_1+1=1$, respectively.}}
	\label{sv_tor_cs_strn_E11_w0.3_h0.4}	
\end{figure}

\begin{figure}[H]
	\centering
	\includegraphics[width=\linewidth]{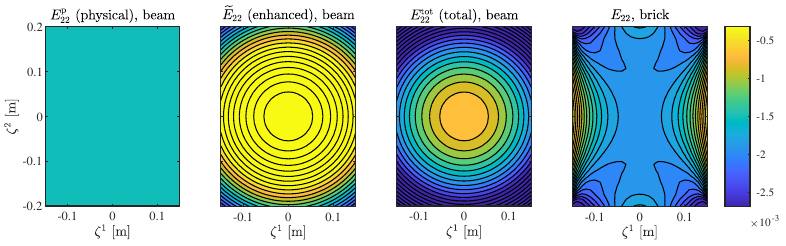}
	\caption{Straight beam under twisting moment for $w=0.3\,\mathrm{m}$ and $h=0.4\,\mathrm{m}$: Comparison of the transverse normal strain, $E_{22}$ at $s=L/2$. For the beam solution, we enhance ${\widetilde E}_{22}$ \textcolor{black}{by ${\bf{w}}_1$ with the maximum and minimum degrees $m_1=4$ and ${\bar m}_1+1=1$, respectively.}}
	\label{sv_tor_cs_strn_E22_w0.3_h0.4}	
\end{figure}

\begin{figure}[H]
	\centering
	\includegraphics[width=\linewidth]{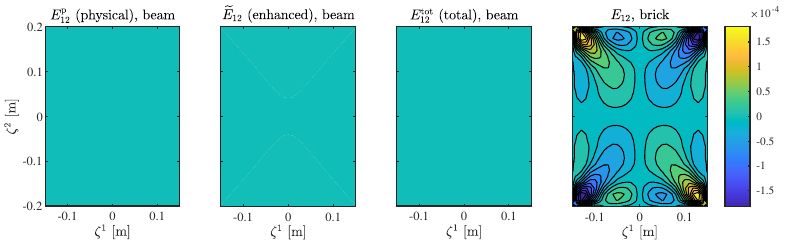}
	\caption{Straight beam under twisting moment for $w=0.3\,\mathrm{m}$ and $h=0.4\,\mathrm{m}$: Comparison of the in-plane shear cross-sectional strain component, $E_{12}$ at $s=L/2$. For the beam solution, we enhance ${\widetilde E}_{12}$ \textcolor{black}{by ${\bf{w}}_2$ with the maximum and minimum degrees $m_2=4$ and ${\bar m}_2+1=1$, respectively.}}
	\label{sv_tor_cs_strn_E12_w0.3_h0.4}	
\end{figure}

We consider a straight beam along the $X$-axis that is subject to a prescribed twisting angle $\bar\theta=2\pi\,[\mathrm{rad}]$, \textcolor{black}{which is applied quasi-statically}, about its center axis \textcolor{black}{at both ends A and B in opposite directions}. The beam has an initial length $L=10\,\mathrm{m}$, and a rectangular cross-section of width $w$ and height $h$, see Fig.\,\ref{sv_torsion_init_geom}. For the following cases of $w$ and $h$, we verify that the present beam formulation can represent the correct torsional stiffness \textcolor{black}{without additional global DOFs}. 
\begin{itemize}
	\item Case 1: $w=0.3\,\mathrm{m}$, and $h=0.4\,\mathrm{m}$,
	\item Case 2: $w=1/3\,\mathrm{m}$, and $h=1\,\mathrm{m}$.
\end{itemize}
A St.\,Venant-Kirchhoff type isotropic hyperelastic material with Young's modulus $E=210\,\mathrm{GPa}$ and Poisson's ratio $\nu=0.3$ is considered. 
\noindent \textcolor{black}{For comparison, we consider two different boundary conditions:
\begin{itemize}
	\item {[free-free]} \textcolor{black}{Prescribe the rotation angles ${\bar \theta}_\mathrm{A}=-{\bar\theta}/2$ and ${\bar \theta}_\mathrm{B}={\bar\theta}/2$.} Note that, here, out-of-plane cross-sectional warping is allowed at both ends, A and B. For the brick solution, the $X$-displacement is constrained at C in order to prevent a rigid-body translation.
	\item {[fixed-free]} Fix the cross-section at A (${\bar \theta}_\mathrm{A}=0$), and prescribe the rotation angle \textcolor{black}{${\bar \theta}_\mathrm{B}={\bar\theta}$} at B. Note that we still allow for out-of-plane cross-sectional warping at the loaded end, B.
\end{itemize}  
It should be noted that, in both cases of the boundary conditions, in-plane cross-sectional strains are not allowed at the ends. For beams, we use the [free-free] only, where the cross-sectional warping due to the enhanced strains is allowed at both ends. An application of boundary conditions to constrain the enhanced strain field remains future work.} Figs.\,\ref{sv_tor_cs_strn_E33_w0.3_h0.4}--\ref{sv_tor_cs_strn_E12_w0.3_h0.4} compare the cross-sectional strains at the middle of the center axis ($s=L/2$) between the present beam and brick solutions for case 1 ($w=0.3\,\mathrm{m}$ and $h=0.4\,\mathrm{m}$). Here, for the beam solution, we use [mix.loc-sr] with $p=3$, ${m_1}={m_2}={m_3}={m_4}=4$, and $n_\mathrm{el}=40$. For the brick solution, we use deg.=$(3,3,3)$ and $n_\mathrm{el}=80\times10\times10$\textcolor{black}{, and the [free-free] condition.}
\vspace{2mm}
\begin{observation}
	\label{rem_num_ex_enrich_e33}
	In Fig.\,\ref{sv_tor_cs_strn_E33_w0.3_h0.4}, it is seen that the enhanced strain ${\widetilde E}_{33}$ vanishes. This implies that the quadratic polynomial basis in ${\bf{A}}(\zeta^1,\zeta^2)$ in the third row of Eq.\,(\ref{mat_A_poly_basis}) is sufficient to represent the axial normal strain. It has been also discussed in \citet{wackerfuss2011nonlinear} that the enhancement of $E_{33}$ is not needed for a correct representation of the three-dimensional strain state. Accordingly, they presented a so called `E2--model' which does not enrich the axial normal strain. This significantly reduces the number of internal DOFs for the enhanced strain field, which makes the formulation more efficient. 
\end{observation}
\vspace{2mm}
\noindent In Figs.\,\ref{sv_tor_cs_strn_E13_w0.3_h0.4} and \ref{sv_tor_cs_strn_E23_w0.3_h0.4}, the physical transverse shear strains $E^\mathrm{p}_{13}$ and $E^\mathrm{p}_{23}$ are limited to linear distributions, since it has polynomial basis up to degree 1. It is noticeable that, by employing the enhanced strains, the total strains show excellent agreement with the brick solution, in both cases, $E^\mathrm{tot}_{13}$ and $E^\mathrm{tot}_{23}$.
\noindent \textcolor{black}{In Fig.\,\ref{sv_tor_cs_strn_E33_w0.3_h0.4_ortho_up_to_lin}, we show the beam solution using ${\bar m}_3=1$ instead of ${\bar m}_3=2$, with all other parameters the same as in the result of Fig.\,\ref{sv_tor_cs_strn_E33_w0.3_h0.4}. This aims at showing that losing the orthogonality between $E^\mathrm{p}_{33}$ and ${\widetilde E}_{33}$ can lead to significant error. Selecting the lower degree ${\bar m}_3=1$ means that the basis ${\bf{w}}_3$ may have quadratic terms in spite of the existing quadratic terms in $E^\mathrm{p}_{33}$. This eventually pollutes the physical strain field, which can be seen by the deviation of the resulting total strain from the brick solution.}


\begin{figure}[H]
	\centering
	\begin{subfigure}[b]{0.49\textwidth}\centering
		\includegraphics[width=\linewidth]{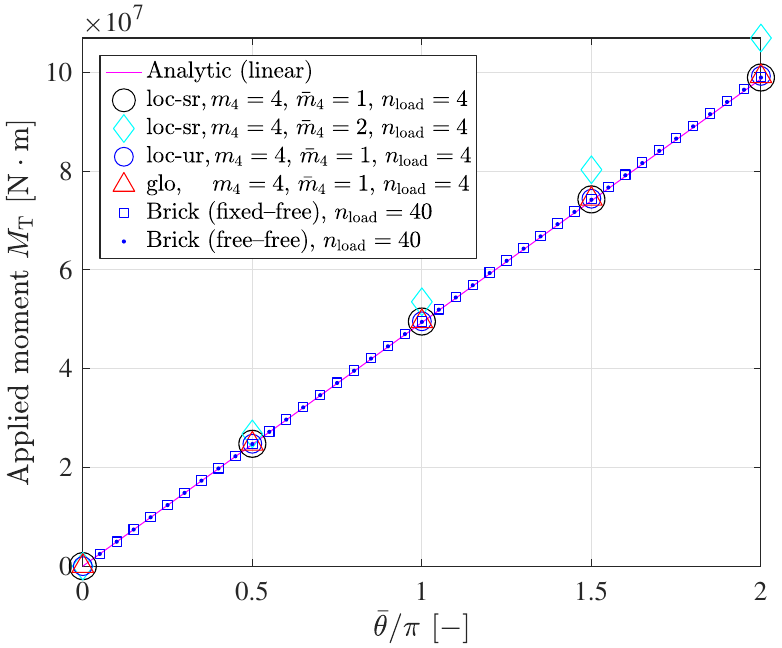}
		\caption{Case 1: $w=0.3\,\mathrm{m}$ and $h=0.4\,\mathrm{m}$}
		\label{sv_torsion_w0.3_h0.4}			
	\end{subfigure}
	\begin{subfigure}[b]{0.49\textwidth}\centering
		\includegraphics[width=\linewidth]{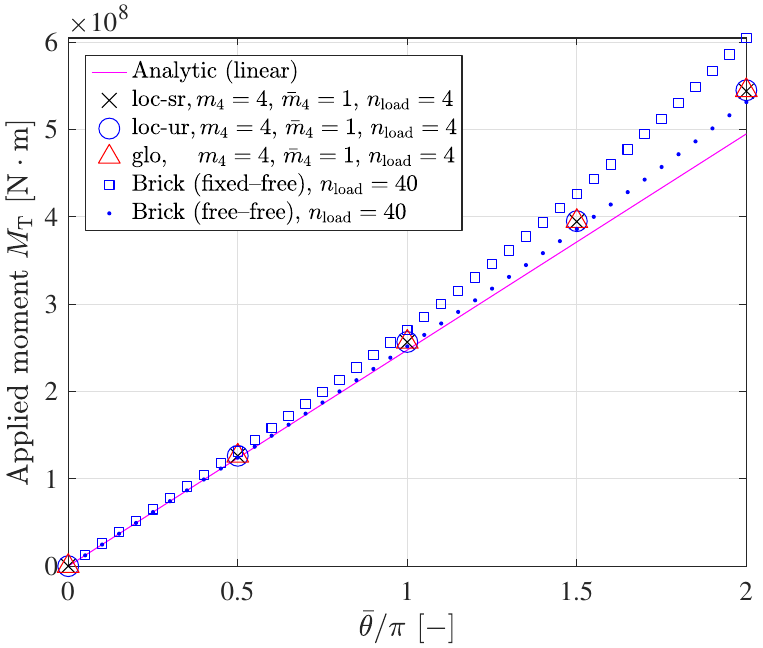}
		\caption{Case 2: $w=1/3\,\mathrm{m}$ and $h=1\,\mathrm{m}$}
		\label{sv_torsion_w0.333_h1}			
	\end{subfigure}
	\caption{Straight beam under twisting moment: Change of the applied moment along the prescribed rotation $\bar \theta$. In [mix.loc-ur] and [mix.loc-sr], we use $p=3$, ${p_\mathrm{d}}=2$. For [mix.glo], we use $p=p_d=3$. All the beam solutions are from using $n_\mathrm{el}=40$, and ${m}_1=m_2=2$, without enriching the axial normal strain, $E_{33}$. For the brick solutions, we use deg.$=(3,3,3)$, and $n_\mathrm{el}=80\times10\times10$.}
	\label{sv_tor_react_mnt}	
\end{figure}

\begin{figure}[H]
	\centering
	\begin{subfigure}[b]{0.495\textwidth}\centering
		\includegraphics[width=\linewidth]{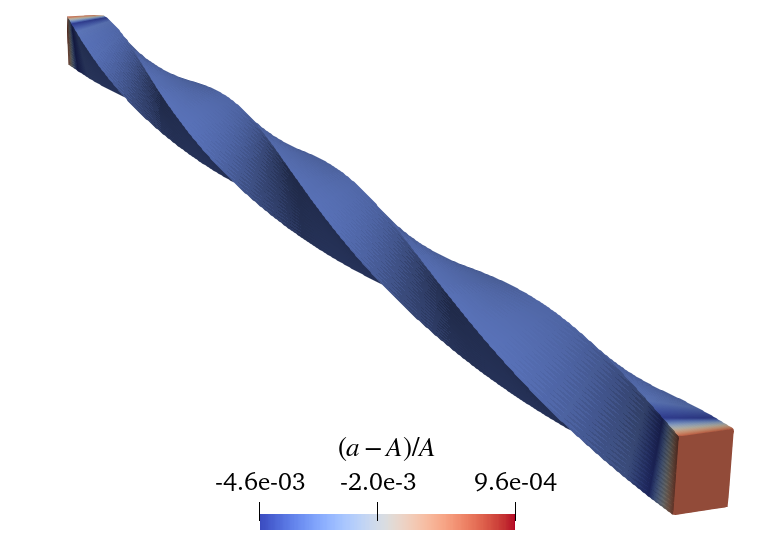}
		\caption{Beam solution (381 DOFs)}
		\label{sv_torsion_deformed_w0.3_h0.4}			
	\end{subfigure}
	\begin{subfigure}[b]{0.495\textwidth}\centering
		\includegraphics[width=\linewidth]{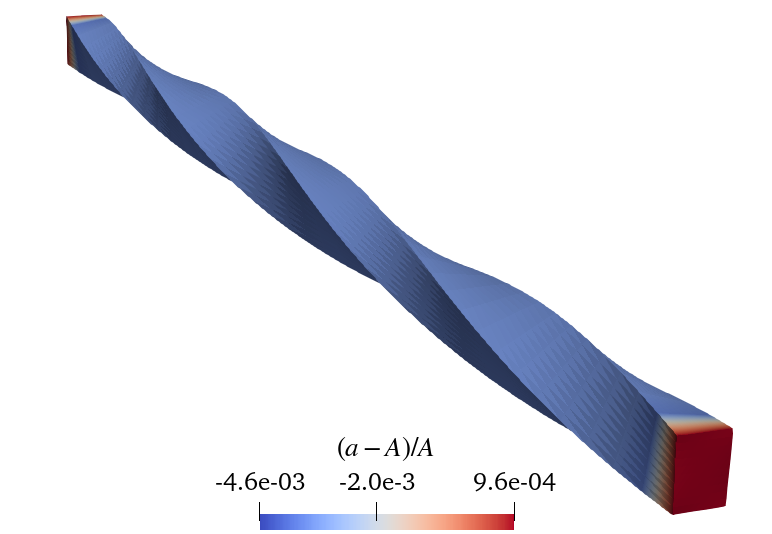}
		\caption{Brick solution (42,081 DOFs)}
		\label{sv_torsion_deformed_w0.3_h0.4_brick}			
	\end{subfigure}
	\caption{Straight beam under twisting moment for $w=0.3\,\mathrm{m}$ and $h=0.4\,\mathrm{m}$ (case 1): \textcolor{black}{Final deformed configuration at ${\bar \theta}_\mathrm{A}=-\pi$ and ${\bar \theta}_\mathrm{B}=\pi$. The colors represent the relative change of the cross-sectional area, where $a$ and $A$ denote the deformed and initial cross-sectional areas, respectively. For the beam solution, we use [mix.loc-sr] with $p=3$, $n_\mathrm{el}=40$, and ${m_1}={m_2}=2$, $m_4=4$, where $E_{33}$ has not been enriched.} For the brick solution, we use deg.$=(3,3,3)$, and $n_\mathrm{el}=80\times10\times10$. There is excellent agreement between beam and brick solutions apart from the two end faces.}
	\label{sv_tor_react_defomred_w0.3_h0.4}	
\end{figure}

\begin{figure}[H]
	\centering
	\begin{subfigure}[b]{0.45\textwidth}\centering
		\includegraphics[width=\linewidth]{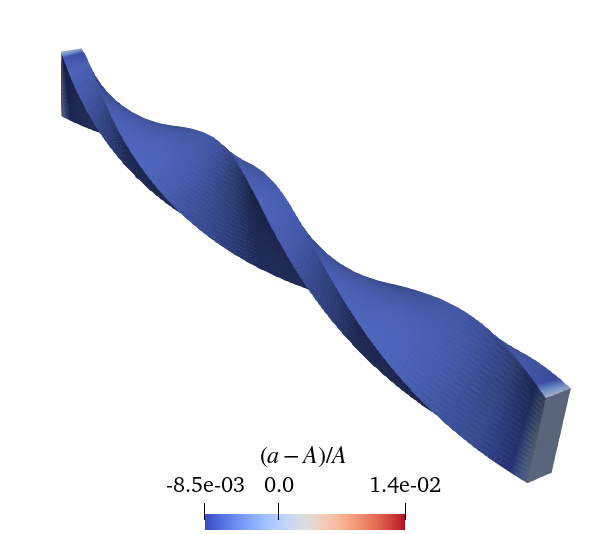}
		\caption{Beam solution (381 DOFs)}
		\label{sv_torsion_deformed_w0.333_h1}			
	\end{subfigure}
	\begin{subfigure}[b]{0.45\textwidth}\centering
		\includegraphics[width=\linewidth]{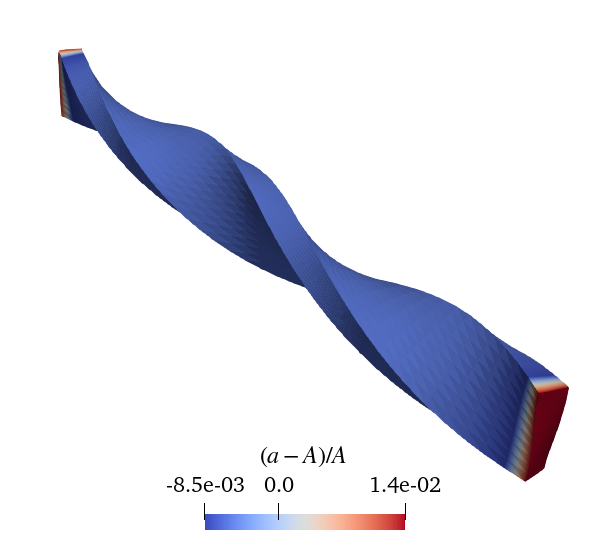}
		\caption{Brick solution (42,081 DOFs)}
		\label{sv_torsion_deformed_w0.333_h1_brick}			
	\end{subfigure}
	\caption{Straight beam under twisting moment for $w=1/3\,\mathrm{m}$ and $h=1\,\mathrm{m}$ (case 2): \textcolor{black}{Final deformed configuration at ${\bar \theta}_\mathrm{A}=-\pi$ and ${\bar \theta}_\mathrm{B}=\pi$. The colors represent the relative change of the cross-sectional area, where $a$ and $A$ denote the deformed and initial cross-sectional areas, respectively. For the beam solution, we use [mix.loc-sr] with $p=3$, $n_\mathrm{el}=40$, and ${m_1}={m_2}=2$, $m_4=4$, where $E_{33}$ has not been enriched.} For the brick solution, we use deg.$=(3,3,3)$, and $n_\mathrm{el}=80\times10\times10$. There is excellent agreement between beam and brick solutions apart from the two end faces.}
	\label{sv_tor_react_defomred_w0.333_h1}	
\end{figure}

\noindent In Figs.\,\ref{sv_tor_cs_strn_E11_w0.3_h0.4} and \ref{sv_tor_cs_strn_E22_w0.3_h0.4}, the physical in-plane cross-sectional strains $E^\mathrm{p}_{11}$ and $E^\mathrm{p}_{22}$, respectively, have negative constant values, which represents a uniform transversal contraction of the cross-section. Although the enhanced strains improve the degree of approximation, they are limited to a quadratic distribution, due to the coupling with the axial normal strain $E^\mathrm{p}_{33}$ by the non-zero Poisson's ratio. Further, in case of the in-plane shear strain, the enhanced strain ${\widetilde E}_{12}$ turns out not to be activated, as Fig.\,\ref{sv_tor_cs_strn_E12_w0.3_h0.4} shows. Further investigation on the selection of the basis ${{\bf{w}}}_1$ and ${\bf{w}}_2$ for improving the enhancement of in-plane cross-sectional strain remains future work.
\noindent In Fig.\,\ref{sv_tor_react_mnt}, we compare the applied moment for a given prescribed rotation angle $\bar \theta$ between the present beam element solution, the reference brick element solution, and the analytical solution in Eq.\,(\ref{sol_sv_tor}). In the first case ($w=0.3\,\mathrm{m}$ and $h=0.4\,\mathrm{m}$), the beam solutions show excellent agreement with both analytic and brick element solutions. \textcolor{black}{In Fig.\,\ref{sv_torsion_w0.3_h0.4}, the beam solution with higher degree ${\bar m_4}=2$ (cyan diamonds) shows overly stiff behavior, due to the missing quadratic terms in the enriched basis ${\bf{w}}_4$}. Overall, the mixed finite element formulation enables much larger load increments, compared to the displacement-based brick element formulation. In the second case ($w=1/3\,\mathrm{m}$ and $h=1\,\mathrm{m}$), a nonlinear effect of increasing torsional stiffness along the prescribed rotation is clearly seen in both beam and brick solutions. \textcolor{black}{In the latter range of the prescribed rotation, it is seen that the brick solution with the [fixed--free] boundary condition exhibits more rapidly increasing torsional stiffness compared to the solution using the [free-free] condition, due to the constrained out-of-plane warping at $s=0$. It is remarkable that the beam and brick solutions under the [free--free] condition show excellent agreement. In Table\,\ref{tab_str_sv_tor}, we compare the total number of global and internal DOFs between beam and brick solutions. It is noticeable that the beam solutions use much smaller number of DOFs than the brick solution does, since the DOFs of the enhanced strain ($\boldsymbol{\alpha}$) can be element-wisely condensed out. In the beam formulation [mix.glo], the DOFs of the physical stress resultants ($\bf{r}$) and physical strain ($\bf{e}$) are considered in the global system of equations, whose number is much smaller than those of the local approaches for the same number of elements, due to the higher-order continuity of B-spline basis functions, see Remark \ref{rem_patch_nofs_compare}.} 	
\begin{table}[h]
	\centering
	\caption{\textcolor{black}{Straight beam under twisting moment: Comparison of total number of DOFs for cases 1 and 2. For our beam formulation, we have DOFs for the configuration ($\boldsymbol{y}$), physical stress resultant ($\boldsymbol{r}_\mathrm{p}$), physical strain ($\boldsymbol{\varepsilon}_\mathrm{p}$), and enhanced strain ($\boldsymbol{\alpha}$). The gray cells represent numbers of internal DOFs.}}
	\begin{tabular}{lcccc}
		\toprule
		& Beam [mix.loc-sr] & Beam [mix.loc-ur] & Beam [mix.glo] & Brick \\
		\midrule
		Degree of basis & $p=3$, $p_\mathrm{d}=2$ & $p=3$, $p_\mathrm{d}=2$ & $p=p_\mathrm{d}=3$ & $\mathrm{deg.}=(3,3,3)$ \\
		\#elements & 40    & 40    & 40    & $80\times10\times10$ \\
		\#DOFs ($\boldsymbol{y}$) & 381   & 381   & 387   & 42,081 \\
		\#DOFs ($\boldsymbol{r}_\mathrm{p}$) & \cellcolor[rgb]{ .851,  .851,  .851}814 & \cellcolor[rgb]{ .851,  .851,  .851}1200 & 630   & -- \\
		\#DOFs ($\boldsymbol{\varepsilon}_\mathrm{p}$) & \cellcolor[rgb]{ .851,  .851,  .851}814 & \cellcolor[rgb]{ .851,  .851,  .851}1200 & 630   & -- \\
		\#DOFs ($\boldsymbol{\alpha}$) & {\cellcolor[rgb]{ .851,  .851,  .851}1,080} & {\cellcolor[rgb]{ .851,  .851,  .851}1,080} &{\cellcolor[rgb]{ .851,  .851,  .851}1,080} & -- \\
		\midrule
		\textbf{\#global DOFs} & \textbf{381} & \textbf{381} & \textbf{1,647} & \textbf{42,081} \\
		\textbf{\#internal DOFs} & \textbf{2,708} & \textbf{3,480} & \textbf{1,080} & -- \\
		\bottomrule
	\end{tabular}
	\label{tab_str_sv_tor}
\end{table}%
\noindent Figs.\,\ref{sv_tor_react_defomred_w0.3_h0.4} and \ref{sv_tor_react_defomred_w0.333_h1} compare the deformed configurations between the beam and brick solutions for the two different dimensions of the initial cross-section. It is noticeable that the same level of decrease of the cross-sectional area is observed in the beam and brick solutions. In the brick solution, it is seen that the cross-sectional area slightly increases at both ends, due to the out-of-plane warping. This is not explicitly shown in the beam solution, since we calculate the cross-sectional area using only two directors ${\boldsymbol{d}}_1$ and ${\boldsymbol{d}}_2$, which can only represent a planar cross-section deformations.

\subsection{Twisting of an elastic ring}
\label{num_ex_twist_ring}
\begin{figure}[h]
	\centering
	\includegraphics[width=0.475\linewidth]{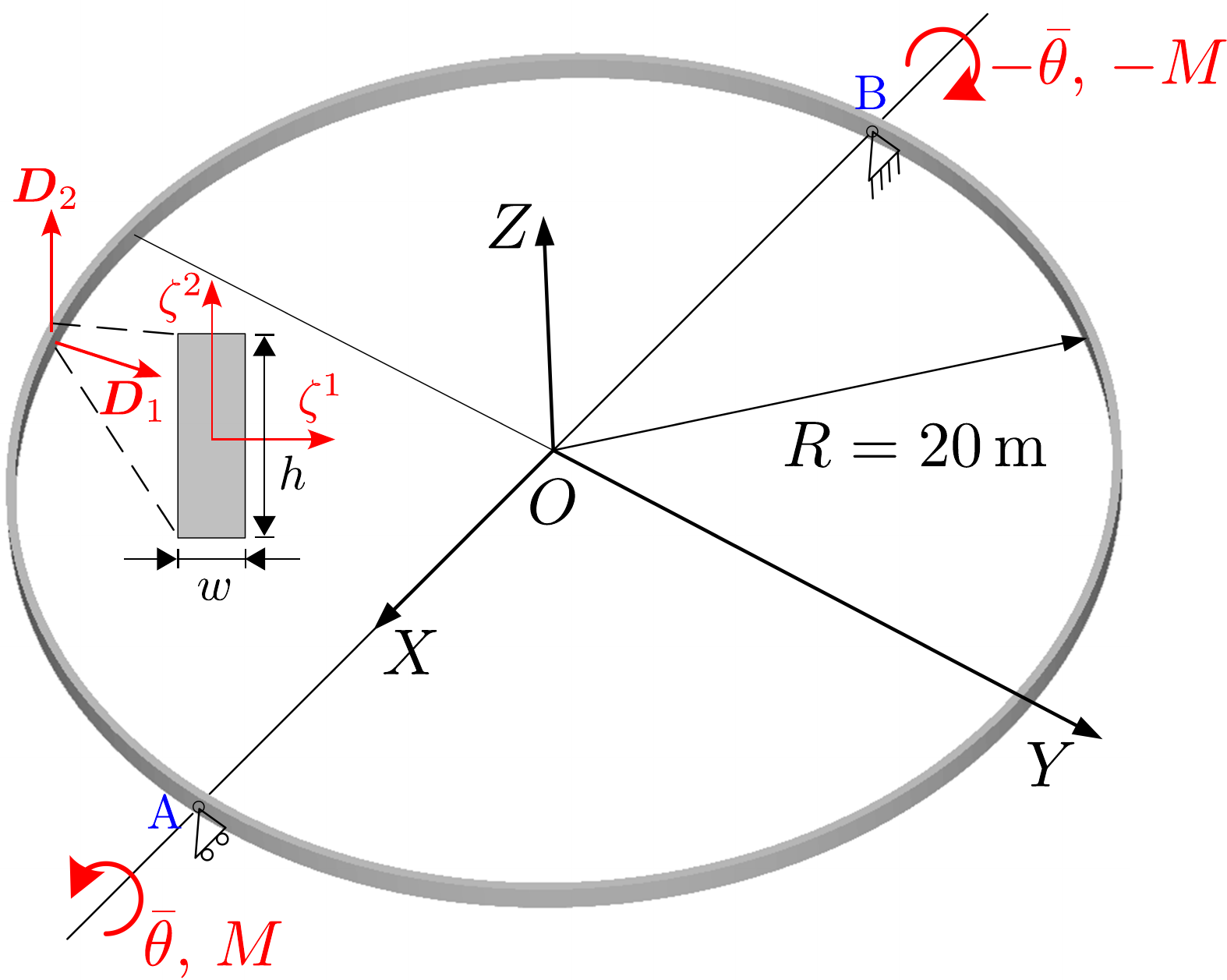}
	\caption{{Twisting of an elastic ring: Undeformed configuration and boundary conditions.}}
	\label{twist_ring_init_config}	
\end{figure}

\begin{figure}[H]
	\centering
	\begin{subfigure}[b]{0.45\textwidth}\centering
		\includegraphics[width=\linewidth]{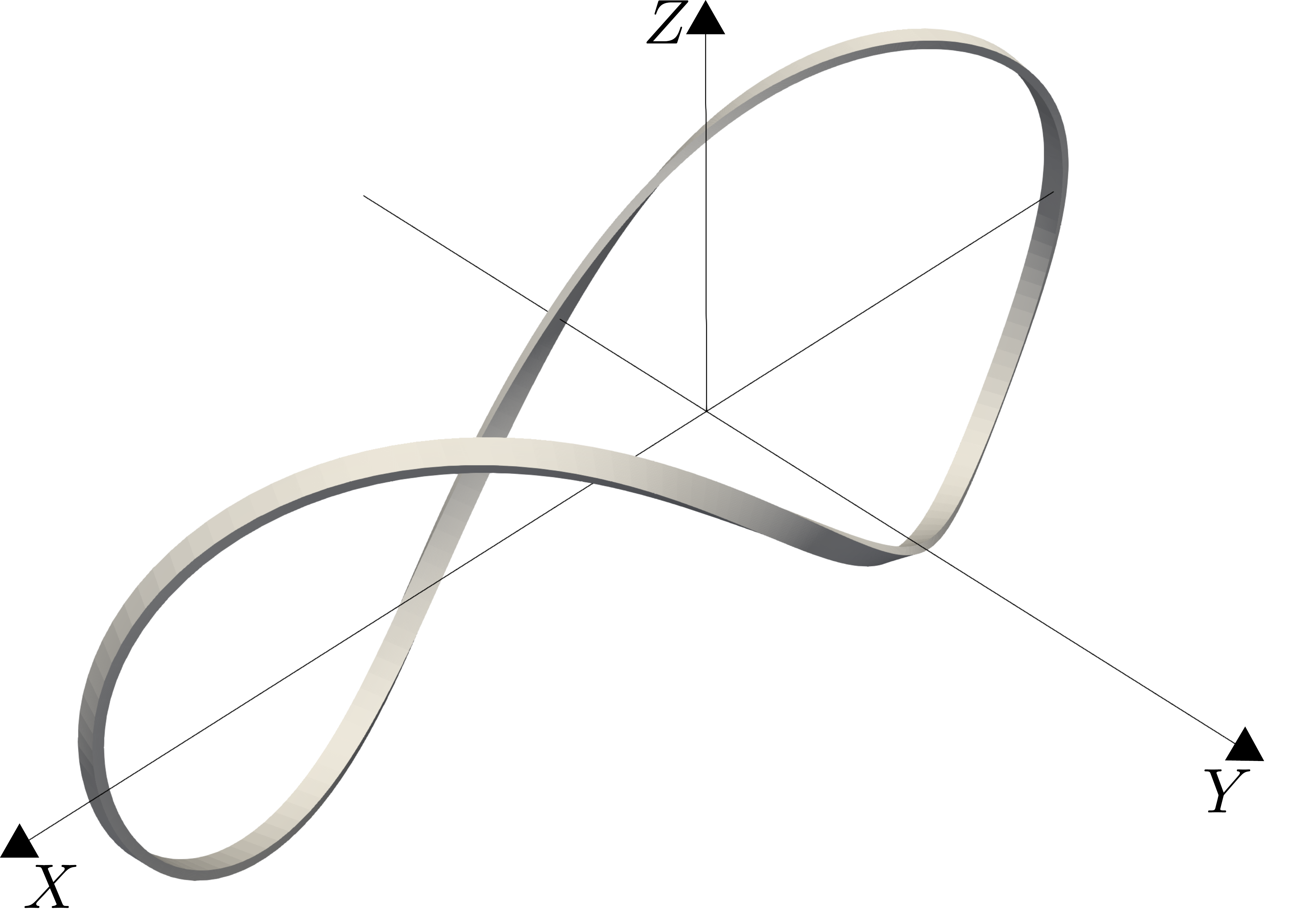}
		\caption{${\bar\theta}=\pi/2$}
		\label{twist_ring_deformed_n4}	
	\end{subfigure}
	\begin{subfigure}[b]{0.45\textwidth}\centering
		\includegraphics[width=\linewidth]{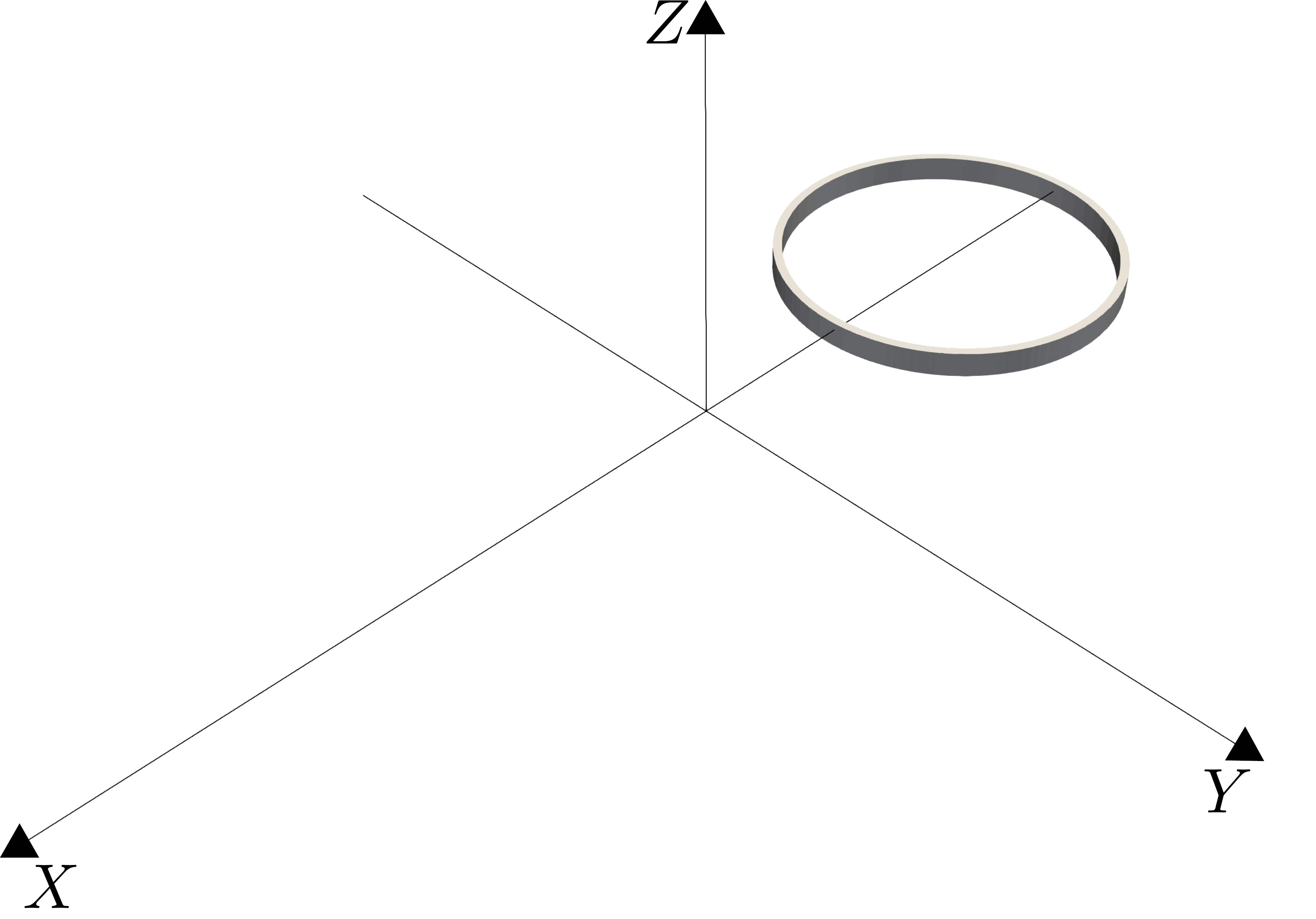}
		\caption{${\bar\theta}=\pi$}
		\label{twist_ring_deformed_n8}	
	\end{subfigure}
	\begin{subfigure}[b]{0.45\textwidth}\centering
		\includegraphics[width=\linewidth]{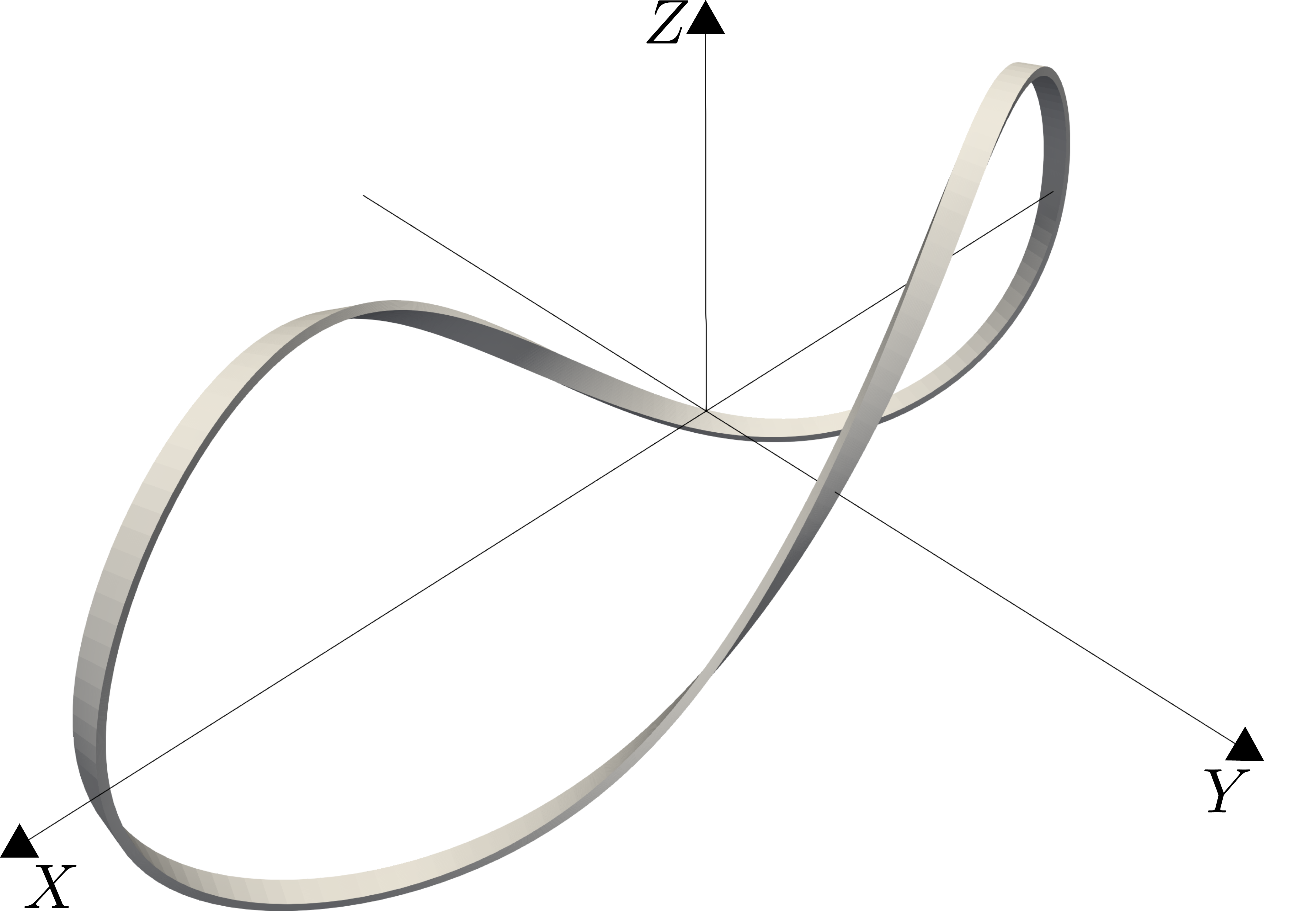}
		\caption{${\bar\theta}=3\pi/2$}
		\label{twist_ring_deformed_n12}	
	\end{subfigure}
	\begin{subfigure}[b]{0.45\textwidth}\centering
		\includegraphics[width=\linewidth]{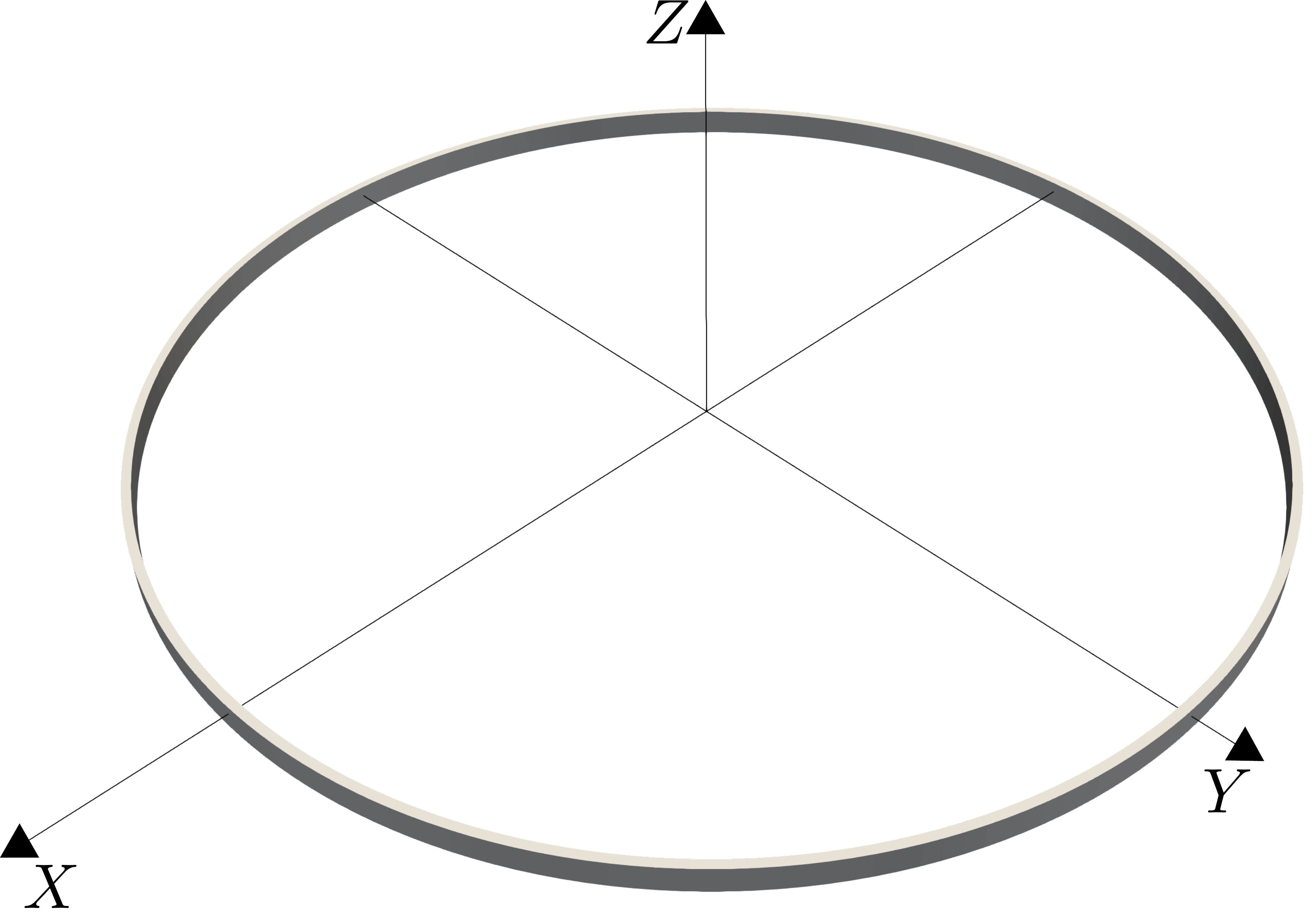}
		\caption{${\bar\theta}=2\pi$}
		\label{twist_ring_deformed_n16}	
	\end{subfigure}	
	\caption{{Twisting of an elastic ring (case 1): Deformed configurations for four different values of $\bar \theta$.}}
	\label{twist_ring_deformed}		
\end{figure}

\begin{figure}[H]
	\centering
	\begin{subfigure}[b]{0.4875\textwidth}\centering
		\includegraphics[width=\linewidth]{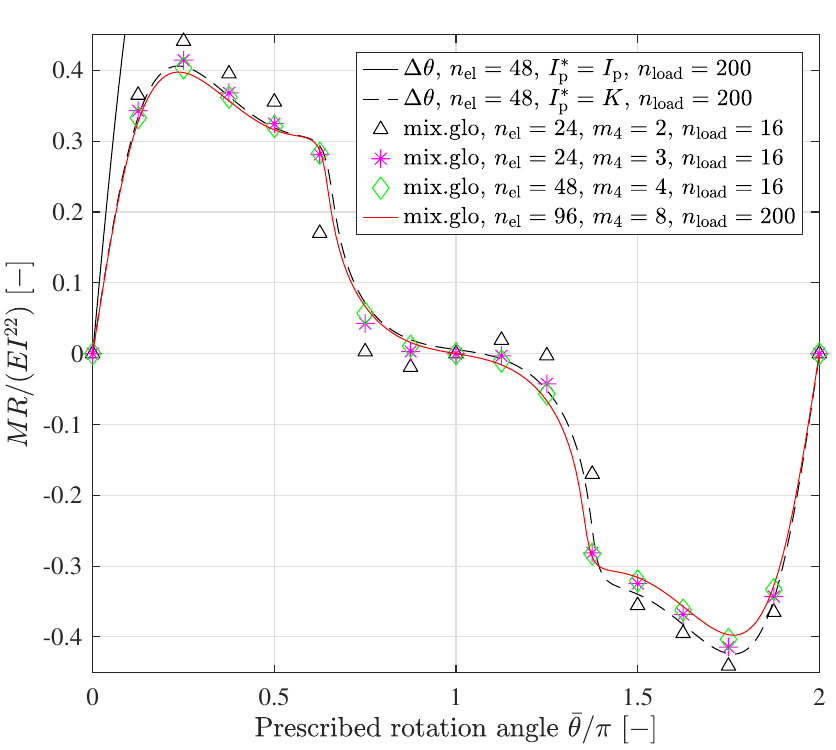}
		\caption{Case 1: $w=1/3\,\mathrm{m}$ and $h=1\,\mathrm{m}$}
		\label{twist_ring_theta-mnt-curve_dim1e-0}			
	\end{subfigure}
	\begin{subfigure}[b]{0.4875\textwidth}\centering
		\includegraphics[width=\linewidth]{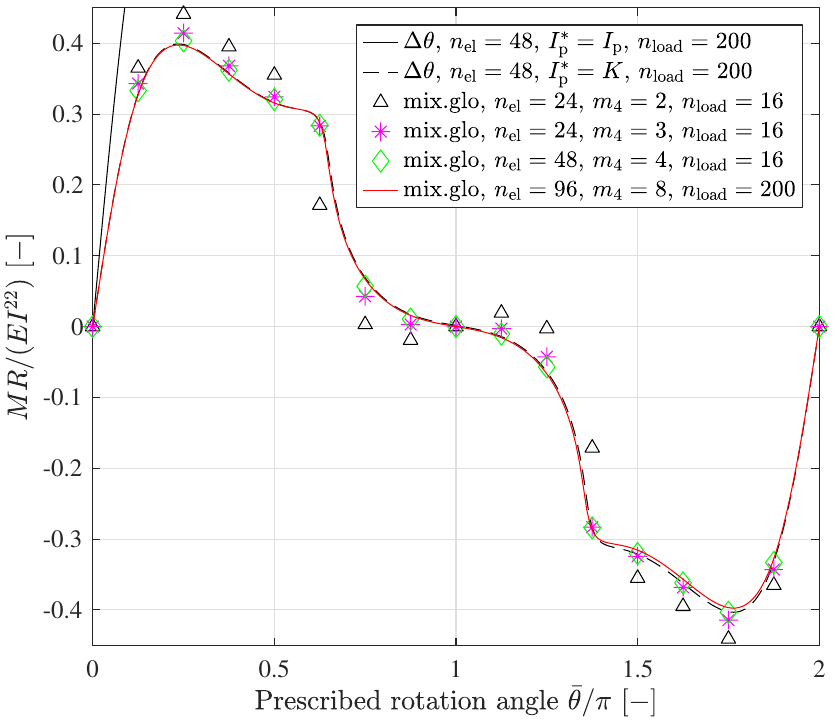}
		\caption{Case 2: $w=1/30\,\mathrm{m}$ and $h=1/10\,\mathrm{m}$}
		\label{twist_ring_theta-mnt-curve_dim1e-1}			
	\end{subfigure}
	\caption{Twisting of an elastic ring: Applied moment at the point A due to the prescribed rotation $\bar \theta$. In both cases, we use $p=1$ and $p=4$ for the formulations [$\Delta \theta$] and [mix.glo], respectively. Further, for [mix.glo], $m_1=m_2=2$, and $E_{33}$ has not been enriched.}
	\label{twist_ring_theta-mnt-curve}
\end{figure}

\begin{figure}[H]
	\centering
	\begin{subfigure}[b]{0.4875\textwidth}\centering
		\includegraphics[width=\linewidth]{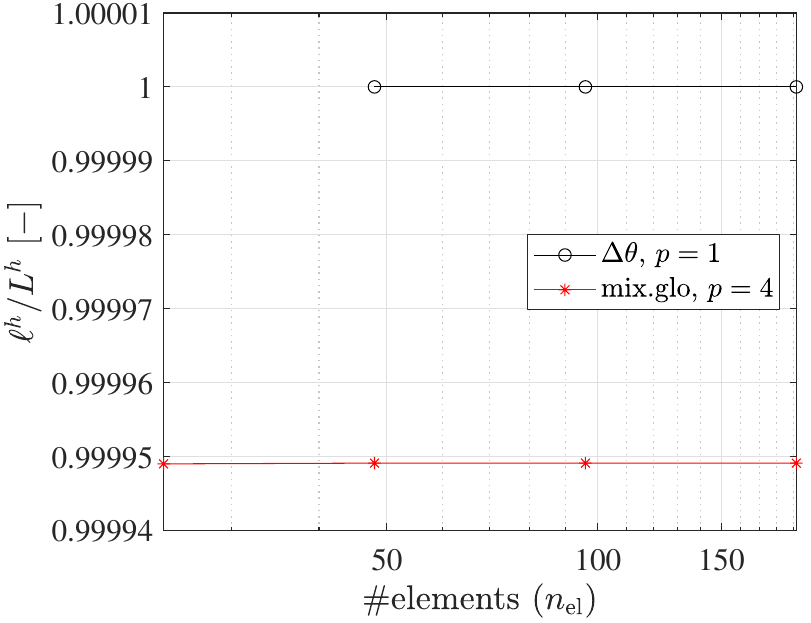}
		\caption{Case 1: $w=1/3$ and $h=1$}
		\label{twist_ring_deformed_length_conv_dim_1e-0}			
	\end{subfigure}
	\begin{subfigure}[b]{0.4875\textwidth}\centering
		\includegraphics[width=\linewidth]{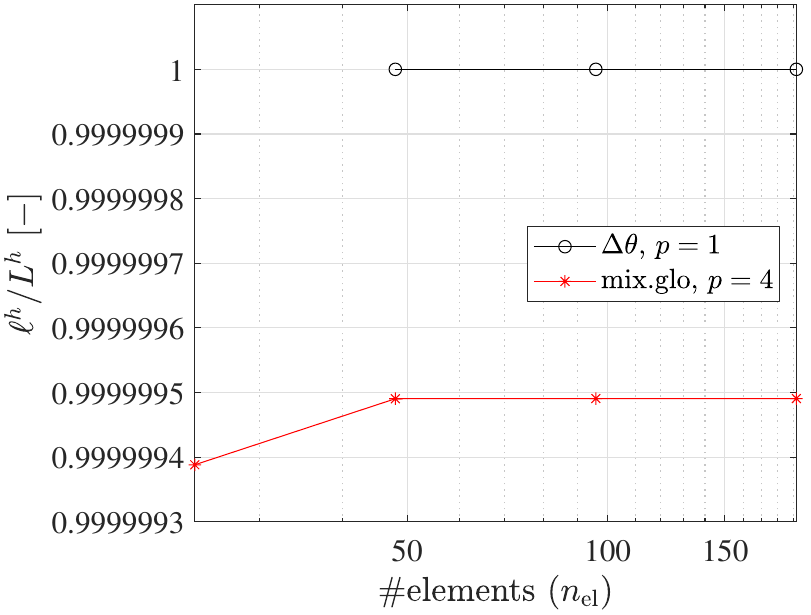}
		\caption{Case 2: $w=1/30$ and $h=1/10$}
		\label{twist_ring_deformed_length_conv_dim_1e-1}			
	\end{subfigure}
	\caption{Twisting of an elastic ring: Change of the ratio between the approximated total lengths of the center axis in the initial ($L^h$) and deformed configuraiton at $\bar \theta=\pi$ ($\ell^h$), for two different cases of $(w,h)$. For [mix.glo], we use $m_1=m_2=2$, $m_4=4$, and $E_{33}$ has not been enriched.}
	\label{twist_ring_deformed_length_conv}	
\end{figure}

\begin{figure}[H]
	\centering
	\begin{subfigure}[b]{0.4875\textwidth}\centering
		\includegraphics[width=\linewidth]{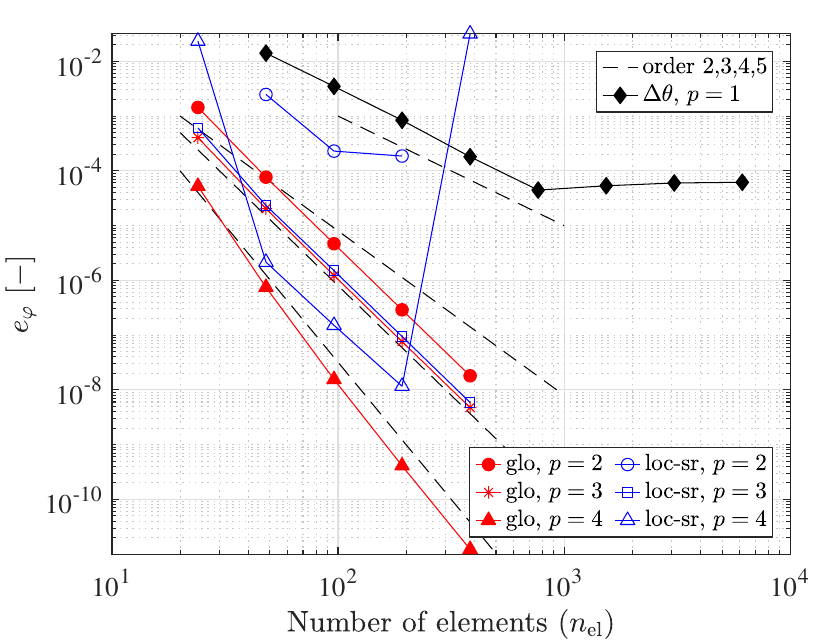}
		\caption{Case 1: $w=1/3$ and $h=1$}
		\label{twist_ring_deformed_rad_conv_dim_1e-0}			
	\end{subfigure}
	\begin{subfigure}[b]{0.4875\textwidth}\centering
		\includegraphics[width=\linewidth]{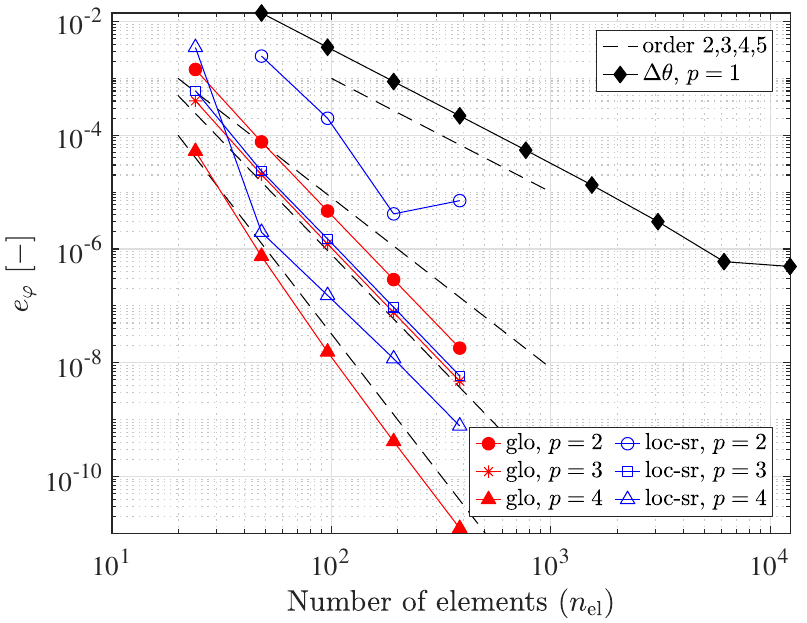}
		\caption{Case 2: $w=1/30$ and $h=1/10$}
		\label{twist_ring_deformed_rad_conv_dim_1e-1}			
	\end{subfigure}
	\caption{Twisting of an elastic ring: Relative difference of the deformed radius from the reference solution $r^*_\mathrm{ref}$ at $\bar\theta=\pi$ for two different cases of $(w,h)$. In both beam formulations [mix.glo] and [mix.loc-sr], we use $m_1=m_2=2$, $m_4=4$, and $E_{33}$ has not been enriched.}
	\label{twist_ring_deformed_rad_conv}	
\end{figure}

\begin{figure}[H]
	\centering
	\begin{subfigure}[b]{0.4875\textwidth}\centering
		\includegraphics[width=\linewidth]{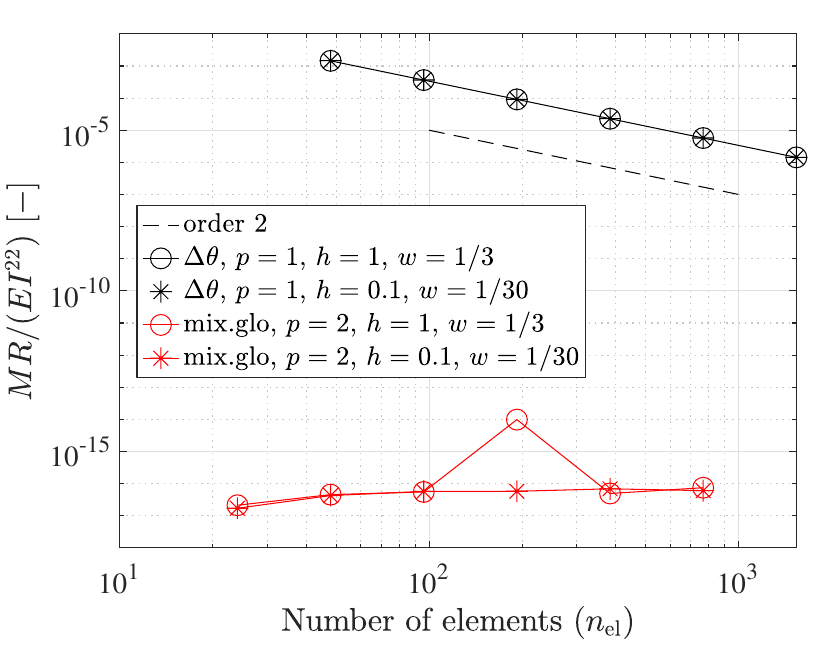}
		\caption{Moment at $\bar \theta=\pi$}
		\label{twist_ring_mnt_tht_pi}			
	\end{subfigure}
	\begin{subfigure}[b]{0.4875\textwidth}\centering
		\includegraphics[width=\linewidth]{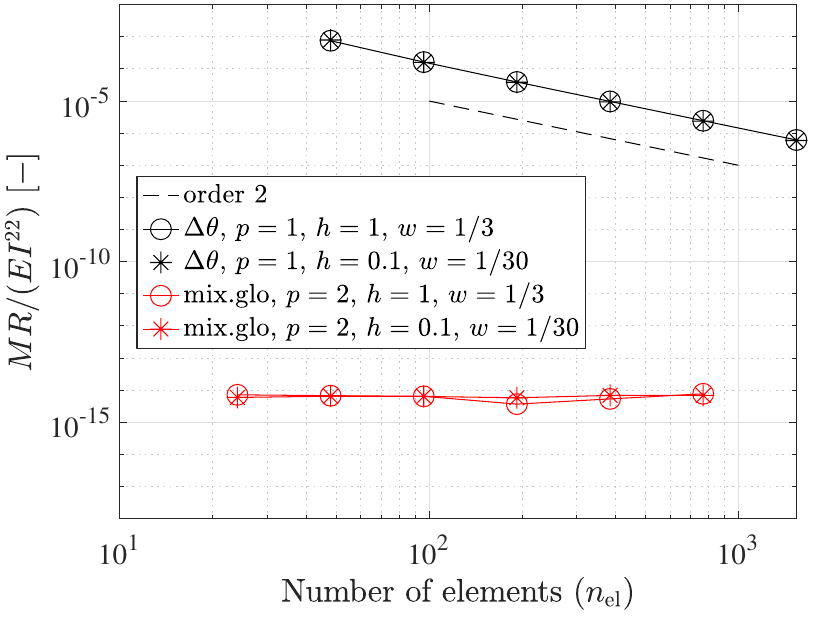}
		\caption{Moment at $\bar \theta=2\pi$}
		\label{twist_ring_mnt_tht_2pi}			
	\end{subfigure}
	\caption{Twisting of an elastic ring: path-independence test. Comparison of the applied moments at $\bar\theta=\pi$ and $\bar\theta=2\pi$. \textcolor{black}{In [mix.glo], we use $m_1=m_2=2$, $m_4=4$, and $E_{33}$ has not been enriched.}}
	\label{twist_ring_conv_mnt_pi_2pi}
\end{figure}

We consider a closed circular ring twisted by a prescribed rotation. Fig.\,\ref{twist_ring_init_config} shows the undeformed configuration and boundary conditions. The initial ring has a radius $R=20\,\mathrm{m}$, and a rectangular cross-section having width $w$ and height $h$. We prescribe a rotation of the cross-sections at the material points of initial coordinates $\left[R,0,0\right]^\mathrm{T}$ (point A) and $\left[-R,0,0\right]^\mathrm{T}$ (point B) by the angle $\bar \theta$ with respect to the $X$-axis in opposite directions. We apply the prescribed rotation angle $\bar \theta = 2\pi$ [rad] in total, such that the configuration finally goes back to the initial (undeformed) one, if no impenetrability condition is imposed, see Fig.\,\ref{twist_ring_deformed} for the deformed configurations at different values of $\bar \theta$. \textcolor{black}{For $\bar \theta=\pi$, the ring is folded into a smaller circular ring. Without axial strain, the radius of the deformed ring is $R/3$, which was first theoretically investigated by \citet{yoshiaki1992elastic}. With axial strain, the deformed configuration is still a folded ring, but has a different radius $r^*$ such that $2\pi\cdot{r}^*=\ell/3$, where $\ell$ denotes the current total length.} In the present paper, this example aims at the following:
\begin{itemize}
	\item Verify the capability of the present beam formulation to follow the correct equilibrium path of the deformation by correcting the torsional stiffness through the enrichment of higher-order cross-sectional strains,\\
	\item \textcolor{black}{Investigate the effect of axial strain coupled with bending strain in the present beam formulation. The coupled axial strain changes the total length $\ell$, and subsequently the radius of the deformed ring becomes $r^*={\ell}/(6\pi)$ at $\bar\theta=\pi$. This additional strain should quadratically decrease with increasing slenderness ratio \citep[Section 6.2]{choi2021isogeometric}. In order to show this, we consider the following two different cases with cross-sectional dimensions}:
	\begin{itemize}
		\item Case 1: $w=1/3\,\mathrm{m}$, and $h=1\,\mathrm{m}$,
		\item Case 2: $w=1/30\,\mathrm{m}$, and $h=1/10\,\mathrm{m}$.
	\end{itemize}
	Note that, here, we always have a ratio $h/w=3$, which is known to produce a deformed configuration of a smaller ring by the prescribe rotation ${\bar \theta}=\pi$ \citep{yoshiaki1992elastic}.\\
	\item Verify the path-independence of the present beam formulation.	
\end{itemize}
A St.\,Venant-Kirchhoff type isotropic hyperelastic material is considered, with Young's modulus $E=21\,\mathrm{MPa}$, and Poisson's ratio $\nu=0.3$. In order to account for the decreased torsional stiffness due to the cross-sectional warping, we employ the following approaches:
\begin{itemize}
	\item In the beam formulation [$\Delta \theta$], we correct the torsional stiffness by replacing the polar moment of inertia with $I^*_\mathrm{p}=K$, see Appendix \ref{app_asol_tors_stiff}. Note that $K$ is much smaller than the polar moment of inertia, $I_\mathrm{p}={I^{11}}+{I^{22}}$. For example, in case of $w=1/3\,\mathrm{m}$ and $h=1\,\mathrm{m}$, we have $I^{11}\approx{3.086\times10^{-3}}\,[\mathrm{m}^4]$ and $I^{22}\approx{2.778\times10^{-2}}\,[\mathrm{m}^4]$, which gives the polar moment of inertia much larger than the corrected value $K=9.753\times10^{-3}\,[\mathrm{m}^4]$. \textcolor{black}{This approach refers to \citet{romero2004interpolation}, \citet{meier2014objective}, and references therein.}\\
	\item For the present beam formulation, instead of directly adjusting the cross-sectional parameters, we enrich the polynomial basis functions up to degree $m_4$ for the transverse shear strains $E_{13}$ and $E_{23}$, and $m_1(=m_2)$ for the in-plane strains $E_{11}$, $E_{22}$, and $E_{12}$. Note that the axial normal strain component $E_{33}$ is not enriched, see Observation\,\ref{rem_num_ex_enrich_e33}.
\end{itemize}
In Fig.\,\ref{twist_ring_theta-mnt-curve}, we compare the equilibrium paths from using the formulation [$\Delta \theta$] and the present one [mix.glo]. The black solid curves represent the result from [$\Delta \theta$] using the corrected torsional stiffness with $I^*_\mathrm{p}=K$. Note that, without this correction (i.e., with $I^*_\mathrm{p}=I_\mathrm{p}$), the torsional stiffness is significantly overestimated, see the much larger initial slope of the black solid curves compared to the others. This eventually leads to the divergence in the equilibrium iteration under the prescribed rotation. \textcolor{black}{The black dashed curves represent the equilibrium paths with the corrected torsional stiffness.} In the results from the present beam formulation [mix.glo], we first enrich the transverse normal strain components $E_{11}$ and $E_{22}$ by the maximum degree of basis $m_1=2$ in order to obtain the correct bending stiffness due to the nonzero Poisson's ratio. Then, as we increase the maximum degree of enriched basis for the transverse shear strain $m_4$, the equilibrium path approaches the result from using $m_4=8$ (red curve). 
\noindent \textcolor{black}{We have an analytical solution of the deformed center axis such that it becomes a folded circular ring of radius $r_\mathrm{ref}=R/3$ at $\bar \theta = \pi$ \citep{yoshiaki1992elastic}. However, it turns out that additional axial compression can be induced by bending \citep[Sec.6.2.1]{choi2021isogeometric}. In order to illustrate this, in Fig.\,\ref{twist_ring_deformed_length_conv}, we compare the approximated total length of the center axis at $\bar \theta = \pi$ between the results from the formulations [$\Delta\theta$] and [mix.glo]. It is noticeable that, in the present beam solution (red curve), the beam is axially compressed, in contrast to the result from using [$\Delta\theta$]. Further, it is seen that the amount of compressive strain decreases by $10^{-2}$ times, with increasing the slenderness ratio by $10$ times from Case 1 to Case 2. This means that as expected, the decrease of the coupled axial strain is in quadratic order with respect to the slenderness ratio. Based on this, we introduce a corrected reference solution for the deformed center axis at $\bar\theta=\pi$, such that it is still a folded circular ring, but with radius $r^{*}_\mathrm{ref}=\ell^h/(6\pi)$, where the (approximated) total length $\ell^h$ is given by an overkill solution from using the present beam formulation, [mix.glo] with $p=4$.} From this, we can define a relative $L^2$-difference of the center axis displacement field, as 
\begin{align}
	{e_\varphi } \coloneqq \sqrt {\dfrac{{\int_0^L {{{\left( {{r^h} - {r^*_{{\rm{ref}}}}} \right)}^2}{\rm{d}}s} }}{{\int_0^L {{r^*_{{\rm{ref}}}}^2\,{\rm{d}}s} }}},\,\,\mathrm{with}\,\,{r^h} \coloneqq \left\| {{{\boldsymbol{\varphi }}^h} - {{\boldsymbol{c}}^*_{{\rm{ref}}}}} \right\|,
\end{align}
where the reference center of the circle is given by $\boldsymbol{c}^*_\mathrm{ref}\coloneqq{[-R+r^*_\mathrm{ref},0,0]^\mathrm{T}}$.
\noindent Fig.\,\ref{twist_ring_deformed_rad_conv} shows the convergence of the relative differences from two different beam formulations, $[\Delta\theta]$ with the degree of basis $p=1$, and the present one [mix.glo] with degree of basis $p=2,3,4$. \textcolor{black}{It is seen that in both cases 1 and 2, the present beam solution from using [mix.glo] converges toward the reference one (a folded circular ring with radius $r^*_\mathrm{ref}$), with increasing number of elements. In contrast, the beam formulation $[\Delta\theta]$ yields a different solution, which deviates from the reference one. This deviation is due to the inability to represent the axial strain coupled with bending, and it should decrease quadratically with increasing slenderness ratio. It is also noticeable that the results from [mix.glo] exhibit optimal convergence rates ($p+1$) or sometimes even faster, which can be clearly seen in both cases. In the results from using the beam formulation [mix.loc-sr], deteriorated convergence is observed due to numerical instability, except for the results from using $p=3$. In Table\,\ref{twist_nlstep_mnt_conv_test}, we compared the total number of load steps and equilibrium iterations between [mix.glo] and [mix.loc-sr]. The present formulation [mix.glo] shows stability in all cases.}
\noindent In \citet{yoshiaki1992elastic}, it is shown that the deformed smaller ring at $\bar \theta=\pi$ can stay without any external load. Further, due to the symmetry, the configuration should go back to the initial (undeformed) one at ${\bar \theta}=2\pi$ by either loading or unloading from the smaller ring at ${\bar \theta}=\pi$. That is, if a formulation is path-independent, the applied moment $M$ should vanish at both $\bar\theta=\pi$ and $2\pi$. Figs.\,\ref{twist_ring_mnt_tht_pi} and \ref{twist_ring_mnt_tht_2pi} plot the convergence of the moments $M$ at ${\bar\theta}=\pi$ and $2\pi$, respectively, with increasing number of elements. It is seen that, for the formulation [$\Delta\theta$], the moments do not vanish, which implies that the formulation suffers from path-dependence. Note that this error decreases in quadratic order, which is consistent with the convergence rate of the displacement field in Fig.\,\ref{twist_ring_deformed_rad_conv}. In contrast, from using the present beam formulation, the moments vanish up to machine precision for any number of elements. 

\noindent\textcolor{black}{The next two examples deal with nonlinear transient dynamics to verify the present time-stepping scheme.}
\subsection{Dynamics of a flying beam}
\begin{figure}[H]
	\centering
	\includegraphics[width=0.55\linewidth]{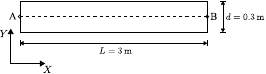}
	\caption{{Flying beam: Initial geometry.}}
	\label{rot_beam_init_geom}	
\end{figure}

\begin{figure}[H]
	\centering
	\begin{subfigure}[b]{0.42\textwidth}\centering
		\includegraphics[width=\linewidth]{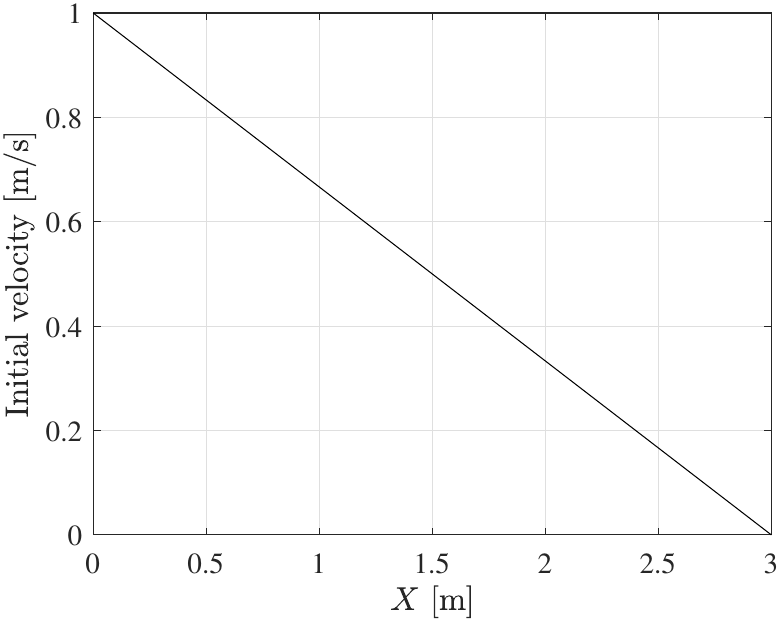}
		\caption{$X$--\,component (${v}^1_{\varphi0}$)}
		\label{rot_beam_init_vel_x}			
	\end{subfigure}
	\begin{subfigure}[b]{0.42\textwidth}\centering
		\includegraphics[width=\linewidth]{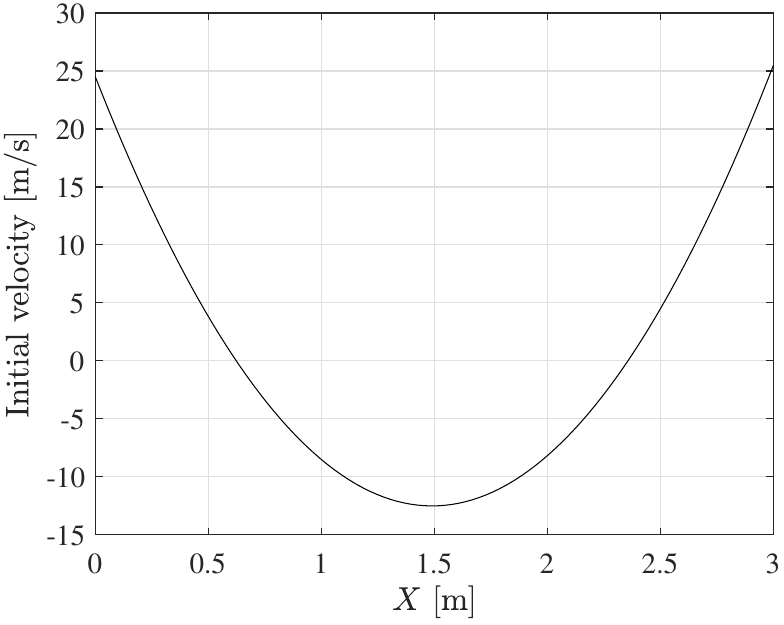}
		\caption{$Y$--\,component (${v}^2_{\varphi0}$)}
		\label{rot_beam_init_vel_y}			
	\end{subfigure}
	\caption{Flying beam: Given initial velocity of the center axis, ${\boldsymbol{v}}_{\varphi0}={{v}}^i_{{\varphi0}}{\boldsymbol{e}}_i$, where ${{v}}^3_{{\varphi0}}=0$.}
	\label{rot_beam_init_vel}
\end{figure}

\begin{figure}[H]
	\centering
	\includegraphics[width=0.65\linewidth]{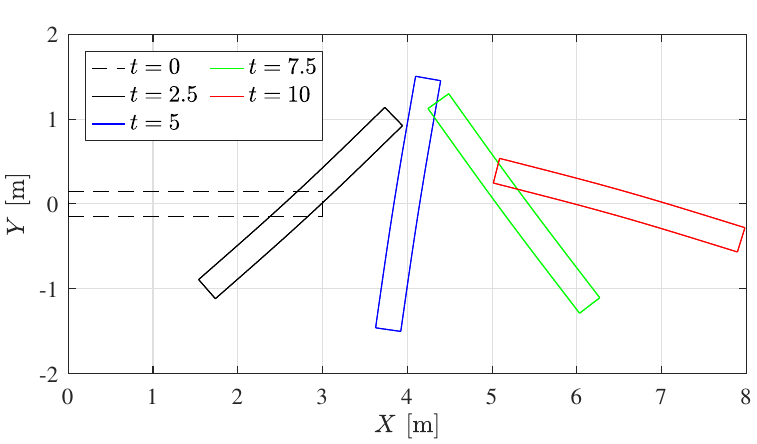}
	\caption{{Flying beam: Deformed configurations at the selected time $t\,[\mathrm{s}]$. \textcolor{black}{The result is from using the present EMC scheme, with the time step size $\Delta t=0.1\,\mathrm{s}$.}}}
	\label{rot_beam_deformed_geom}	
\end{figure}

\begin{figure}[H]
	\centering
	\begin{subfigure}[b]{0.45\textwidth}\centering
		\includegraphics[width=\linewidth]{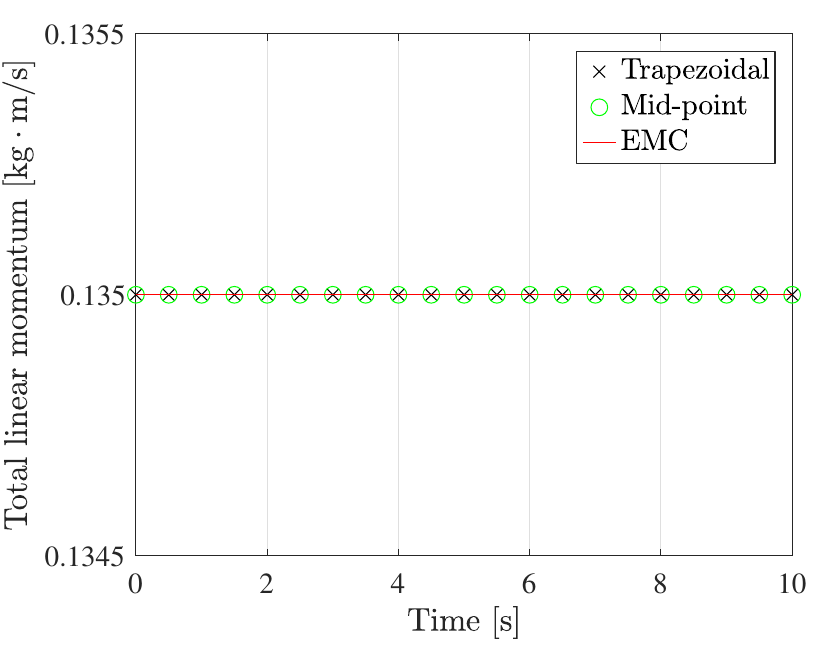}
		\caption{Total linear momentum}
		\label{rot_beam_lin_momentum}			
	\end{subfigure}
	\begin{subfigure}[b]{0.45\textwidth}\centering
		\includegraphics[width=\linewidth]{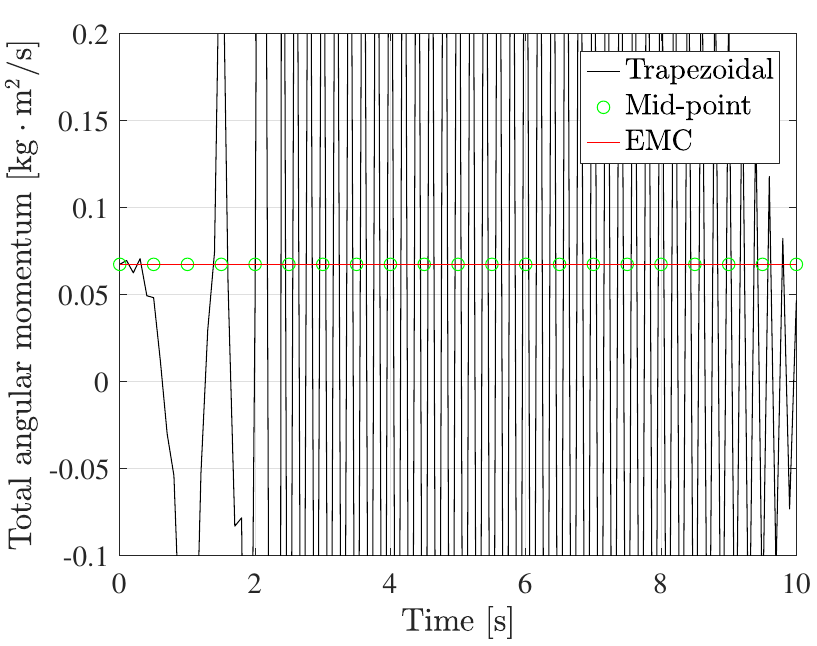}
		\caption{Total angular momentum}
		\label{rot_beam_ang_momentum}			
	\end{subfigure}
	\begin{subfigure}[b]{0.45\textwidth}\centering
		\includegraphics[width=\linewidth]{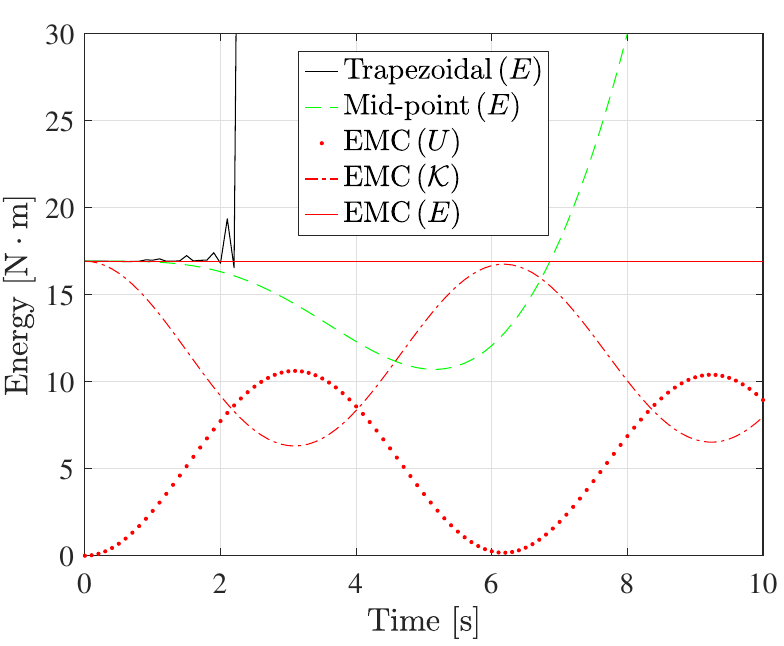}
		\caption{Total energy}		
		\label{rot_beam_tot_e}			
	\end{subfigure}
	\caption{Flying beam: Time history of the energy and momentum. In (a) and (b), the markers are plotted at every five time steps for clearer visualization. In all the results, we use the same time step size, $\Delta t=0.1\,\mathrm{s}$. In (c), the total energy $E$ is obtained from the summation of the kinetic energy ($\mathcal{K}$) and the strain energy ($U$).}
	\label{rot_beam_conserv}
\end{figure}

We consider an initially straight beam of length $L=3\,\mathrm{m}$ which has a square cross-section of dimension $d=0.3\,\mathrm{m}$, see Fig.\,\ref{rot_beam_init_geom}. A St.\,Venant-Kirchhoff type hyperelastic material with Young's modulus $E=21\,\mathrm{MPa}$ and Poisson's ratio $\nu=0.3$ is considered. \textcolor{black}{The initial mass density is $\rho_0=1\,\mathrm{kg}/\mathrm{m}^3$.} We apply an initial velocity condition in $X$-- and $Y$--\,directions along the initial center axis, as
\begin{subequations}
	\begin{align}
		{v^1_{\varphi0}} &=  - \frac{{{{\bar v}^1_\varphi}}}{L}\left( {X - L} \right),\label{rot_beam_init_vel_eq_x}\\
		{v^2_{\varphi0}} &= {\bar v^2_\varphi}\left\{ {\frac{X}{L}\left( {\frac{X}{L} - 1} \right) + \frac{1}{6}} \right\} + \frac{{{{\bar v}^1_\varphi}}}{L}\left( {X - \frac{L}{2}} \right),\label{rot_beam_init_vel_eq_y}
	\end{align}
\end{subequations}
respectively, with ${\bar v}^1_\varphi=1\,[\mathrm{m}/\mathrm{s}]$ and ${\bar v}^2_\varphi=150\,[\mathrm{m}/\mathrm{s}]$, see also Fig.\,\ref{rot_beam_init_vel}. The directors are assumed initially stationary, which means that all the material points in each cross-section have the same initial velocity. The initial velocity component in Eq.\,(\ref{rot_beam_init_vel_eq_x}) leads to a translation in $X$--\,direction. In Eq.\,(\ref{rot_beam_init_vel_eq_y}), the quadratic $Y$--\,directional velocity by the first term leads to the bending and transversal shear strains in the $X$--$Y$ plane, as well as the transverse normal strain due to the non-zero Poisson's ratio. The second term in Eq.\,(\ref{rot_beam_init_vel_eq_y}) is introduced to have a rotational motion in the $X$--$Y$ plane. Therefore, the total linear and angular momentum vectors may have non-zero components in $X$--, and $Z$--\,directions, respectively. Fig.\,\ref{rot_beam_deformed_geom} shows the deformed configurations at selected time instants.
\noindent In Fig.\,\ref{rot_beam_conserv}, we investigate the conservation of energy and momentum. It is seen that the standard trapezoidal rule conserves the linear momentum, but neither the angular momentum nor the total energy. For the given time step size $\Delta t=0.1\,\mathrm{s}$, the total energy blows-up. The total angular momentum initially has non-zero value due to the initial velocity condition, and is maintained by both the standard mid-point rule and the present EMC schemes. Initially, the total energy only has a kinetic part, and it is exactly conserved by the present EMC scheme. This implies the symmetry between the strain energy ($U$) and the kinetic energy ($\mathcal{K}$). In contrast, when using the mid-point rule, the total energy is not conserved, but blows-up. \textcolor{black}{Here, for all the beam solutions, we use [mix.glo] with $p=2$, $n_\mathrm{el}=10$, and $m_1=m_2=2$, $m_4=4$, with no enrichment for $E_{33}$.}
\subsection{Slotted ring with a rectangular cross-section}
\begin{figure}[H]
	\centering
	\includegraphics[width=0.35\linewidth]{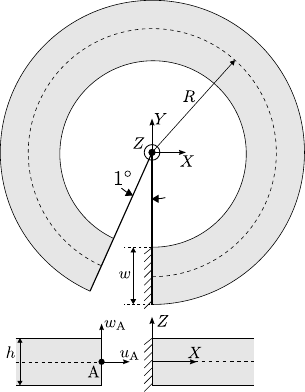}
	\caption{{Slotted ring: Undeformed configuration and boundary conditions. Dashed line represents the initial center axis.}}
	\label{slot_ring_init_config}	
\end{figure}

An initially circular ring of radius $R=1.3\,\mathrm{m}$ with a slot of angle $1^\circ$ is considered. It has a rectangular cross-section of height $h=0.4\,\mathrm{m}$ and width $w=0.3\,\mathrm{m}$, see Fig.\,\ref{slot_ring_init_config}. One end is fixed, where no cross-sectional deformation is allowed, and a distributed force is applied in $Z$-direction on the other end (cross-section at point A). In this example, we investigate two different types of loading:
\begin{itemize}
	\item Case 1 (elastostatics): A quasi-statically applied surface traction (force per unit undeformed area), ${\boldsymbol{T}}={F_0}\,\boldsymbol{e}_3$, where ${F_0}=5\times10^4\,\mathrm{N}/\mathrm{m}^2$,\\
	\item Case 2 (elastodynamics): A (piece-wise linearly) time-dependent surface traction, $\boldsymbol{T}(t)=F(t)\,\boldsymbol{e}_3$, where $F(t)$ (force per unit undeformed area) is given by
\end{itemize}
\begin{align}
	\label{slotted_ring_force_area}
	F(t) = \left\{ {\renewcommand{\arraystretch}{2.5}\begin{array}{*{20}{c}}
			{{F_0}\,t}/{T_0}&{{\rm{if}}\,\,0 \le t \le \frac{T}{2}},\\
			{{F_0}\left( {T - t} \right)/{T_0}}&\,{{\rm{if}}\,\,\frac{T}{2} < t \le T}.
	\end{array}} \right.
\end{align}
\textcolor{black}{Here, ${F_0}=5\times10^4\,\mathrm{N}/\mathrm{m}^2$, $T_0=1\,\mathrm{s}$ denotes a unit time, and $T$ denotes the terminal time.}
\subsubsection{Case 1: Comparison with the brick element solution}
\begin{figure}[H]
	\centering
	\begin{subfigure}[b]{0.35\textwidth}\centering
		\includegraphics[width=\linewidth]{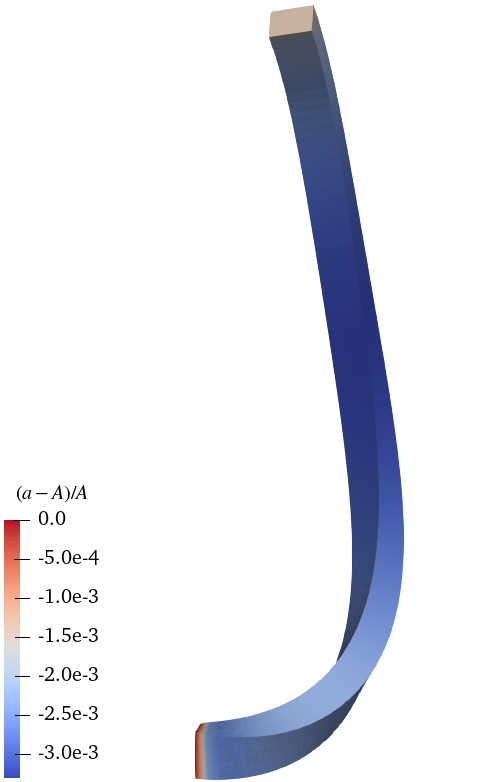}
		\caption{Brick element (32,448 DOFs)}		
		\label{deformed_brick}			
	\end{subfigure}
	\begin{subfigure}[b]{0.35\textwidth}\centering
		\includegraphics[width=\linewidth]{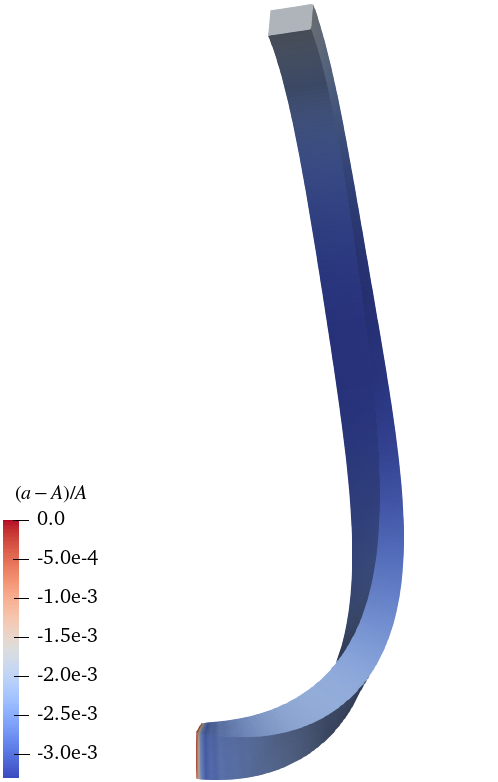}
		\caption{Beam element (6,561 DOFs)}		
		\label{deformed_beam}			
	\end{subfigure}
	\caption{Slotted ring: Comparison of the final deformed configurations from the reference brick and the present beam formulations. The contour represents the relative change of the cross-sectional area, where $A$ and $a$ denote the initial and current cross-sectional areas, respectively. For the brick solution, we use deg.=$(3,3,3)$, and $n_\mathrm{el}=160\times5\times5$. For the beam solution, we use [mix.glo] with $p=3$, $n_\mathrm{el}=160$, and ${m_1}={m_2}=2$, $m_4 = 4$, where $E_{33}$ has not been enriched.}
	\label{slotted_ring_deformed_config_ca_contour}
\end{figure}

\begin{figure}[H]
	\centering
	\begin{subfigure}[b]{0.4875\textwidth}\centering
		\includegraphics[width=\linewidth]{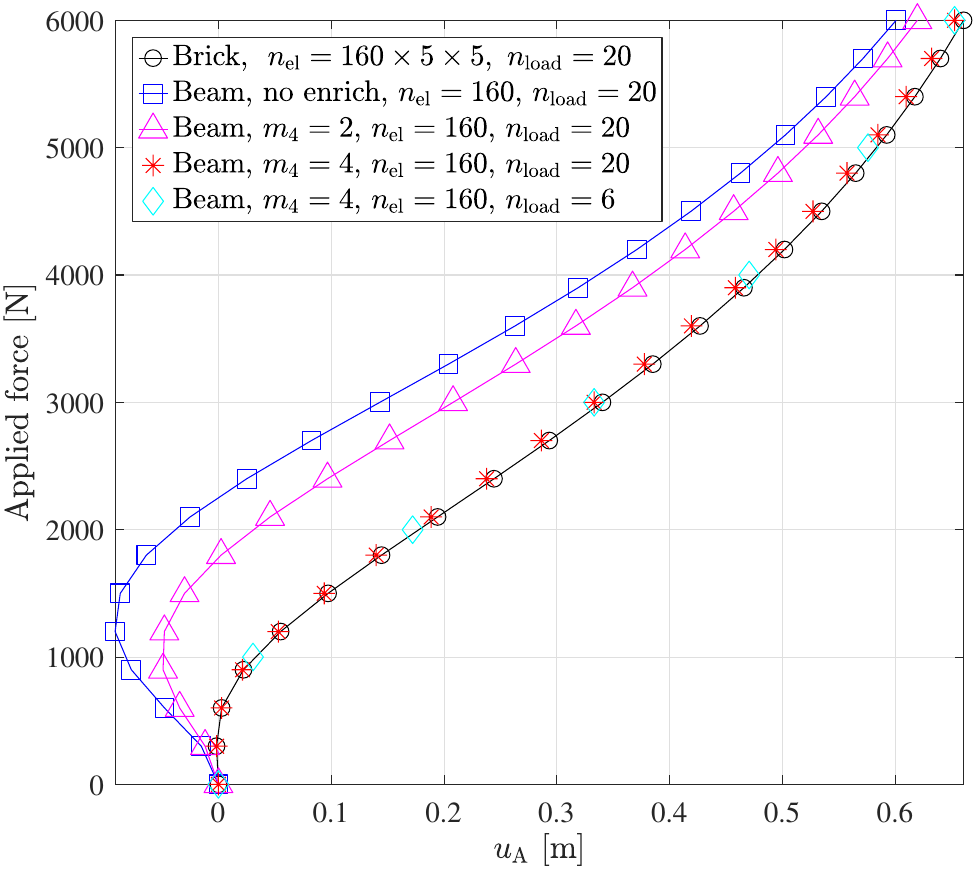}
		\caption{$X$--displacement}				
		\label{slot_ring_ux}			
	\end{subfigure}	
	\begin{subfigure}[b]{0.4875\textwidth}\centering
		\includegraphics[width=\linewidth]{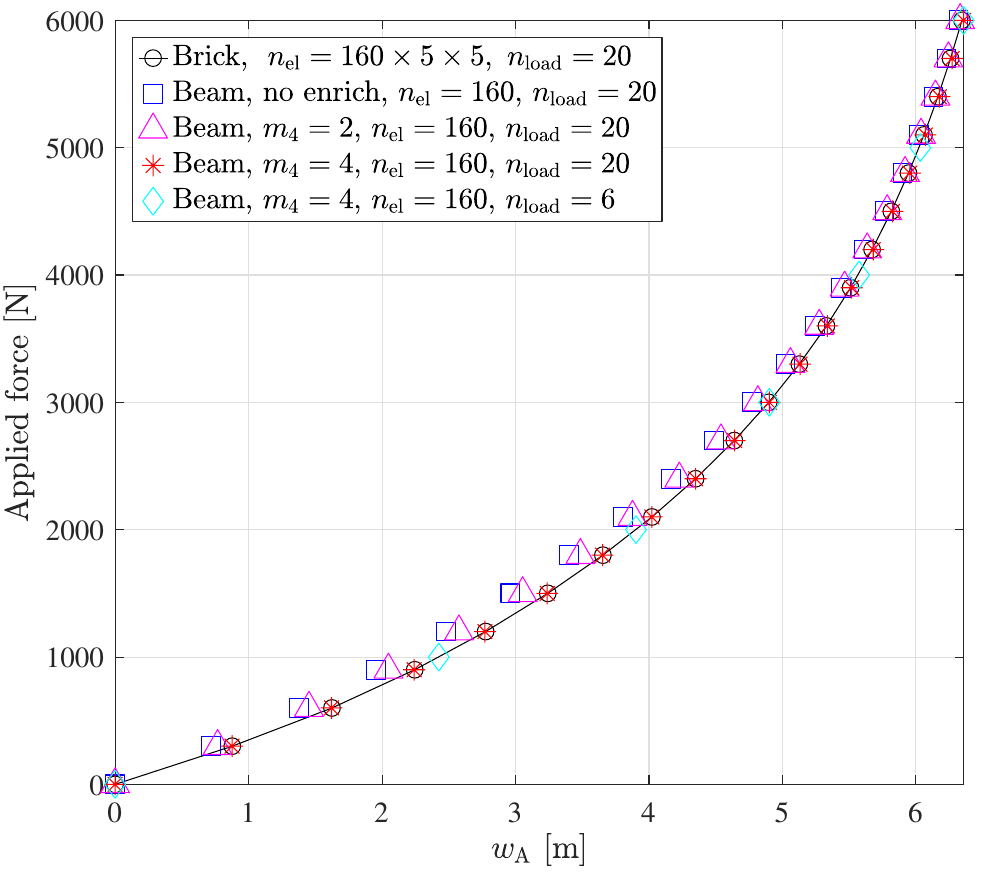}
		\caption{$Z$--displacement}			
		\label{slot_ring_uz}			
	\end{subfigure}	
	\caption{Slotted ring: Comparison of the $X$ and $Z$-directional displacements at the tip (point A) with increasing applied force. \textcolor{black}{In all the beam solutions, we use [mix.glo] with $p=3$, $m_1=m_2=2$, and $E_{33}$ has not been enriched.}}
	\label{slotted_ring_fd_curve}
\end{figure}

We consider a Neo-Hookean type isotropic hyperelastic material with Young's modulus $E=1.12\times{10^7}\,\mathrm{Pa}$, and Poisson's ratio $\nu=0.4$. For the given load, the ring undergoes a large deformation with torsion-induced cross-sectional warping. Here, we verify the corrected stiffness from using the present EAS method. In Fig.\,\ref{slotted_ring_deformed_config_ca_contour}, we compare the deformed configurations between the present beam and the brick element solutions, where it is seen that the cross-section contracts due to the positive Poisson's ratio, except for the fixed end. It is seen that the change of cross-sectional area can be properly represented by only two directors in the present beam formulation.
\noindent Figs.\,\ref{slot_ring_ux} and \ref{slot_ring_uz} show the $X$-- and $Z$--displacements of the ring at point A, which are denoted by $u_\mathrm{A}$ and $w_\mathrm{A}$, respectively. In the beam formulation, we compare different degrees of the polynomial basis $m_4=2,3,4$ for the transverse shear strains, $E_{13}$ and $E_{23}$. Note that we also use $m_1=m_2=2$, and the axial normal strain $E_{33}$ has not been enriched. It is seen that without the enrichment of the transverse shear strain, the beam solution (blue curve) deviates signifcantly from the brick solution (black curve). It is noticeable that with increasing degree $m_4$, the solution approaches the brick solution (magenta and red marks). Further, a much larger load step is allowed in the present beam formulation, compared to the displacement-based brick formulation, which can be seen by the results of the beam element using 6 uniform load increments (cyan mark). 
\subsubsection{Case 2: Verification of the consistency in the EMC scheme}
\begin{figure}[H]
	\centering
	\begin{subfigure}[b]{0.4875\textwidth}\centering
		\includegraphics[width=\linewidth]{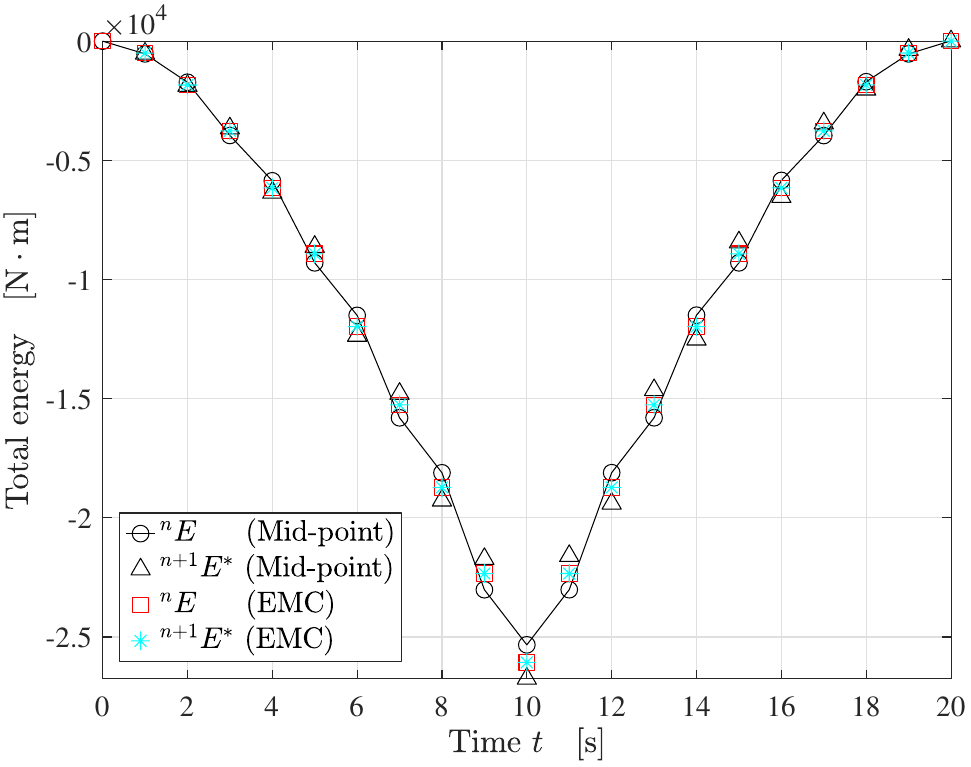}
		\caption{Total energy}				
		\label{slot_ring_hamt_svk}			
	\end{subfigure}	
	\begin{subfigure}[b]{0.4875\textwidth}\centering
		\includegraphics[width=\linewidth]{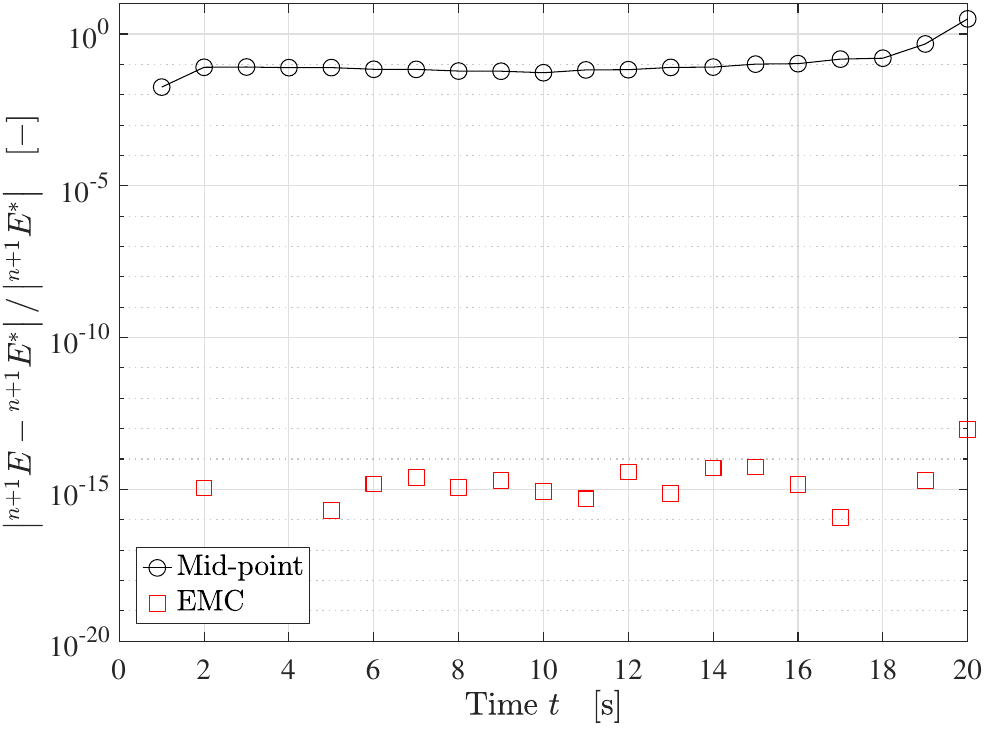}
		\caption{Relative difference}			
		\label{slot_ring_hamt_diff_svk}			
	\end{subfigure}	
	\caption{Slotted ring: Verification of the energy--consistency for a linear constitutive law (St.\,Venant-Kirchhoff material). Missing data points (red rectangles) in (b) are due to exact zero difference in the numerator.}
	\label{slotted_ring_time_hamt_curve_SVK}
\end{figure}

\begin{figure}[H]
	\centering
	\begin{subfigure}[b]{0.4875\textwidth}\centering
		\includegraphics[width=\linewidth]{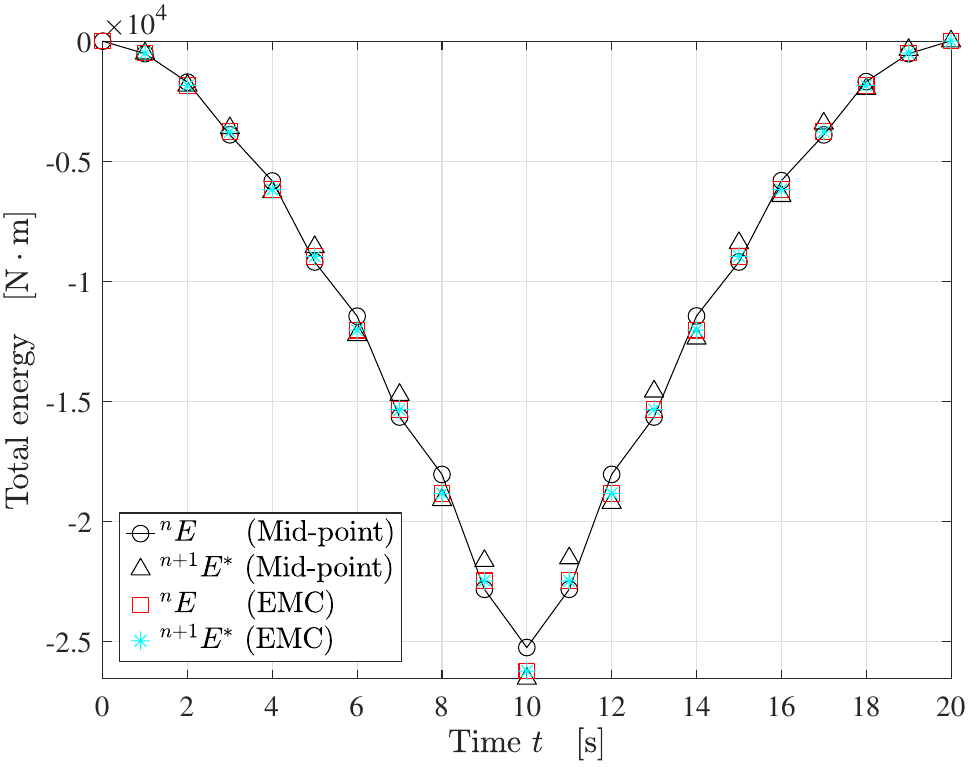}
		\caption{Total energy}				
		\label{slot_ring_hamt_nh}			
	\end{subfigure}	
	\begin{subfigure}[b]{0.4875\textwidth}\centering
		\includegraphics[width=\linewidth]{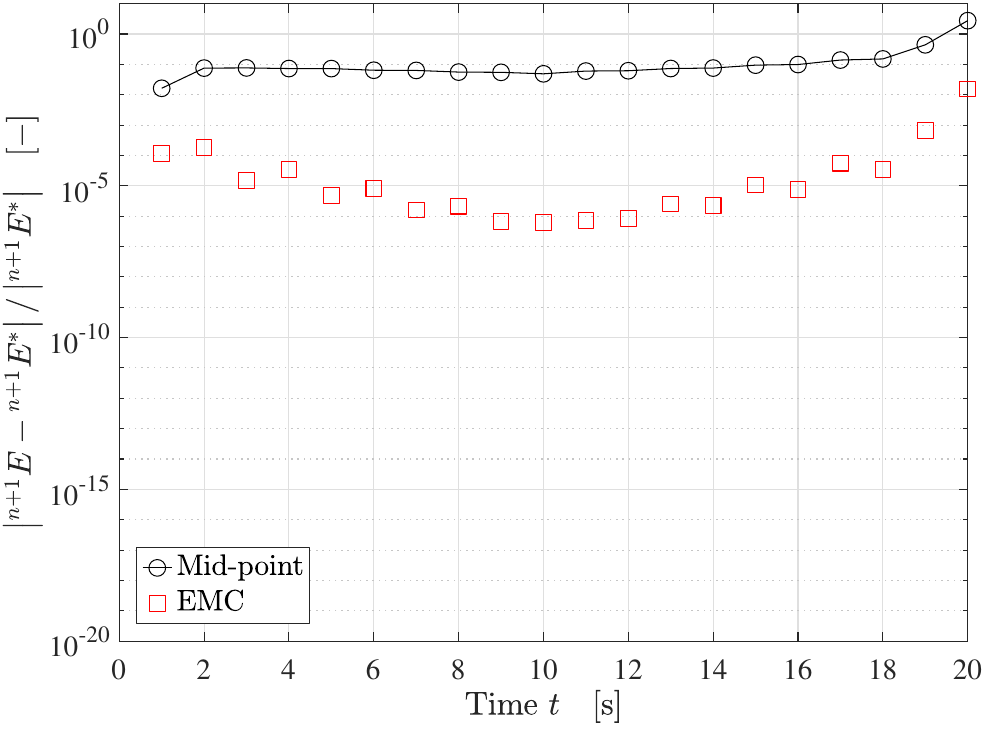}
		\caption{Relative difference}			
		\label{slot_ring_hamt_diff_nh}			
	\end{subfigure}	
	\caption{Slotted ring: Verification of the energy--consistency for a nonlinear constitutive law (Neo-Hookean material).}
	\label{slotted_ring_time_hamt_curve_NH}
\end{figure}

In the second case, we investigate a transient dynamic response of the slotted ring, under the time-dependent load given by Eq.\,(\ref{slotted_ring_force_area}). The total simulation time is $T=20\,\mathrm{s}$. Figs.\,\ref{slotted_ring_time_hamt_curve_SVK} and \ref{slotted_ring_time_hamt_curve_NH} investigate whether a time-stepping scheme satisfies the following identity
\begin{align}
	\label{identity_recall_energy_consistency}
	{}^{n+1}E={}^{n+1}E^*,\,\,n=0,1,\cdots,
\end{align}
which is recalled from Eq.\,(\ref{total_energy_consistency_pie_lin_f}). For constitutive laws, we consider two different models, St.\,Venant-Kirchhoff (linear) and Neo-Hookean (nonlinear) types, with the same values of the parameters $E$ and $\nu$ as in the case 1. \textcolor{black}{The initial mass density is $\rho_0=1\,\mathrm{kg}/\mathrm{m}^3$.} From using the standard mid-point rule, it is seen that the total energy ${}^{n+1}E$ deviates from ${}^{n+1}E^{*}$ for both linear and nonlinear constitutive laws. In contrast, the EMC scheme satisfies the identity exactly (up to machine precision) for the linear constitutive law, see Fig.\,\ref{slot_ring_hamt_diff_svk}. For the nonlinear constitutive law, the EMC scheme has also lost the consistency, but the relative difference is much smaller than that from the mid-point rule, see Fig.\,\ref{slot_ring_hamt_diff_nh}.
\section{Conclusions} 
\label{conclusions_sec}
We present a stable and efficient implicit \textcolor{black}{isogeometric} finite element method for geometrically and materially nonlinear beams in transient dynamics. The superior numerical stability of the developed method comes from the following aspects:
\begin{itemize}
	\item Objectivity and path-independence of the isogeoemtric finite element formulation that directly approximates the director field, which turns out to hold for any degree $p_\mathrm{d}$, which is verified in several numerical examples. \textcolor{black}{This formulation is based on the approximation of the initial director field, which is, in the present paper, denoted by [$D$-disc.].}
	\item An energy--momentum consistent (EMC) time-stepping scheme, which gives exact energy conservation property for linear constitutive laws.
	\item Employing the independent DOFs for the stress resultants in the mixed finite element formulation based on the Hu-Washizu variational principle alleviates the overestimation of the internal force, which improves the stability of the equilibrium iteration for larger time (load) steps in the Newton-Raphson method.	
	\item \textcolor{black}{In the isogeometric approach, we have observed numerical instabilities of the local approach [mix.loc-sr] in some cases. The present global approach [mix.glo] shows numerical stability in all those cases. Further, compared to the local approach, the present global approach requires much less DOFs for the physical stress resultants and strains due to the higher-order continuity of the B-spline basis functions.}
\end{itemize}
Further, the developed method provides superior computational efficiency, compared to brick element formulations, due to the following aspects:
\begin{itemize}
	\item The present 9-DOF formulation with extensible (unconstrained) directors allows for an efficient representation of a constant in-plane cross-sectional strain, e.g., the change of cross-sectional area. 
	\item The higher order cross-sectional strains like torsion-induced warping have been further enriched by an EAS method. Due to the allowed inter-element discontinuity in the enhanced strain field, one can eliminate those additional DOFs from the global system of linear equations. 
\end{itemize}
One can extend the present work in the following aspects:
\begin{itemize}
	\item \textcolor{black}{In the inertia terms,} we have considered the cross-sectional deformation only up to the constant strains, represented by two extensible directors.
	\item We have limited our scope of investigation to convex-shaped cross-sections like a rectangular one.
	\item The present EAS method corrects stiffness very effectively; however, it is seen that the enrichment of in-plane cross-sectional strains is limited.
\end{itemize}
\backmatter

\bmhead{Acknowledgements}
	M.-J. Choi was supported by the Alexander von Humboldt Foundation (postdoctoral research fellowship) in Germany, and the Deutsche Forschungsgemeinschaft (DFG, German Research Foundation) - Project number 523829370.

\bmhead{Author contributions} \textbf{Myung-Jin Choi:} Writing - review \& editing, writing - original draft, investigation, conceptualization, methodology, software, validation, visualization, formal analysis, data curation, resources. \textbf{Sven Klinkel:} Writing - review \& editing, investigation, conceptualization, resources, supervision. \textbf{Simon Klarmann:} Writing - review \& editing, investigation, conceptualization. \textbf{Roger A. Sauer:} Writing - review \& editing, investigation, conceptualization, supervision.
	
\section*{Declarations}
The authors declare that they have no conflict of interest.

\begin{appendices}

\section{Beam formulation}
\subsection{Objectivity of the geometric strain in continuous form}
\label{verif_objectivity_strn_cont}
The \textit{invariance} of the beam strain components under a superposed rigid body motion can be analytically verified. We superpose an arbitrary constant rigid body rotation $\boldsymbol{\Lambda}_\theta\in{\rm{SO}(3)}$ and translation ${\boldsymbol{c}}_\varphi\in\mathbb{R}^3$ to the current configuration $\boldsymbol{y}$, as \citep{crisfield1999objectivity}
\begin{subequations}
	\label{superposed_config_y}
	\begin{align}
		{{\boldsymbol{\varphi }}^*} &= {{\boldsymbol{\Lambda }}_{\theta}}\left( {{\boldsymbol{\varphi }} + {{\boldsymbol{c}}_{\varphi}}} \right),\label{eq_rigid_mo_phi}\\
		{{\boldsymbol{d}}_\alpha^*} &= {{\boldsymbol{\Lambda }}_{\theta}}\,{{\boldsymbol{d }}_\alpha},\,\,\alpha\in\left\{1,2\right\},\label{eq_rigid_mo_d}
	\end{align}
\end{subequations}
where $\mathrm{SO(3)}$ denotes a three-dimensional rotation group, defined as
\begin{align}
	{\rm{SO}}(3) \coloneqq \left\{ {\left. {{\boldsymbol{\Lambda }} \in {{\mathbb{R}}^{3 \times 3}}} \right|{{\boldsymbol{\Lambda }}^{\rm{T}}}{\boldsymbol{\Lambda }} = {\bf{1}},\,{\rm{and}}\,\,{\rm{det}}\,{\boldsymbol{\Lambda }} = 1} \right\}.\!
\end{align}
\textcolor{black}{By inserting Eqs.\,(\ref{eq_rigid_mo_phi}) and (\ref{eq_rigid_mo_d}) into Eqs.\,(\ref{axial_strn_comp})--(\ref{strn_high_bend})}, and using ${{\boldsymbol{\Lambda }}_{\theta,s}} = {\bf{0}}$, ${{\boldsymbol{c}}_{\varphi {\rm{,}}s}} = {\bf{0}}$, and ${{\boldsymbol{\Lambda }}^\mathrm{T}_{\theta}}{{\boldsymbol{\Lambda }}_{\theta}} = {\boldsymbol{1}}$, we have
\begin{equation}
	{\boldsymbol{\varepsilon }}({{\boldsymbol{y}}^ * }) = {\boldsymbol{\varepsilon }}({\boldsymbol{y}}).
\end{equation} 
\subsection{Objectivity of the geometric strain in discrete form}
\label{verif_objectivity_strn_disc}
We investigate the invariance of the strain measures approximated by Eqs.\,(\ref{approx_caxis_pos}) and (\ref{tot_dirv_disp}) under superposed rigid body translation and rotation. We superpose an arbitrary constant rigid body rotation $\boldsymbol{\Lambda}_\theta\in{\rm{SO}(3)}$ and translation ${\boldsymbol{c}}_{{\varphi}I}\in\mathbb{R}^3$ to the current control coefficients $\boldsymbol{\varphi}_I$ and $\boldsymbol{d}_{\alpha{I}}$, as	
\begin{subequations}
\begin{align}
	{{\boldsymbol{\varphi }}^{h * }}(\xi ) &= \sum\limits_{I = 1}^{{n_{{\rm{cp}}}}} {N_I^p(\xi ){{\boldsymbol{\Lambda }}_\theta }\left( {{{\boldsymbol{\varphi }}_I} + {{\boldsymbol{c}}_\varphi }} \right)} \nonumber\\
	&= {{\boldsymbol{\Lambda }}_\theta }\left(\sum\limits_{I = 1}^{{n_{{\rm{cp}}}}} {N_I^p(\xi ){{\boldsymbol{\varphi }}_I}}  + {{\boldsymbol{c}}_\varphi }\right)\nonumber\\
	&= {{\boldsymbol{\Lambda }}_\theta }\,\left({{\boldsymbol{\varphi }}^h}(\xi ) + {{\boldsymbol{c}}_\varphi }\right), \label{obj_verif_ptb_caxis_approx}
\end{align}
and
\begin{align}
	{\boldsymbol{d}}_\alpha ^{h * }(\xi ) &= \sum\limits_{I = 1}^{{n^\mathrm{d}_{{\rm{cp}}}}} {N_I^{p_\mathrm{d}}(\xi ){{\boldsymbol{\Lambda }}_\theta }{{\bf{d}}_{\alpha I}}} \nonumber\\
	&= {{\boldsymbol{\Lambda }}_\theta }\sum\limits_{I = 1}^{{n^\mathrm{d}_{{\rm{cp}}}}} {N_I^{p_\mathrm{d}}(\xi ){{\bf{d}}_{\alpha I}}} \nonumber\\
	&= {{\boldsymbol{\Lambda }}_\theta }\,{\boldsymbol{d}}_\alpha ^h(\xi ), \label{obj_verif_ptb_dir_approx}
\end{align}
\end{subequations}
$\alpha\in\left\{1,2\right\}$, where we utilize the partition of unity property of the NURBS basis functions. Therefore, \textcolor{black}{substituting Eqs.\,(\ref{obj_verif_ptb_caxis_approx}) and (\ref{obj_verif_ptb_dir_approx}) into Eqs.\,(\ref{axial_strn_comp})--(\ref{strn_high_bend})}, and using ${{\boldsymbol{\Lambda }}_{{\theta},s}} = {\boldsymbol{0}}$, ${{\boldsymbol{c}}_{\varphi {\rm{,}}s}} = {\boldsymbol{0}}$, and ${\boldsymbol{\Lambda }}_\theta^{\rm{T}}{{\boldsymbol{\Lambda }}_\theta} = {\boldsymbol{1}}$, the beam strain at the superposed configuration ${{\boldsymbol{y}}^{h*}}$, we have
\begin{align}
	{\boldsymbol{\varepsilon }}\left({{\boldsymbol{
				y}}^{h*}}\right) = {\boldsymbol{\varepsilon }}({{\boldsymbol{y}}^h}),
\end{align}
which represents the frame-invariance (objectivity) of the approximated beam strain.
\subsection{Construction of warping basis functions}
\label{app_construct_warp_basis}
\textcolor{black}{Here we discuss the determination of the polynomial basis ${\bf{w}}_1$, ${\bf{w}}_2$, ${\bf{w}}_3$, and ${\bf{w}}_4$ in Eq.\,(\ref{wg_basis_matrix_gamma}). Note that this is a generalization of the formulation of \citet{wackerfuss2009mixed} in order to account for the higher order strains in Eq.\,(\ref{E_vgt_decomp_comp}) from the beam kinematics.
\subsubsection{Orthogonality condition}
\label{app_ortho_cond}
Here we construct a complete set of polynomials \textcolor{black}{from degree $\bar m + 1$ to} $M$ that are orthogonal to an arbitrary polynomial of degree $\bar m$. We first introduce the following arrays \citep{wackerfuss2011nonlinear}
\begin{align}
	{\boldsymbol{\pi} _0} &\coloneqq \left[ {{P_{(0,0)}}} \right],\nonumber\\	
	{\boldsymbol{\pi} _1} &\coloneqq \left[ {{P_{(1,0)}},{P_{(0,1)}}} \right],\nonumber\\
	&\quad\vdots\nonumber\\
	{\boldsymbol{\pi}_M} &\coloneqq \left[ {{P_{(M,0)}},{P_{(M-1,1)}},\cdots,{P_{(1,M-1)}},{P_{(0,M)}}} \right],
\end{align}
with ${P_{(n,m)}} \coloneqq {\left( {{\zeta ^1}} \right)^{\!n}}{\left( {{\zeta ^2}} \right)^{\!m}}$, $(n,m)\in{\Bbb{S}}_M\coloneqq{\Bbb{P}_0}\cup{\Bbb{P}_1}\cdots\cup{\Bbb{P}_{M}}$, where $M$ denotes the maximum degree, and we define
\begin{align*}
{\Bbb{P}}_M\coloneqq\Bigl\{\Bigl.(n,m)\in{\Bbb{N}_{0}\times\Bbb{N}_{0}} \Bigr|\,n+m=M\Bigr\}.
\end{align*}
Here, $\Bbb{N}_{0}$ denotes the set of positive integers, and the number of elements in ${\Bbb{P}}_M$ is $\left|{\Bbb{P}}_M\right|=M+1$. Further, we introduce a modified polynomial basis from degree $\bar m + 1$, as 
\begin{align}
	{\boldsymbol{\pi} ^*_{\bar m + 1}} &\coloneqq \left[ {{P^*_{(\bar m + 1,0)}},{P^*_{(\bar m,1)}},\cdots,{P^*_{(1,\bar m)}},{P^*_{(0,\bar m + 1)}}} \right],\nonumber\\	
	{\boldsymbol{\pi} ^*_{\bar m + 2}} &\coloneqq \left[ {{P^*_{(\bar m + 2,0)}},{P^*_{(\bar m + 1,1)}},\cdots,{P^*_{(1,\bar m + 1)}},{P^*_{(0,\bar m + 2)}}} \right],\nonumber\\
	&\quad\vdots\nonumber\\
	{\boldsymbol{\pi}^*_M} &\coloneqq \left[ {{P^*_{(M,0)}},{P^*_{(M-1,1)}},\cdots,{P^*_{(1,M-1)}},{P^*_{(0,M)}}} \right].
\end{align}
Here we define, for each $(n,m)\in{\Bbb{S}}_M - {\Bbb{S}}_{\bar m}={\Bbb{P}_{\bar m + 1}}\cup{\Bbb{P}_{\bar m + 2}}\cdots\cup{\Bbb{P}_{M}}$, 
\begin{align}
	\label{mod_P_poly_N_n}	
	P^*_{(n,m)} &\coloneqq {P_{(n,m)}} + \sum\limits_{k = 0}^{\bar m} {{\boldsymbol{\pi }}_k{\boldsymbol{\beta }}_{(n,m)}^k},
\end{align}
with the coefficients ${\boldsymbol{\beta }}_{(n,m)}^k \in {\Bbb{R}^{\left| {{\Bbb{P}_k}} \right|}}$,
\begin{align}
	\label{coeff_beta_mod_poly_basis}
	{\boldsymbol{\beta }}_{\left( {n,m} \right)}^0 &\coloneqq \left[ {\beta _{\left( {n,m} \right)}^{\left( {0,0} \right)}} \right]^{\!\mathrm{T}},\nonumber\\
	{\boldsymbol{\beta }}_{\left( {n,m} \right)}^1 &\coloneqq \left[ {\beta _{\left( {n,m} \right)}^{\left( {1,0} \right)},\beta _{\left( {n,m} \right)}^{\left( {0,1} \right)}} \right]^{\!\mathrm{T}},\nonumber\\
	&\quad\vdots\nonumber\\
	{\boldsymbol{\beta }}_{\left( {n,m} \right)}^{\bar m} &\coloneqq \left[ {\beta _{\left( {n,m} \right)}^{\left( {{\bar m},0} \right)},\beta _{\left( {n,m} \right)}^{\left( {{\bar m}-1,1} \right)},\cdots,\beta _{\left( {n,m} \right)}^{\left( {0,{\bar m}} \right)}} \right]^{\!\mathrm{T}},
\end{align}
such that the following orthogonality to arbitrary polynomials up to degree $\bar m$ is satisfied:
\begin{align}
	\label{ortho_pbasis_mod_deter_b}
	\int_\mathcal{A} P_{(p,q)}{P_{(n,m)}^{*}\,{\rm{d}}\mathcal{A}}  = 0,\quad\forall(p,q)\in{\Bbb{S}}_{\bar m},
\end{align}
where ${\Bbb{S}}_{\bar m}\coloneqq{\Bbb{P}_0}\cup{\Bbb{P}_1}\cdots\cup{\Bbb{P}_{\bar m}}$ whose number of elements is $\left|{{\Bbb{S}}_{\bar m}}\right|=(\bar m +1)(\bar m +2)/2$. Equation\,(\ref{mod_P_poly_N_n}) can be rewritten in the compact form
\begin{align}
	\label{mod_P_poly_N_n_mod_form}
	P_{(n,m)}^* = {P_{(n,m)}} + {\bf{\bar w}}{{\boldsymbol{\beta }}_{(n,m)}},
\end{align}
with ${\bf{\bar w}}\coloneqq{\left[ {{{\boldsymbol{\pi }}_0},{{\boldsymbol{\pi }}_1}, \cdots ,{{\boldsymbol{\pi }}_{\bar m}}} \right]}\in\Bbb{R}^{1\times{\left|{\Bbb{S}}_{\bar m}\right|}}$ and ${\boldsymbol{\beta }}_{(n,m)}\coloneqq\left[{\boldsymbol{\beta }}_{\left( {n,m} \right)}^{0\,\mathrm{T}},{\boldsymbol{\beta }}_{\left( {n,m} \right)}^{1\,\mathrm{T}},\cdots,{\boldsymbol{\beta }}_{\left( {n,m} \right)}^{{\bar m}\,\mathrm{T}}\right]^{\!\mathrm{T}}\in\Bbb{R}^{{\left|{\Bbb{S}}_{\bar m}\right|}}$.
Substituting Eq.\,(\ref{mod_P_poly_N_n_mod_form}) into Eq.\,(\ref{ortho_pbasis_mod_deter_b}), we obtain the following system of linear equations to determine $\boldsymbol{\beta}_{(n,m)}$ for each $(n,m)\in{\Bbb{S}}_{M}-{\Bbb{S}}_{\bar m}$:
\begin{align}
	\label{deter_beta_nm_lin}
	{{\bar{\Bbb{W}}}}\,\boldsymbol{\beta}_{(n,m)} =  - \int_\mathcal{A} {{{\bf{\bar w}}^{\rm{T}}}{P_{(n,m)}}\,{\rm{d}}\mathcal{A}},
\end{align}
where the system matrix is given by
\begin{align}
	{\bar{\Bbb{W}}} \coloneqq\int_\mathcal{A} {{\bf{\bar w}}^\mathrm{T}{\bf{\bar w}}\,{\rm{d}}\mathcal{A}}\in\Bbb{R}^{{\left|{\Bbb{S}}_{\bar m}\right|}\times{\left|{\Bbb{S}}_{\bar m}\right|}}.
\end{align}
\begin{remark}
	\label{rem_determin_beta_lin}
Since $\bar {\Bbb{W}}$ is invertible, the coefficients $\boldsymbol{\beta}_{(n,m)}$ for each $(n,m)\in{\Bbb{S}_M}-{\Bbb{S}_{\bar m}}$ can be uniquely determined by solving Eq.\,(\ref{deter_beta_nm_lin}). It should be noted that Eq.\,(\ref{deter_beta_nm_lin}) depends solely on the cross-section's initial shape and dimension. Therefore, this process to contruct a complete set of polynomial basis functions from degree $\bar m + 1$ to $M$, ${\bf{w}}\coloneqq{\left[{\boldsymbol{\pi}}^*_{\bar m + 1},{\boldsymbol{\pi}}^*_{\bar m + 2},\cdots,{\boldsymbol{\pi}}^*_M\right]}$ can be done in \textit{pre-processing}. For a given (fixed) value of $\bar m$, the total time for determining the whole coefficients $\boldsymbol{\beta}_{(n,m)}$, $(n,m)\in{\Bbb{S}_M}-{\Bbb{S}_{\bar m}}$ is linearly proportional to the number of elements in ${\Bbb{S}}_M - {\Bbb{S}}_{\bar m}$, $\left|{\Bbb{S}}_M\right| - \left|{\Bbb{S}}_{\bar m}\right|=(M+1)(M+2)/2 - ({\bar m}+1)({\bar m}+2)/2$, i.e., the time complexity of constructing ${\bf{w}}$ is $O(M^2)$.
\end{remark}}
\subsubsection{Polynomial basis for enhanced strains}
\textcolor{black}{Here, we show the construction of polynomial basis ${\bf{w}}_1$, ${\bf{w}}_2$, ${\bf{w}}_3$, and ${\bf{w}}_4$ in Eq.\,(\ref{wg_basis_matrix_gamma}), based on the method in Section\,\ref{app_ortho_cond}. The enhanced basis can be expressed by
\begin{align}
	{{\bf{w}}_i} = \left[ {{\boldsymbol{\pi }}_{{{\bar m}_i} + 1}^ * ,{\boldsymbol{\pi }}_{{{\bar m}_i} + 2}^ * , \cdots ,{\boldsymbol{\pi }}^ *_{{m_i}} } \right],\quad{i\in\left\{1,2,3,4\right\}},
\end{align}
where $m_i$ denotes the maximum degree of the polynomial, and ${\bar m}_i<{m_i}$ denotes the maximum degree in the corresponding compatible strain component, i.e., from Eq.\,(\ref{mat_A_poly_basis}), we have ${\bar m}_1={\bar m}_2=0$, ${\bar m}_3=2$, and ${\bar m}_4=1$. We first show that the orthogonality condition in Eq.\,(\ref{S_E_ortho_cond}) can be simplified to the form of Eq.\,(\ref{ortho_pbasis_mod_deter_b}). Since the orthogonality condition in Eq.\,(\ref{S_E_ortho_cond}) should hold for all $\boldsymbol{\alpha}\in\Bbb{R}^{d_\mathrm{a}}$, we obtain
\begin{align}
	\label{rew_ortho_b_alp_arb}
	{\int_\mathcal{A} {{\boldsymbol{\Gamma}}(\zeta^1,\zeta^2)^\mathrm{T}\underline{{\boldsymbol{S}}_{\rm{p}}}\,{j_0}\,{\rm{d}}\mathcal{A}} = \boldsymbol{0}}, 
\end{align}
where $\underline{{\boldsymbol{S}}_\mathrm{p}} \coloneqq {\left[ {{S_\mathrm{p}^{11}},{S_\mathrm{p}^{22}},{S_\mathrm{p}^{33}},{S_\mathrm{p}^{12}},{S_\mathrm{p}^{13}},{S_\mathrm{p}^{23}}} \right]^{\rm{T}}}$, such that ${{\boldsymbol{S}}_{\rm{p}}} = S_{\rm{p}}^{ij}{{\boldsymbol{G}}_i} \otimes {{\boldsymbol{G}}_j}$. Further, we make the following assumptions in the evaluation of the orthogonality condition in Eq.\,(\ref{rew_ortho_b_alp_arb}):
\begin{itemize}
	\item ${S^{ij}_\mathrm{p}}$ has the same polynomial order in $\zeta^1$ and $\zeta^2$ as the conjugate strain component ${E^{ij}_\mathrm{p}}$ in Eq.\,(\ref{phy_kin_dec}).
	\item $j_0\approx{1}$, which holds exactly only for initially straight beams or linear finite elements, see Remark\,\ref{rem_init_curv_dir}. This makes the orthogonality condition solely depend on the cross-section's initial shape.
\end{itemize}
Then, we obtain the following orthogonality condition: For each $(n,m)\in\Bbb{S}_{{m}_i}-\Bbb{S}_{{\bar m}_i}$, we have
\begin{align}
	\label{ortho_cond_w11_22_12_lin}
	\int_\mathcal{A}{P_{(p,q)}\,P^*_{(n,m)}\,{\rm{d}}\mathcal{A}}  = \boldsymbol{0},\quad \forall(p,q)\in{\Bbb{S}}_{{\bar m}_i},\,\,i\in\left\{1,2,3,4\right\},
\end{align}
which has the same form as Eq.\,(\ref{ortho_pbasis_mod_deter_b}). Subsequently, from Eq.\,(\ref{deter_beta_nm_lin}), we obtain
\begin{align}
	\label{deter_beta_nm_lin_case_i}
	{{\bar{\Bbb{W}}}_i}\,\boldsymbol{\beta}_{(n,m)} =  - \int_\mathcal{A} {{{\bf{\bar w}}_i^{\rm{T}}}{P_{(n,m)}}\,{\rm{d}}\mathcal{A}}
\end{align}
to determine the unknown coefficients ${\boldsymbol{\beta }}_{(n,m)}\coloneqq\left[{\boldsymbol{\beta }}_{\left( {n,m} \right)}^{0\,\mathrm{T}},\cdots,{\boldsymbol{\beta }}_{\left( {n,m} \right)}^{{\bar m}_i\,\mathrm{T}}\right]^{\!\mathrm{T}}\in\Bbb{R}^{{\left|{\Bbb{S}}_{{\bar m}_i}\right|}}$ in
\begin{align}
	\label{mod_P_poly_N_n_mod_form_case_i}
	P_{(n,m)}^* = {P_{(n,m)}} + {\bf{\bar w}}_i\,{{\boldsymbol{\beta }}_{(n,m)}},\,\,(n,m)\in{\Bbb{S}}_{m_i}-{\Bbb{S}}_{{\bar m}_i},
\end{align}
in each ${\bf{w}}_i\in\Bbb{R}^{1\times{d_i}}$, $i\in\left\{1,2,3,4\right\}$, where ${\bf{\bar w}}_i\coloneqq{\left[ {{{\boldsymbol{\pi }}_0}, \cdots ,{{\boldsymbol{\pi }}_{{\bar m}_i}}} \right]}\in\Bbb{R}^{1\times{\left|{\Bbb{S}}_{{\bar m}_i}\right|}}$, and the system matrix is given by
\begin{align}
	\label{deter_beta_nm_lin_case_i_system_mat}	
	{\bar{\Bbb{W}}}_i \coloneqq\int_\mathcal{A} {{\bf{\bar w}}_i^\mathrm{T}{\bf{\bar w}}_i\,{\rm{d}}\mathcal{A}}\in\Bbb{R}^{{\left|{\Bbb{S}}_{{\bar m}_i}\right|}\times{\left|{\Bbb{S}}_{{\bar m}_i}\right|}}.
\end{align}
Here, $d_i$ denotes the number of elements in the row array ${\bf{w}}_i$, and it is obtained by
\begin{align}
	\label{dimension_w_i}
	{d_i}\coloneqq{\left|{\Bbb{S}}_{m_i}\right|-\left|{\Bbb{S}}_{{\bar m}_i}\right|}=\dfrac{1}{2}(m_i+1)(m_i+2)-\dfrac{1}{2}({\bar m}_i+1)({\bar m}_i+2).
\end{align}	
Therefore, from Remark\,\ref{rem_determin_beta_lin}, the time complexity to construct ${\bf{w}}_i$ is $O({m_i}^2)$, $i\in\left\{1,2,3,4\right\}$.}
\section{Time-stepping scheme}
\subsection{Verification of the energy--momentum conservation for continuous time}
\label{sub_sec_verif_tcont_emc}
Here we verify that energy--momentum conservation laws can be deduced from Eqs.\,(\ref{el_eq_var})-(\ref{el_eq_var_a}) under appropriate boundary conditions. 
\subsubsection{Total linear momentum}
Let $\delta {\boldsymbol{y}} = {\left[ {{{\boldsymbol{c}}_\varphi^{\rm{T}}},{{\boldsymbol{0}}^{\rm{T}}},{{\boldsymbol{0}}^{\rm{T}}}} \right]^{\rm{T}}}\in{\bar{\mathcal{V}}}$ with a constant vector ${{\boldsymbol{c}}_\varphi}\in\mathbb{R}^3$, which describes a superposed rigid body translation. Here we assume no displacement boundary conditions, i.e., $\Gamma_\mathrm{D}=\emptyset$, so that ${\boldsymbol{c}}_\varphi\ne\boldsymbol{0}$. Then, the internal virtual work vanishes, and from Eq.\,(\ref{el_eq_var}), we obtain
\begin{equation}
	\label{cont_tot_lin_mnt_balance}
	{\boldsymbol{c}}_\varphi \cdot \left\{ {\frac{{\rm{d}}}{{{\rm{d}}t}}{\boldsymbol{L}}({\boldsymbol{V}}) - {{\boldsymbol{f}}_{{\rm{ext}}}}(t)} \right\} = 0,\,\,\forall {\boldsymbol{c}}_\varphi \in {{\mathbb{R}}^3},
\end{equation}
with the time-dependent total external force,
\begin{equation}
	{{\boldsymbol{f}}_{{\rm{ext}}}}(t) \coloneqq \int_0^L {{\renewcommand{\arraystretch}{1.25}\begin{array}{*{20}{c}}
				{{\boldsymbol{\bar n}}}
		\end{array}}{\rm{d}}s + {\renewcommand{\arraystretch}{1.25}\begin{array}{*{20}{c}}
				{{{{\boldsymbol{\bar n}}}_0}}
		\end{array}}}.
\end{equation}
Therefore, without any displacement boundary conditions, we have the \textit{balance of total linear momentum}
\begin{equation}
	\label{cont_tot_lin_mnt_balance_fext}	
	\frac{{\rm{d}}}{{{\rm{d}}t}}{\boldsymbol{L}}({\boldsymbol{V}}) = {{\boldsymbol{f}}_{{\rm{ext}}}}(t).
\end{equation}
This implies that the total linear momentum is conserved if no external force and displacement boundary condition are applied. 
\subsubsection{Total angular momentum}
Let $\delta {\boldsymbol{y}}\!\coloneqq\!{\left[ {{{(\boldsymbol{c}_\theta\times{\boldsymbol{\varphi}})}^{\rm{T}}},{{(\boldsymbol{c}_\theta\times{\boldsymbol{d}}_1)}^{\rm{T}}},{{(\boldsymbol{c}_\theta\times{\boldsymbol{d}}_2)}^{\rm{T}}}}\right]^\mathrm{T}} \in \bar {\mathcal{V}}$ with a constant vector ${\boldsymbol{c}}_\theta\in\mathbb{R}^3$, which represents a superposed (infinitesimal) rigid body rotation, and ${\boldsymbol{c}}_\theta$ represents the rotation axis vector. Here, we also assume no displacement boundary conditions, i.e., $\Gamma_\mathrm{D}=\emptyset$, so that ${\boldsymbol{c}}_\theta\ne\boldsymbol{0}$. Then, the internal virtual work vanishes, and from Eq.\,(\ref{el_eq_var}), we obtain
\begin{align}
	{{\boldsymbol{c}}_\theta } \cdot \left[ {\frac{{\rm{d}}}{{{\rm{d}}t}}{\boldsymbol{J}}({\boldsymbol{y}},{\boldsymbol{V}}) - {{\boldsymbol{m}}_{{\rm{ext}}}}(\boldsymbol{y},t)} \right] = 0,\,\,\forall {{\boldsymbol{c}}_\theta } \in {\mathbb{R}}^3,
\end{align}
with the total external moment defined by
\begin{align}
	{{\boldsymbol{m}}_{{\rm{ext}}}}(\boldsymbol{y},t) &\coloneqq \int_0^L \left( {{\boldsymbol{\varphi }} \times {\boldsymbol{\bar n}} + {{\boldsymbol{d}}_\alpha } \times {{{\boldsymbol{\bar {\tilde m}}}}\:\!\!^\alpha }} \right){\rm{d}}s + {{\left[ {{\boldsymbol{\varphi }} \times {{{\boldsymbol{\bar n}}}_0} + {{\boldsymbol{d}}_\alpha } \times {\boldsymbol{\bar {\tilde m}}}\:\!_0^\alpha } \right]}_{s \in {\Gamma _\mathrm{N}}}}.
\end{align}
Therefore, without any displacement boundary conditions, we have the \textit{balance of total angular momentum}
\begin{equation}
	\label{cont_proof_ang_mnt_bal}
	\frac{{\rm{d}}}{{{\rm{d}}t}}{\boldsymbol{J}}({\boldsymbol{y}},{\boldsymbol{V}}) = {{\boldsymbol{m}}_{{\rm{ext}}}}(\boldsymbol{y},t).
\end{equation}
\subsubsection{Total energy}
\label{tot_e_conserv_cont_f}
We define the total energy at an equilibrium solution of Eqs.\,(\ref{el_eq_var})-(\ref{el_eq_var_a}), by
\begin{equation}
	\label{tot_mech_e}
	E(\boldsymbol{y},\boldsymbol{V},{\boldsymbol{\varepsilon}}_\mathrm{p},\boldsymbol{\alpha},t)\coloneqq \mathcal{K}({\boldsymbol{V}}) + U({{\boldsymbol{\varepsilon }}_{\rm{p}}},\boldsymbol{\alpha}) - {W_{{\rm{ext}}}}({\boldsymbol{y}},t),
\end{equation}
where the internal energy can be evaluated, from Eq.\,(\ref{hw_int_ener}), by
\begin{equation}
	U({\boldsymbol{\varepsilon}_\mathrm{p}},{\boldsymbol{\alpha}})\coloneqq \int_0^L{\psi({\boldsymbol{\varepsilon}}_\mathrm{p},\boldsymbol{\alpha})}\,\mathrm{d}s.
\end{equation}
We choose $\delta {\boldsymbol{y}} = {{\boldsymbol{\dot y}}} \in {\bar {\mathcal{V}}}$, and have the compatibility condition $\boldsymbol{V}={\boldsymbol{\dot y}}$ from Eq.\,(\ref{compat_integ_v_y}). Then, from Eq.\,(\ref{el_eq_var}), we obtain
\begin{equation}
	\label{tderiv_tot_e_wfunc_m}
	\frac{{\rm{d}}}{{{\rm{d}}t}}E(\boldsymbol{y},\boldsymbol{V},{\boldsymbol{\varepsilon}}_\mathrm{p},\boldsymbol{\alpha},t)=-{\frac{\partial}{\partial t}{W_\mathrm{ext}}(\boldsymbol{y},t)}.
\end{equation}
This means that the total energy is conserved unless the external load has an explicit dependence on time. That is, for a conservative system, we have
\begin{equation}						\frac{{\rm{d}}}{{{\rm{d}}t}}E(\boldsymbol{y},\boldsymbol{V},{\boldsymbol{\varepsilon}}_\mathrm{p},\boldsymbol{\alpha})=0,
\end{equation}
where the total energy is
\begin{equation}
	E(\boldsymbol{y},\boldsymbol{V},{\boldsymbol{\varepsilon}}_\mathrm{p},\boldsymbol{\alpha})={\mathcal{K}({\boldsymbol{V}}) + U({{\boldsymbol{\varepsilon }}_{\rm{p}}},\boldsymbol{\alpha}) + {U_{{\rm{ext}}}}({\boldsymbol{y}})}.
\end{equation}
\subsection{Verification of the energy--momentum conservation in discrete-time intervals}
\label{app_verif_emc_discrete_time}
We show that the time-stepping scheme in Eqs.\,(\ref{eq_var_eq_y})-(\ref{t_emc_enh_a_cond0}) preserves the conservation properties of the continuous form in Section \ref{sub_sec_verif_tcont_emc}.
\subsubsection{Total linear momentum}
Let $\delta {\boldsymbol{y}} = {\left[ {{\boldsymbol{c}}_\varphi ^{\rm{T}},{{\boldsymbol{0}}^{\rm{T}}},{{\boldsymbol{0}}^{\rm{T}}}} \right]^{\rm{T}}}\in{\bar{\mathcal{V}}}$ with a constant vector ${\boldsymbol{c}}_\varphi\in\mathbb{R}^3$. Here we assume no displacement boundary condition, i.e., $\Gamma_\mathrm{D}=\emptyset$, so that ${\boldsymbol{c}}_\varphi\ne\boldsymbol{0}$. Then, we obtain the discrete form of the linear momentum balance from Eq.\,(\ref{eq_var_eq_y}), as
\begin{align}
	\frac{1}{{\Delta t}}\!\left\{ {{\boldsymbol{L}}\left( {{}^{n + 1}{\boldsymbol{V}}} \right)\!-\! {\boldsymbol{L}}\left( {{}^n{\boldsymbol{V}}} \right)} \right\}\!=\! {{\boldsymbol{f}}_{{\rm{ext}}}}\!\left( {{t_n} + \frac{1}{2}\Delta t} \right),
\end{align}
which is consistent with Eq.\,(\ref{cont_tot_lin_mnt_balance}). Therefore, the total linear momentum is conserved, if no external force and displacement boundary conditions are applied. 
\subsubsection{Total angular momentum}
We choose
\begin{equation}
	\delta {\boldsymbol{y}} = \left\{ {\begin{array}{*{20}{c}}
			{{{\boldsymbol{c}}_\theta } \times {}^{n + \frac{1}{2}}{\boldsymbol{\varphi }}}\\
			{{{\boldsymbol{c}}_\theta } \times {}^{n + \frac{1}{2}}{{\boldsymbol{d}}_1}}\\
			{{{\boldsymbol{c}}_\theta } \times {}^{n + \frac{1}{2}}{{\boldsymbol{d}}_2}}
	\end{array}} \right\}\in{\bar{\mathcal{V}}},
\end{equation}
with a constant vector ${\boldsymbol{c}}_\theta\in\mathbb{R}^3$, and assume no displacement boundary conditions, i.e., $\Gamma_\mathrm{D}=\emptyset$, so that ${\boldsymbol{c}_\theta}\ne\boldsymbol{0}$. Then, from Eq.\,(\ref{eq_var_eq_y}), we obtain 
\begin{align}
	\frac{1}{{\Delta t}}\left\{ {{\boldsymbol{J}}\left( {{}^{n + 1}{\boldsymbol{y}},{}^{n + 1}{\boldsymbol{V}}} \right) - {\boldsymbol{J}}\left( {{}^n{\boldsymbol{y}},{}^n{\boldsymbol{V}}} \right)} \right\} = {{\boldsymbol{m}}_{{\rm{ext}}}}\left( {{}^{n + \frac{1}{2}}{\boldsymbol{y}},{t_n} + \frac{1}{2}\Delta t} \right),
\end{align}
which is consistent with Eq.\,(\ref{cont_proof_ang_mnt_bal}). Therefore, the angular momentum is conserved if no external loads and boundary conditions are applied. Note that this holds for both the conventional mid-point rule and the present EMC scheme.
\subsubsection{Total energy}
\label{tdisc_conserv_tot_e}
We choose $\delta {\boldsymbol{y}}={}^{n+1}{\boldsymbol{y}}-{}^{n}{\boldsymbol{y}}\in{\bar{\mathcal{V}}}$. Then, in Eq.\,(\ref{eq_var_eq_y}), using Eq.\,(\ref{rep_cur_V_conf}), we have
\begin{align}
	\label{preserv_tot_e_kin}
	&\frac{2}{{\Delta {t}}}{G_{{\rm{iner}}}}\left( \frac{1}{\Delta t}\left({{}^{n + 1}{\boldsymbol{y}}-{}^{n}{\boldsymbol{y}}}\right)-{}^n\boldsymbol{V},{{}^{n + 1}{\boldsymbol{y}}-{}^{n}{\boldsymbol{y}}}\right)\nonumber\\
	&= \frac{1}{2}\int_0^L {\left( {{}^{n + 1}{\boldsymbol{V}} + {}^n{\boldsymbol{V}}} \right) \cdot {\boldsymbol{\mathcal{M}}}\left( {{}^{n + 1}{\boldsymbol{V}} - {}^n{\boldsymbol{V}}} \right)\,\,{\rm{d}}s} \nonumber\\
	&= {\mathcal{K}}({}^{n + 1}{\boldsymbol{V}}) - {\mathcal{K}}({}^n{\boldsymbol{V}}).
\end{align}
Further, due to the linearity of the operator $\mathbb{B}(\boldsymbol{y})$ with respect to $\boldsymbol{y}$, we can obtain
\begin{align}
	\label{btot_phy_strn_e_dif}
	{{\mathbb{B}}}\!\left( {{}^{n + \frac{1}{2}}{\boldsymbol{y}}} \right)\delta\boldsymbol{y} &= {\boldsymbol{\varepsilon }}\!\left( {{}^{n + 1}{\boldsymbol{y}}} \right) - {\boldsymbol{\varepsilon }}\!\left( {{}^n{\boldsymbol{y}}} \right).
\end{align}
Then, from Eq.\,(\ref{def_int_vw_y}), we have
\begin{align}
	\label{integ_identity_Bop_int_u}	
	G_{{\mathop{\rm int}} }^{\rm{y}}\left( {{}^{n + \frac{1}{2}}{\boldsymbol{y}},{}^{n + \frac{1}{2}}{{\boldsymbol{r}}_{\rm{p}}},{}^{n + 1}{\boldsymbol{y}} - {}^n{\boldsymbol{y}}} \right) 
	= \int_0^L {\left\{ {{\boldsymbol{\varepsilon }}\left( {{}^{n + 1}{\boldsymbol{y}}} \right) - {\boldsymbol{\varepsilon }}\left( {{}^n{\boldsymbol{y}}} \right)} \right\} \cdot {}^{n + \frac{1}{2}}{{\boldsymbol{r}}_{\rm{p}}}\,{\rm{d}}s}.
\end{align}
By choosing $\delta {\boldsymbol{\varepsilon}}_\mathrm{p}={}^{n+1}{\boldsymbol{\varepsilon}}_\mathrm{p}-{}^{n}{\boldsymbol{\varepsilon}}_\mathrm{p}\in{\mathcal{V}}_\mathrm{p}$ and $\delta {\boldsymbol{\alpha}}={}^{n+1}{\boldsymbol{\alpha}}-{}^{n}{\boldsymbol{\alpha}}\in{\mathcal{V}}_\mathrm{a}$ in Eqs.\,(\ref{t_emc_constitutiv_cond0}) and (\ref{t_emc_enh_a_cond0}), respectively, and combining the resulting equations, we obtain the identity,
\begin{align}
	\label{comb_const_stres_cond_e_a}
	\int_0^L {\left\{ {\psi \left( {{}^{n + 1}{{\boldsymbol{\varepsilon }}_{\rm{p}}},{}^{n + 1}{\boldsymbol{\alpha }}} \right) - \psi \left( {{}^n{{\boldsymbol{\varepsilon }}_{\rm{p}}},{}^n{\boldsymbol{\alpha }}} \right)} \right\}\,{\rm{d}}s} = \int_0^L {\left( {{}^{n + 1}{{\boldsymbol{\varepsilon }}_{\rm{p}}} - {}^n{{\boldsymbol{\varepsilon }}_{\rm{p}}}} \right) \cdot {}^{n + \frac{1}{2}}{{\boldsymbol{r}}_{\rm{p}}}\,{\rm{d}}s}.
\end{align}
Note that, here, we assume linear constitutive laws. Evaluating Eq.\,(\ref{t_emc_compat_cond}) at time $t=t_{n+1}$ and $t=t_{n}$ with $\delta{\boldsymbol{r}}_\mathrm{p}={}^{n+\frac{1}{2}}{\boldsymbol{r}}_\mathrm{p}\in{\mathcal{V}}_\mathrm{p}$, we have  
\begin{subequations}	
	\begin{align}
		\int_0^L {{}^{n + \frac{1}{2}}{{\boldsymbol{r}}_{\rm{p}}} \cdot \left\{ {{\boldsymbol{\varepsilon }}\left( {{}^{n + 1}{\boldsymbol{y}}} \right) - {}^{n + 1}{{\boldsymbol{\varepsilon }}_{\rm{p}}}} \right\}{\rm{d}}s}  &= 0,\label{compat_eval_at_time_tn}\\	
		\int_0^L {{}^{n + \frac{1}{2}}{{\boldsymbol{r}}_{\rm{p}}} \cdot \left\{ {{\boldsymbol{\varepsilon }}\left( {{}^n{\boldsymbol{y}}} \right) - {}^n{{\boldsymbol{\varepsilon }}_{\rm{p}}}} \right\}{\rm{d}}s}  &= 0,\label{compat_eval_at_time_tn_1}
	\end{align}
\end{subequations}
respectively. Equations\,(\ref{compat_eval_at_time_tn}) and (\ref{compat_eval_at_time_tn_1}) imply that Eqs.\,(\ref{integ_identity_Bop_int_u}) and (\ref{comb_const_stres_cond_e_a}) are the same, so that we obtain the identity
\begin{align}
	\label{te_pres_eq_Gy_int}
	G_{{\mathop{\rm int}} }^{\rm{y}}\left( {{}^{n + \frac{1}{2}}{\boldsymbol{y}},{}^{n + \frac{1}{2}}{{\boldsymbol{r}}_{\rm{p}}},{}^{n + 1}{\boldsymbol{y}} - {}^n{\boldsymbol{y}}} \right)
	&=\int_0^L {\left\{ {\psi \left( {{}^{n + 1}{{\boldsymbol{\varepsilon }}_{\rm{p}}},{}^{n + 1}{\boldsymbol{\alpha }}} \right) - \psi \left( {{}^n{{\boldsymbol{\varepsilon }}_{\rm{p}}},{}^n{\boldsymbol{\alpha }}} \right)} \right\}\,{\rm{d}}s}\nonumber\\
	&=U({{}^{n + 1}{{\boldsymbol{\varepsilon }}_{\rm{p}}},{}^{n + 1}{\boldsymbol{\alpha }}}) - U({{}^{n}{{\boldsymbol{\varepsilon }}_{\rm{p}}},{}^{n}{\boldsymbol{\alpha }}}).
\end{align}
Substituting Eqs.\,(\ref{preserv_tot_e_kin}) and (\ref{te_pres_eq_Gy_int}) into Eq.\,(\ref{eq_var_eq_y}), we have
\begin{align}
	&\mathcal{K}({}^{n + 1}{\boldsymbol{V}}) + {U}({}^{n + 1}{\boldsymbol{\varepsilon}_{\rm{p}}},{}^{n + 1}{\boldsymbol{\alpha }}) - {W_{{\rm{ext}}}}\left( {{}^{n + 1}{\boldsymbol{y}}} ,{t_n+\frac{1}{2}\Delta t}\right)\nonumber\\
	&= \mathcal{K}({}^n{\boldsymbol{V}}) + {U}({}^n{\boldsymbol{\varepsilon}_{\rm{p}}},{}^n{\boldsymbol{\alpha }}) - {W_{{\rm{ext}}}}\left( {{}^n{\boldsymbol{y}}},{t_n+\frac{1}{2}\Delta t}\right).
\end{align}
This implies the following:
\begin{itemize}
	\item If the applied load is (piece-wise) linearly proportional to time, we obtain the following identity 
	\begin{align}
		\label{total_energy_consistency_pie_lin_f}
		{}^{n+1}E={}^{n+1}E^*,\quad n=0,1,\cdots,
	\end{align}
	where we define
	\begin{align*}
		{}^{n + 1}E&\coloneqq E\left( {{}^{n + 1}{\boldsymbol{y}},{}^{n + 1}{\boldsymbol{V}},{}^{n + 1}{{\boldsymbol{\varepsilon }}_{\rm{p}}},{}^{n + 1}{\boldsymbol{\alpha }},{t_{n + 1}}} \right),
	\end{align*}
	\begin{align*}
		{}^{n+1}E^*\coloneqq{}^{n}E-\overline {{\Delta _t}{W_{{\rm{ext}}}}}\,,
	\end{align*}
	with
	\begin{align*}
		{}^{n}E&\coloneqq E\left( {{}^{n}{\boldsymbol{y}},{}^{n}{\boldsymbol{V}},{}^{n}{{\boldsymbol{\varepsilon }}_{\rm{p}}},{}^{n}{\boldsymbol{\alpha }},{t_{n}}} \right),		
	\end{align*}	
	and
	\begin{align*}
		\overline {{\Delta _t}{W_{{\rm{ext}}}}}\coloneqq\frac{1}{2}\left( {{{\left. {\frac{{\partial {W_{{\rm{ext}}}}}}{{\partial t}}} \right|}_{{\boldsymbol{y}} = {}^{n + 1}{\boldsymbol{y}},\,t = {t_{n + 1}}}} + {{\left. {\frac{{\partial {W_{{\rm{ext}}}}}}{{\partial t}}} \right|}_{{\boldsymbol{y}} = {}^n{\boldsymbol{y}},\,t = {t_n}}}} \right)\Delta t.
	\end{align*}		
	\item If the external load has no explicit time-dependence, the total energy is conserved, i.e., ${}^{n + 1}E = {}^{n}E$.
\end{itemize}
\section{Numerical examples}
\subsection{An analytical solution for a pure torsion problem}
\label{app_asol_tors_stiff}
We consider a pure torsion of a straight bar of length $L$, and a shear modulus $G=E/(2(1+\nu))$. In the isotropic linear elasticity, the torsional angle (radians) of the bar under a twisting moment ${M_\mathrm{T}}$ can be expressed by \citep[page 350]{budynas2020roark}
\begin{equation}
	\label{sol_sv_tor}
	\theta = \frac{{{M_\mathrm{T}}L}}{{KG}},
\end{equation}
where $K$ is a factor which depends on the cross-sectional shape and dimensions. For a circular cross-section, the cross-section remains plane under the moment load, and $K$ is equal to the polar moment of inertia, $I_\mathrm{p}$. For other cross-sections, the cross-section warps, and $K$ is less than $I_\mathrm{p}$. For a (solid) rectangular cross-section of dimensions $2a$ and $2b$ ($a \ge b$), we have a solution for $K$ in \citet[page 366, Table 10.7]{budynas2020roark}, 
\begin{equation}
	\label{fac_K_sv}
	K = a{b^3}\left\{ {\frac{{16}}{3} - 3.36\frac{b}{a}\left( {1 - \frac{{{b^4}}}{{12{a^4}}}} \right)} \right\},
\end{equation}
obtained from a truncation of an infinite series solution of the St.\,Venant torsion problem, which can be found in \citet[page 288]{timoshenko1951theory}.
From Eq.\,(\ref{sol_sv_tor}), $KG/L={M_\mathrm{T}}/\theta$ represents the torsional stiffness of the bar. Note that Eq.\,(\ref{fac_K_sv}) is based on the assumptions \citep[page 349]{budynas2020roark}: (i) The bar is straight, (ii) the material is homogeneous, isotropic, and linear elastic, (iii) the deformation is purely torsional under equal and opposite twisting moments at the ends.
\subsection{Total number of load steps in example 4}
\textcolor{black}{
Here we present total number of load steps and iterations in the results of Fig.\,\ref{twist_ring_deformed_rad_conv}. For the formulation [$\Delta \theta$], we use in total 200 and 2,000 load steps with uniform increments in cases 1 and 2, respectively. Table\,\ref{twist_nlstep_mnt_conv_test} shows the total number of load steps and iterations in the results from using [mix.glo] and [mix.loc-sr]. We have observed that, using [mix.loc-sr] with $p=2$ and $n_\mathrm{el}=24$, the ring does not deform into a smaller one, and the equilibrium iteration eventually diverges. 
}
%
\begin{table}[htbp]
	\centering
	\caption{Twisting of an elastic ring: Total number of load steps ($n_\mathrm{load}$) with uniform increments, and total number of iterations ($n_\mathrm{iter}$), for $\bar \theta=2\pi$ in the results of Fig.\,\ref{twist_ring_deformed_rad_conv}. Here, the symbol $\times$ represents failure in convergence of the equilibrium iteration.}
    \begin{tabular}{ccccrccrccrcc}
		\toprule
		&       & \multicolumn{5}{c}{Case 1: $w=1/3$, $h=1$} &       & \multicolumn{5}{c}{Case 2: $w=1/30$, $h=1/10$} \\
		\cmidrule{3-7}\cmidrule{9-13}          &       & \multicolumn{2}{c}{[mix.loc-sr]} &       & \multicolumn{2}{c}{[mix.glo]} &       & \multicolumn{2}{c}{[mix.loc-sr]} &       & \multicolumn{2}{c}{[mix.glo]} \\
		\cmidrule{3-4}\cmidrule{6-7}\cmidrule{9-10}\cmidrule{12-13}    \multicolumn{1}{l}{$p$} & $n_\mathrm{el}$ & $n_\mathrm{load}$ & $n_\mathrm{iter}$ &       & $n_\mathrm{load}$ & $n_\mathrm{iter}$ &       & $n_\mathrm{load}$ & $n_\mathrm{iter}$ &       & $n_\mathrm{load}$ & $n_\mathrm{iter}$ \\
		\midrule
		\multirow{5}[2]{*}{2} & 24    & 96    & $\times$     &       & 16    & 105   &       & 96    & $\times$     &       & 16    & 103 \\
		& 48    & 96    & 641   &       & 16    & 103   &       & 96    & 523   &       & 16    & 101 \\
		& 96    & 96    & 673   &       & 16    & 103   &       & 192   & 973   &       & 16    & 100 \\
		& 192   & 16    & 224   &       & 16    & 102   &       & 16    & 183   &       & 16    & 100 \\
		& 384   & 384   & $\times$ &       & 20    & 161   &       & 16    & 208   &       & 16    & 100 \\
		\midrule
		\multirow{5}[2]{*}{3} & 24    & 16    & 105   &       & 16    & 105   &       & 16    & 105   &       & 16    & 103 \\
		& 48    & 16    & 104   &       & 16    & 103   &       & 16    & 103   &       & 16    & 101 \\
		& 96    & 16    & 104   &       & 16    & 103   &       & 16    & 103   &       & 16    & 100 \\
		& 192   & 16    & 116   &       & 16    & 114   &       & 16    & 103   &       & 16    & 100 \\
		& 384   & 16    & 144   &       & 20    & 175   &       & 16    & 103   &       & 16    & 100 \\
		\midrule
		\multirow{5}[2]{*}{4} & 24    & 16    & 150   &       & 16    & 105   &       & 16    & 109   &       & 16    & 103 \\
		& 48    & 16    & 140   &       & 16    & 103   &       & 16    & 104   &       & 16    & 101 \\
		& 96    & 16    & 141   &       & 16    & 103   &       & 16    & 103   &       & 16    & 100 \\
		& 192   & 64    & 480   &       & 16    & 157   &       & 16    & 104   &       & 16    & 100 \\
		& 384   & 64    & 445   &       & 36    & 253   &       & 16    & 104   &       & 16    & 100 \\
		\bottomrule
	\end{tabular}%
	\label{twist_nlstep_mnt_conv_test}
\end{table}%
\end{appendices}


\bibliography{sn-bibliography}

\end{document}